\newcommand{\etal}{\textit{et al.}}
\documentclass[final,12pt,letterpaper]{elsarticle}

\usepackage{amssymb}
\usepackage{amsthm}
\usepackage{amsmath,mathrsfs}

\biboptions{sort&compress}
\newdefinition{rmk}{Remark}

\journal{Journal of Computational Physics}

\newcommand{\p}[1]{\ensuremath{\left(#1\right)}}
\newcommand{\fst}{1$^{\textrm{st}}$}
\newcommand{\snd}{2$^{\textrm{nd}}$}
\newcommand{\trd}{3$^{\textrm{rd}}$}
\newcommand{\fth}{4$^{\textrm{th}}$}
\newcommand{\ie}{\emph{i.e.}}
\newcommand{\eg}{\emph{e.g.}}
\newcommand{\myremarkend}{~\hfill$\clubsuit\/$}

\begin{document}

\begin{frontmatter}



\title{A Correction Function Method for Poisson Problems with Interface Jump Conditions.}


\author[mit-aa]{Alexandre Noll Marques}
\author[mcgill]{Jean-Christophe Nave}
\author[mit-math]{Rodolfo Ruben Rosales}

\address[mit-aa]{Department of Aeronautics and Astronautics, Massachusetts Institute of Technology\\Cambridge, MA 02139-4307}
\address[mcgill]{Department of Mathematics and Statistics, McGill University\\Montreal, Quebec H3A 2K6, Canada}
\address[mit-math]{Department of Mathematics, Massachusetts Institute of Technology\\Cambridge, MA 02139-4307}

\begin{abstract}
In this paper we present a method to treat interface jump conditions for constant
coefficients Poisson problems that allows the use of standard ``black box'' solvers,
without compromising accuracy. The basic idea of the new approach is similar to the
Ghost Fluid Method (GFM). The GFM relies on corrections applied on nodes located
across the interface for discretization stencils that straddle the interface.
If the corrections are solution-independent, they can be moved to the
right-hand-side (RHS) of the equations, producing a problem with the same linear
system as if there were no jumps, only with a different RHS. However, achieving
high accuracy is very hard (if not impossible) with the ``standard'' approaches
used to compute the GFM correction terms.

In this paper we generalize the GFM correction terms to a correction function,
defined on a band around the interface. This function is then shown to be
characterized as the solution to a PDE, with appropriate boundary conditions.
This PDE can, in principle, be solved to any desired order of accuracy. As an
example, we apply this new method to devise a \fth\/ order accurate scheme for
the constant coefficients Poisson equation with discontinuities in 2D. This scheme
is based on (i) the standard 9-point stencil discretization of the Poisson equation,
(ii) a representation of the correction function in terms of  bicubics, and
(iii) a solution of the correction function PDE by a least squares minimization.
Several applications of the method are presented to illustrate its robustness
dealing with a variety of interface geometries, its capability to capture
sharp discontinuities, and its high convergence rate.
\end{abstract}


\begin{keyword}
 Poisson equation \sep interface jump condition \sep Ghost Fluid method \sep Gradient augmented level set method \sep High accuracy \sep Hermite cubic spline

\PACS 47.11-j \sep 47.11.Bc

\MSC[2010] 76M20 \sep 35N06



\end{keyword}

\end{frontmatter}


\section{Introduction.} \label{sec:intro}
%
\subsection{Motivation and background information.} \label{sub:int:MB}
In this paper we present a new and efficient method to solve the constant coefficients
Poisson equation in the presence of discontinuities across an interface, with a high
order of accuracy. Solutions of the Poisson equation with discontinuities are
of fundamental importance in the description of fluid flows separated by
interfaces (\eg\/~the contact surfaces for immiscible multiphase fluids,
or fluids separated by a membrane) and other multiphase diffusion phenomena.
Over the last three decades, several methods have been developed to solve problems
of this type numerically~\cite{peskin:77, sussman:94, leveque:94, leveque:97,
   johansen:98, fedkiw:99, fedkiwetal:99, kang:00, liu:00, lai:00, li:01,
   nguyen:01, lee:03, gibou:07, gong:08, dolbow:09, bedrossian:10}.
However, obtaining a high order of accuracy still poses great challenges in
terms of complexity and computational efficiency.

When the solution is known to be smooth, it is easy to obtain highly accurate
finite-difference discretizations of the Poisson equation on a regular grid.
Furthermore, these discretizations commonly yield symmetric and banded linear
systems, which can be inverted efficiently \cite{trefthen:97}. On the
other hand, when singularities occur (\eg\/~discontinuities) across
internal interfaces, some of the regular discretization stencils will straddle
the interface, which renders the whole procedure invalid.

Several strategies have been proposed to tackle this issue.
Peskin~\cite{peskin:77} introduced the Immersed Boundary Method
(IBM)~\cite{peskin:77, lai:00}, in which the discontinuities are re-interpreted
as additional (singular) source terms concentrated on the interface. These
singular terms are then ``regularized'' and appropriately spread out over the
regular grid --- in a ``thin'' band enclosing the interface. The result is a
first order scheme that smears discontinuities. In order to avoid this
smearing of the interface information, LeVeque and Li~\cite{leveque:94}
developed the Immersed Interface Method
(IIM)~\cite{leveque:94, leveque:97, li:01, lee:03}, which is a methodology to
modify the discretization stencils, taking into consideration the
discontinuities at their actual locations. The IIM guarantees second order
accuracy and sharp discontinuities, but at the cost of added discretization
complexity and loss of symmetry.

The new method advanced in this paper builds on the ideas introduced by the
Ghost Fluid Method
(GFM)~\cite{mayo:84, fedkiw:99, fedkiwetal:99, liu:00, kang:00, nguyen:01, gibou:07}.
The GFM is based on defining both actual and ``ghost'' fluid variables at
every node on a narrow band enclosing the interface. The ghost variables work
as extensions of the actual variables across the interface --- the solution on
each side of the interface is assumed to have a smooth extension into the other
side. This approach allows the use of standard discretizations everywhere in the
domain. In most GFM versions, the ghost values are written as the actual
values, plus corrections that are independent of the underlying solution to the
Poisson problem. Hence, the corrections can be pre-computed, and moved into the
source term for the equation. In this fashion the GFM yields the same linear
system as the one produced by the problem without an interface, except for
changes in the right-hand-side (sources) only, which can then be inverted just
as efficiently.

The key difficulty in the GFM is the calculation of the correction terms, since
the overall accuracy of the scheme depends heavily on the quality of the
assigned ghost values. In~\cite{fedkiw:99, fedkiwetal:99, liu:00, kang:00, nguyen:01}
the authors develop first order accurate approaches to deal with discontinuities.
In the present work, we show that for the constant coefficients Poisson equation 
we can generalize the GFM correction term (at each ghost point) concept to that
of a correction function defined on a narrow band enclosing the interface.
Hence we call this new approach the \emph{Correction Function Method} (CFM). This
correction function is then shown to be characterized as the solution to a PDE,
with appropriate boundary conditions on the interface --- see \S~\ref{sec:general}.
Thus, at least in principle, one can calculate the correction function to any order
of accuracy, by designing algorithms to solve the PDE that defines it. In this
paper we present examples of \snd\/ and \fth\/ order accurate schemes (to solve the
constant coefficients Poisson equation, with discontinuities across interfaces,
in 2D) developed using this general framework. 

A key point (see \S~\ref{sec:scheme}) in the scheme developed here is the way we
solve the PDE defining the correction function. This PDE is solved in a weak fashion
using a least squares minimization procedure. This provides a flexible approach that
allows the development of a robust scheme that can deal with the geometrical
complications of the placement of the regular grid stencils relative to the
interface. Furthermore, this approach is easy to generalize to 3D, or to 
higher orders of accuracy.

\subsection{Other related work.} \label{sub:int:RW}
It is relevant to note other developments to solve the Poisson equation under
similar circumstances --- multiple phases separated by interfaces --- but with
different interface conditions. The Poisson problem with Dirichlet boundary
conditions, on an irregular boundary embedded in a regular grid, has been
solved to second order of accuracy using fast Poisson solvers~\cite{mayo:84},
finite volume~\cite{johansen:98}, and finite
differences~\cite{udaykumar:99, gibou:02, jomaa:05, gibou:05}
approaches. In particular, Gibou \etal\/~\cite{gibou:02} and Jomaa and
Macaskill~\cite{jomaa:05} have shown that it is possible to obtain symmetric
discretizations of the embedded Dirichlet problem, up to second order of
accuracy. Gibou and Fedkiw~\cite{gibou:05} have developed a fourth order
accurate discretization of the problem, at the cost of giving up symmetry.
More recently, the same problem has also been solved, to second order of
accuracy, in non-graded adaptive Cartesian grids by Chen \etal~\cite{chen:07}.
Furthermore, the embedded Dirichlet problem is closely
related to the Stefan problem modeling dendritic growth, as described
in~\cite{sethian:92, chen:97}.

The finite-element community has also made significant progress in incorporating
the IIM and similar techniques to solve the Poisson equation using
embedded grids. In particular the works by Gong \etal\/~\cite{gong:08}, Dolbow and
Harari~\cite{dolbow:09}, and Bedrossian \etal\/~\cite{bedrossian:10} describe
second order accurate finite-element discretizations that result in symmetric
linear systems. Moreover, in these works the interface (or boundary) conditions
are imposed in a weak fashion, which bears some conceptual similarities with the
CFM presented here, although the execution is rather different.

\subsection{Interface representation.} \label{sub:int:IR}
Another issue of primary importance to multiphase problems is the
representation of the interface (and its tracking in unsteady cases). Some
authors (see~\cite{peskin:77, mayo:84, udaykumar:99}) choose to
represent the interface explicitly, by tracking interface particles. The
location of the neighboring particles is then used to produce local
interpolations (\eg\/~splines), which are then applied to compute
geometric information --- such as curvature and normal directions. Although
this approach can be quite accurate, it requires special treatment when the
interface undergoes either large deformations or topological changes --- such
as mergers or splits. Even though we are not concerned with these issues in this
paper, we elected to adopt an implicit representation, to avoid complications
in future applications. In an implicit representation, the interface is given
as the zero level of a function that is defined everywhere in the regular grid
--- the level set function~\cite{osher:88}. In particular, we adopted the
Gradient-Augmented Level Set (GA-LS) method~\cite{nave:10}. With this
extension of the level set method, we can obtain highly accurate
representations of the interface, and other geometric information, with the
additional advantage that this method uses only local grid information.
We discuss the question of interface representation in more detail in
\S~\ref{sub:interface}.

\subsection{Organization of the paper.} \label{sub:int:OP}
The remainder of the paper is organized as follows. In \S~\ref{sec:problem}
we introduce the Poisson problem that we seek to solve. In
\S~\ref{sec:idea} the basic idea behind the solution method, and its
relationship to the GFM, are explored. Next, in \S~\ref{sec:general}, we
introduce the concept of the correction function and show how it is defined
by a PDE problem. In \S~\ref{sec:scheme} we apply this new framework to
build a \fth\/ order accurate scheme in 2D. Since the emphasis of this paper
is on high-order schemes, we describe the \snd\/ order accurate scheme in
appendix~\ref{ap:scheme2}. Next, in \S~\ref{sec:results} we demonstrate the
robustness and accuracy of the 2D scheme by applying it to several example
problems. The conclusions are in \S~\ref{sec:conclusion}. In appendix~\ref{ap:bicubic}
we review some background material, and notation, on bicubic interpolation. Finally,
in appendix~\ref{ap:omega} we discuss some technical issues regarding the
construction of the sets where the correction function is solved for.

\section{Definition of the problem.} \label{sec:problem}
\begin{figure}[htb!]
 \begin{center}
  \includegraphics[width=2.8in]{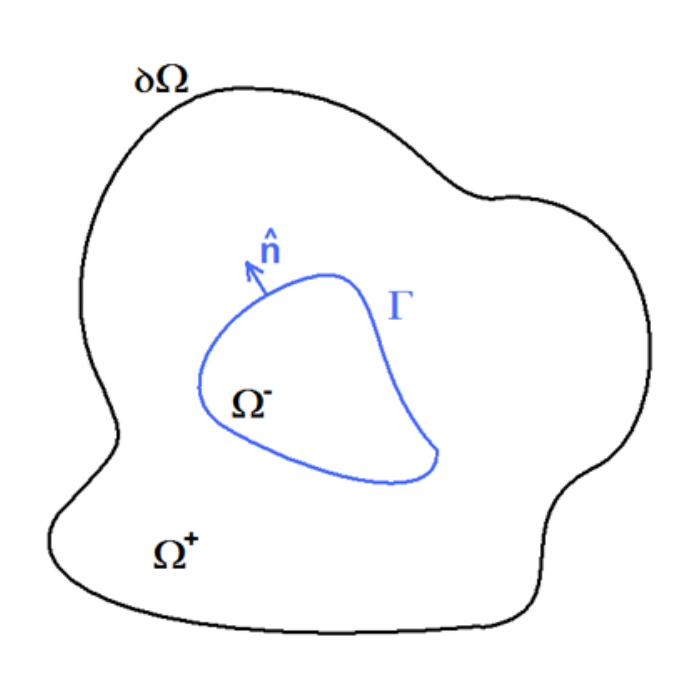}
 \end{center}
 \caption{Example of a domain $\Omega$, with an interface $\Gamma$.}
 \label{fig:problem}
\end{figure}

Our objective is to solve the constant coefficients Poisson's equation in a domain
$\Omega$ in which the
solution is discontinuous across a co-dimension 1 interface $\Gamma$, which
divides the domain into the subdomains $\Omega^+$ and $\Omega^-$, as illustrated
in figure~\ref{fig:problem}. We use the notation $u^+$ and $u^-$ to denote the
solution in each of the subdomains. Let the discontinuities across $\Gamma$ be
given in terms of two functions defined on the interface: $a = a(\vec{x})$ for
the jump in the function values, and $b = b(\vec{x})$ for the jump in the
normal derivatives. Furthermore, assume Dirichlet boundary conditions on the
``outer'' boundary $\partial \Omega$ (see figure~\ref{fig:problem}). Thus the
problem to be solved is
\begin{subequations}\label{eq:poisson}
 \begin{align}
   \nabla^2u\p{\vec{x}} &= f\p{\vec{x}} &\mathrm{for}\;\; \vec{x}
   &\in \Omega\/, \label{eq:poisson-eq}\\
   [u]_{\Gamma} &= a\p{\vec{x}} &\mathrm{for}\;\; \vec{x}
   &\in \Gamma\/, \label{eq:a}\\
   \left[u_n\right]_{\Gamma} &= b\p{\vec{x}} &\mathrm{for}\;\; \vec{x}
   &\in \Gamma\/,  \label{eq:b}\\
   u\p{\vec{x}} &= g\p{\vec{x}} &\mathrm{for}\;\; \vec{x}
   &\in \partial\Omega\/, \label{eq:dirichlet}
 \end{align}
\end{subequations}
where
\begin{subequations}\label{eq:jump}
 \begin{align}
 [u]_{\Gamma} &= u^+\p{\vec{x}} - u^-\p{\vec{x}}
    &\mathrm{for}\;\; \vec{x} &\in \Gamma\/,\\
 \left[u_n\right]_{\Gamma} &= u_n^+\p{\vec{x}} - u_n^-\p{\vec{x}}
    &\mathrm{for}\;\; \vec{x} &\in \Gamma\/.
 \end{align}
\end{subequations}
Throughout this paper, $\vec{x} = (x_1\/,\,x_2\/,\,\dots) \in \mathbb{R}^\nu\/$
is the spatial vector (where $\nu = 2\/$, or $\nu = 3\/$), and $\nabla^2\/$ is
the Laplacian operator defined by
\begin{equation}
 \nabla^2 = \sum_{i=1}^{\nu}\dfrac{\partial^2}{\partial x_i^2}\/.
\end{equation}
Furthermore,
\begin{equation}
 u_n = \hat{n}\cdot\vec{\nabla}u =
       \hat{n}\cdot\p{u_{x_1}\/,\,u_{x_2}\/,\,\dots}
\end{equation}
denotes the derivative of $u$ in the direction of $\hat{n}$, the unit vector
normal to the interface $\Gamma$ pointing towards $\Omega^+$
(see figure~\ref{fig:problem}).

It is important to note that our method focuses on the discretization of the
problem in the vicinity of the interface only. Thus, the method is compatible
with any set of boundary conditions on $\partial \Omega$, not just Dirichlet.

\section{Solution method -- the basic idea.} \label{sec:idea}
To achieve the goal of efficient high order discretizations of Poisson
equation in the presence of discontinuities, we build on the idea of the
Ghost Fluid Method (GFM). In essence, we use a standard discretization of the
Laplace operator (on a domain without an interface $\Gamma\/$) and modify the
right-hand-side (RHS) to incorporate the jump conditions across $\Gamma\/$.
Thus, the resulting linear system can be inverted as efficiently as in the
case of a solution without discontinuities.
\begin{figure}[htb!]
 \begin{center}
  \includegraphics[width=4in]{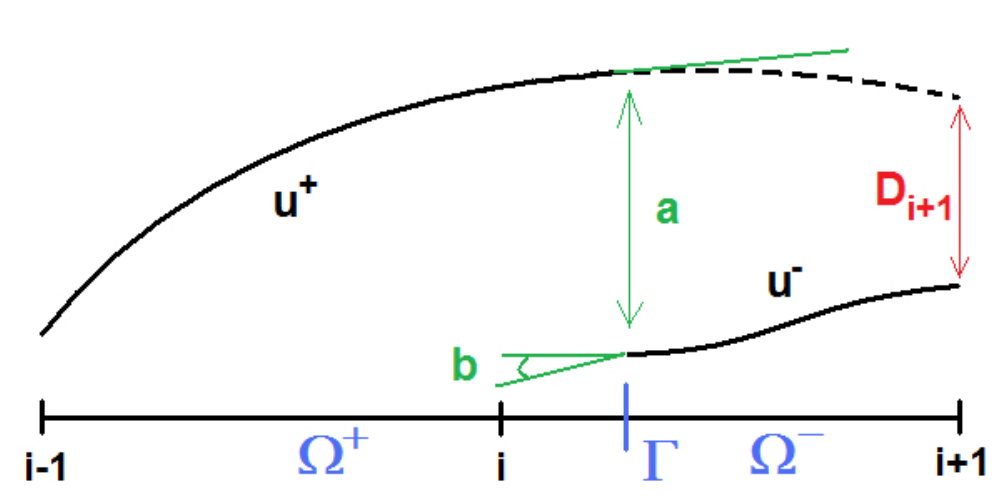}
 \end{center}
 \caption{Example in 1D of a solution with a jump discontinuity.}
 \label{fig:uxx}
\end{figure}

Let us first illustrate the key concept in the GFM with a simple example,
involving a regular grid and the standard second order discretization of the
1D analog of the problem we are interested in. Thus, consider the problem of
discretizing the equation $u_{xx} = f(x)$ in some interval
$\Omega = \{x : x_L < x < x_R\}$, where $u$ is discontinuous across some point
$x_{\Gamma}$ --- hence $\Omega^+$ (respectively $\Omega^-$) is the domain
$x_L < x < x_{\Gamma}$ (respectively $x_{\Gamma} < x < x_R$). Then, see
figure~\ref{fig:uxx}, when trying to approximate $u_{xx}$ at a grid point $x_i$
such that $x_i < x_{\Gamma} < x_{i+1}$, we would like to write
\begin{equation}\label{eq:uxx}
 u_{xx_i}^+ \approx \dfrac{u_{i-1}^+ -2u_i^+ + u_{i+1}^+}{h^2}\/,
\end{equation}
where $h = x_{j+1}-x_j$ is the grid spacing. However, we do not have information
on $u_{i+1}^+$, but rather on $u_{i+1}^-$. Thus, the idea is to estimate a
correction for $u_{i+1}^-$, to recover $u_{i+1}^+$, such that Eq.~\eqref{eq:uxx}
can be applied:
\begin{equation}
 u_{xx_i}^+ \approx \dfrac{u_{i-1}^+ -2u_i^+ +
 \overbrace{\p{u_{i+1}^- + D_{i+1}}}^{u_{i+1}^+}}{h^2},
\end{equation}
where $D_{i+1} = u_{i+1}^+ - u_{i+1}^-$ is the correction term. Now we note that,
if $D_{i+1}$ can be written as a correction that is independent on the solution
$u$, then it can be moved to the RHS of the equation, and absorbed into $f$.
That is
\begin{equation}
 \dfrac{u_{i-1}^+ -2u_i^+ + u_{i+1}^-}{h^2} = f_i - \dfrac{D_{i+1}}{h^2}.
\end{equation}
This allows the solution of the problem with prescribed discontinuities using
the same discretization as the one employed to solve the simpler problem
without an interface --- which leads to a great efficiency gain.

\begin{rmk}\label{rem:3:1}
 The error in estimating 
 $D\/$ is crucial in determining the accuracy of the final discretization. Liu,
 Fedkiw, and Kang~\cite{liu:00} introduced a dimension-by-dimension linear
 extrapolation of the interface jump conditions, to get a first order
 approximation for $D\/$. Our new
 method is based on generalizing the idea of a correction term to that of a
 \emph{correction function}, for which we can write an equation. One can then
 obtain high accuracy representations for $D\/$ by solving this equation,
 without the complications into which dimension-by-dimension (with Taylor
 expansions) approaches run into.\myremarkend
\end{rmk}

\begin{rmk}\label{rem:3:2}
 An additional advantage of the correction function approach is that $D\/$ can
 be calculated at any point near the interface $\Gamma\/$. Hence it can be used
 with any finite differences discretization of the Poisson equation, without
 regard to the particulars of its stencil (as would be the case with any
 approach based on Taylor expansions).\myremarkend
\end{rmk}

\section{The correction function and the equation defining it.}
\label{sec:general}
As mentioned earlier, the aim here is to generalize the correction term concept
to that of a \emph{correction function}, and then to find an equation (a PDE, with
appropriate boundary conditions) that uniquely characterizes the \emph{correction
function}. Then, at least in principle, one can design algorithms to solve the
PDE in order to obtain solutions to the \emph{correction function} of any desired
order of accuracy.

Let us begin by considering a small region $\Omega_\Gamma$ enclosing the
interface $\Gamma$, defined as the set of all the points within some
distance $\mathcal{R}$ of $\Gamma$, where $\mathcal{R}$ is of the order of the grid
size $h$. As we will see below, we would like to have $\mathcal{R}$ as small as
possible. On the other hand, $\Omega_{\Gamma}$ has to include all the points
where the CFM requires corrections to be computed, which means\footnote{For the
particular discretization of the Laplace operator that we use in this paper.} that
$\mathcal{R}$ cannot be smaller than $\sqrt{2}h$. In addition, algorithmic
considerations (to be seen later) force $\mathcal{R}$ to be slightly larger than
this last value.

Next, we assume we that can extrapolate both $u^+$ and $u^-$, so that they are
valid everywhere within $\Omega_{\Gamma}$, in such a way that they satisfy the
Poisson equations
\begin{subequations}\label{eq:up&um}
 \begin{align}
  \nabla^2u^+\p{\vec{x}} &= f^+\p{\vec{x}}
  &\mathrm{for}\;\; \vec{x} &\in \Omega_{\Gamma}\/,\\
  \nabla^2u^-\p{\vec{x}} &= f^-\p{\vec{x}}
  &\mathrm{for}\;\; \vec{x} &\in \Omega_{\Gamma}\/,
 \end{align}
\end{subequations}
where $f^+$ and $f^-$ are smooth enough (see remark~\ref{rem:4:1} below)
extensions of the source term $f$ to $\Omega_{\Gamma}$. In particular, notice
that the introduction of $f^+$ and $f^-$ allows the possibility of the
source term changing (\ie\/ a discontinuous source term) across
$\Gamma$. The \emph{correction function} is then defined by
$D(\vec{x}) = u^+(\vec{x})-u^-(\vec{x})$.

Taking the difference between the equations in \eqref{eq:up&um}, and using the
jump conditions (\ref{eq:a}-\ref{eq:b}), yields
\begin{subequations}\label{eq:D}
 \begin{align}
  \nabla^2D\p{\vec{x}} &= f^+\p{\vec{x}} - f^-\p{\vec{x}} =
  f_D\p{\vec{x}} &\mathrm{for}\;\; \vec{x} &\in \Omega_{\Gamma},
  \label{eq:D-poisson}\\
  D\p{\vec{x}} &= a\p{\vec{x}} &\mathrm{for}\;\; \vec{x}
  &\in \Gamma,\label{eq:D-a}\\
  D_n\p{\vec{x}} &= b\p{\vec{x}} &\mathrm{for}\;\; \vec{x}
  &\in \Gamma.\label{eq:D-b}
 \end{align}
\end{subequations}
This achieves the aim of having the \emph{correction function} defined by a set of
equations, with some provisos --- see remark~\ref{rem:4:2} below. Note that:
\begin{enumerate}
 \item If $f^+\p{\vec{x}} = f^-\p{\vec{x}}\/$, for
 $\vec{x} \in \Omega_{\Gamma}\/$, then $f_D\p{\vec{x}} = 0\/$, for
 $\vec{x} \in \Omega_{\Gamma}\/$.
 \item Equation~\eqref{eq:D-b} imposes the true jump condition
 in the normal direction, whereas some versions of the GFM rely on a
 dimension-by-dimension approximation of this condition (see
 Ref.~\cite{liu:00}).
\end{enumerate}

\begin{rmk}\label{rem:4:1}
 The smoothness requirement on $f^+$ and $f^-$ is tied up to how accurate an
 approximation to the correction term $D$ is needed. For example, if a \fth\/
 order algorithm is used to find $D$, this will (generally) necessitate $D$ to
 be at least $C^4$ for the errors to actually be \fth\/ order. Hence, in this
 case, $f_D=f^+-f^-$ must be $C^2\/$.\myremarkend
\end{rmk}

\begin{rmk}\label{rem:4:2}
 Equation~\eqref{eq:D} is an elliptic Cauchy problem for $D$ in
 $\Omega_{\Gamma}$. In general, such problems are ill-posed. However, we are
 seeking for solutions within the context of a numerical approximation where
 \begin{enumerate}[(a)]
  \item There is a frequency cut-off in both the data $a = a(\vec{x})$ and
  $b = b(\vec{x})$, and the description of the curve $\Gamma$.
  \item We are interested in the solution only a small distance away from the
  interface $\Gamma$, where this distance vanishes simultaneously with the
  inverse of the cut-off frequency in point (a).
 \end{enumerate}
 What (a) and (b) mean is that the arbitrarily large growth rate for
 arbitrarily small perturbations, which is responsible for the ill-posedness
 of the Cauchy problem in Eq.~\eqref{eq:D}, does not occur within the special
 context where we need to solve the problem. This large growth rate does not
 occur because, for the solutions of the Poisson equation, the growth rate
 for a perturbation of wave number $0 < k < \infty$ along some straight line,
 is given by $e^{2\/\pi\/k\/d}$ --- where $d$ is the distance from the line.
 However, by construction, in the case of interest to us $k\/d$ is bounded.\myremarkend
\end{rmk}

\begin{rmk}\label{rem:4:3}
 Let us be more precise, and define a number characterizing how well posed the
 discretized version of Eq.~\eqref{eq:D} is, by
 \begin{equation*}
  \alpha = \text{largest growth rate possible,}
 \end{equation*}
 where growth is defined relative to the size of a perturbation to the solution
 on the interface. This number is determined by $\mathcal{R}$ (the ``radius'' of
 $\Omega_\Gamma$) as the following calculation shows: First of all, there is no
 loss of generality in assuming that the interface is flat, provided that the
 numerical grid is fine enough to resolve $\Gamma$. In this case, let us
 introduce an orthogonal coordinate system $\vec{y}$ on $\Gamma$, and let $d$
 be the signed distance to $\Gamma$ (say, $d > 0$ in $\Omega^-$). Expanding the
 perturbations in Fourier modes along the interface, the typical mode has the
 form
 \begin{equation*}
  \varphi_{\vec{k}} = e^{2\/\pi\/i\/\vec{k}\cdot\vec{y} \pm 2\/\pi\/k\/d},
 \end{equation*}
 where $\vec{k}$ is the Fourier wave vector, and $k = |\vec{k}|$. The shortest
 wave-length that can be represented on a grid with mesh size $0 < h \ll 1$
 corresponds to $k = k_{\max} = 1/(2\/h)$. Hence, we obtain the estimate
 \begin{equation*}
  \alpha \approx e^{\pi\/\mathcal{R}_c/h}
 \end{equation*}
 for the maximum growth rate. \myremarkend
\end{rmk}

\begin{rmk}\label{rem:4:4}
 Clearly, $\alpha$ is intimately related to the condition number for the
 discretized problem --- see \S~\ref{sec:scheme}. In fact, at leading
 order, the two numbers should be (roughly) proportional to each other ---
 with a proportionality constant that depends on the details of the
 discretization. For the discretization used in this paper (described further
 below), $\sqrt{2}h \le \mathcal{R} \le 2\sqrt{2}\/h$, which leads to the rough estimate
 $85 < \alpha < 7,200$. On the other hand, the observed condition numbers vary between
 5,000 and 10,000. Hence, the actual condition numbers are only slightly higher than
 $\alpha$ for the ranges of grid sizes $h\/$ that we used (we did not explore the asymptotic limit $h \to 0\/$).
\myremarkend
\end{rmk}

\begin{rmk}\label{rem:4:5}
 Equation~\eqref{eq:D} depends on the known inputs for the problem only. Namely:
 $f^+$, $f^-$, $a$, and $b$. Consequently  $D$ does not depend on the solution
 $u$. Hence, after solving for $D$, we can use a discretization for $u$ that
 does not involve the interface: Whenever $u$ is discontinuous, we evaluate
 $D$ where the correction is needed, and transfer these values to the RHS.\myremarkend
\end{rmk}

\begin{rmk}\label{rem:4:6}
 When developing an algorithm for a linear Cauchy problem, such as the one in
 Eq.~\eqref{eq:D}, the two key requirements are consistency and stability. In
 particular, when the solution depends on the ``initial
 conditions'' globally, stability (typically) imposes stringent constraints on
 the ``time'' step for any local (explicit) scheme. This would seem to suggest
 that, in order to solve Eq.~\eqref{eq:D}, a ``global'' (involving the whole
 domain $\Omega_\Gamma\/$) method will be needed. This, however, is not true:
 because we need to solve Eq.~\eqref{eq:D} for one ``time'' step only ---
 \ie\/~within an $\mathcal{O}(h)\/$ distance from $\Gamma\/$, stability
 is not relevant. Hence, consistency is enough, and a fully local scheme is
 possible. In the algorithm described in \S~\ref{sec:scheme} we found that,
 for (local) quadrangular patches, the Cauchy problem leads to a well behaved
 algorithm when the length of the interface contained in each patch is of the
 same order as the diagonal length of the patch. This result is in line with
 the calculation in remark~\ref{rem:4:3}: we want to keep the ``wavelength''
 (along $\Gamma\/$) of the perturbations introduced by the discretization as
 long as possible. In particular, this should then minimize the condition
 number for the local problems --- see remark~\ref{rem:4:4}.\myremarkend
\end{rmk}

\section{A \fth\/ Order Accurate Scheme in 2D.} \label{sec:scheme}
%
\subsection{Overview.}
In this section we use the general ideas presented earlier to develop a
specific example of a \fth\/ order accurate scheme in 2D. Before proceeding with
an in-depth description of the scheme, we highlight a few key points:
\begin{enumerate}[(a)]
 \item We discretize Poisson's equation using a compact 9-point stencil.
 Compactness is important since it is directly related to the size of
 $\mathcal{R}_c$, which has a direct impact on the problem's conditioning
 --- see remarks~\ref{rem:4:2} -- \ref{rem:4:4}.
 \item We approximate $D$ using bicubic interpolations (bicubics), each
 valid in a small neighborhood $\Omega_{\Gamma}^{i,j}$ of the interface. This
 guarantees local \fth\/ order accuracy with only 12 interpolation parameters
 --- see~\cite{nave:10}. Each $\Omega_{\Gamma}^{i,j}$ corresponds to a point
 in the grid at which the standard discretization of Poisson's equation
 involves a stencil that straddles the interface $\Gamma\/$.
 \item The domains $\Omega_{\Gamma}^{i,j}$ are rectangular regions, each
 enclosing a portion of $\Gamma$, and all the nodes where $D$ is needed
 to complete the discretization of the Poisson equation at the
 ($i\/,\,j$)-th stencil. Each is a sub-domain of $\Omega_\Gamma\/$.
 \item Starting from (b) and (c), we design a local solver that provides an
 approximation to $D$ inside each domain $\Omega_{\Gamma}^{i,j}$.
 \item The interface $\Gamma$ is represented using the Gradient-Augmented
 Level Set approach --- see~\cite{nave:10}. This guarantees a local
 \fth\/ order representation of the interface, as required to keep the overall
 accuracy of the scheme.
 \item In each $\Omega_{\Gamma}^{i,j}$, we solve the PDE in \eqref{eq:D} in a
 least squares sense. Namely: First we define an appropriate positive quadratic
 integral quantity $J_P\/$ --- Eq.~\eqref{eq:jp} --- for which the solution is a minimum (actually, zero).
 Next we substitute the bicubic approximation for the solution into $J_P\/$,
 and discretize the integrals using Gaussian quadrature. Finally, we find
 the bicubic parameters by minimizing the discretized $J_P\/$.
\end{enumerate}

\begin{rmk}\label{rem:5:1}
 Solving the PDE in a least squares sense is crucial, since an algorithm is
 needed that can deal with the myriad ways in which the interface $\Gamma\/$
 can be placed relative to the fixed rectangular grid used to discretize
 Poisson's equation. This approach provides a scheme that (i) is robust with
 respect to the details of the interface geometry, (ii) has a formulation
 that is (essentially) dimension independent --- there are no fundamental changes from
 2D to 3D, and (iii) has a clear theoretical underpinning that allows extensions
 to higher orders, or to other discretizations of the Poisson equation.\myremarkend
\end{rmk}

\subsection{Standard Stencil.}
We use the standard \fth\/ order accurate 9-point discretization of Poisson's
equation\footnote{Notice that here we allow for the possibility of different
   grid spacings in each direction.}:
\begin{equation}\label{eq:9p}
 L^5 u_{i,j} + \dfrac{1}{12}\p{h_x^2 + h_y^2}\hat{\partial}_{xx}\hat{\partial}_{yy} u_{i,j} =
 f_{i,j} + \dfrac{1}{12}\p{h_x^2\,(f_{xx})_{i,j} +  h_y^2\,(f_{yy})_{i,j}}\/,
\end{equation}
where $L^5$ is the second order 5-point discretization of the Laplace operator:
\begin{equation}\label{eq:L5}
 L^5 u_{i,j} = \hat{\partial}_{xx} u_{i,j} + \hat{\partial}_{yy} u_{i,j},
\end{equation}
and
\begin{align}
 \hat{\partial}_{xx} u_{i,j} &= \dfrac{u_{i+1,j} - 2u_{i,j} + u_{i-1,j}}{h_x^2},\label{eq:partialxx}\\
 \hat{\partial}_{yy} u_{i,j} &= \dfrac{u_{i,j+1} - 2u_{i,j} + u_{i,j-1}}{h_y^2}.\label{eq:partialyy}
\end{align}
The terms $(f_{xx})_{i,j}\/$ and  $(f_{yy})_{i,j}\/$ may be given analytically
(if known), or computed using appropriate second order discretizations.

In the absence of discontinuities, Eq.~\eqref{eq:9p} provides a compact
\fth\/ order accurate representation of Poisson's equation. In the vicinity
of the discontinuities at the interface $\Gamma\/$, we define an appropriate
domain $\Omega_{\Gamma}^{i,j}$, and compute the correction terms necessary to
Eq.~\eqref{eq:9p} --- as described in detail next.

\begin{figure}[htb!]
 \begin{center}
  \includegraphics[width=2.5in]{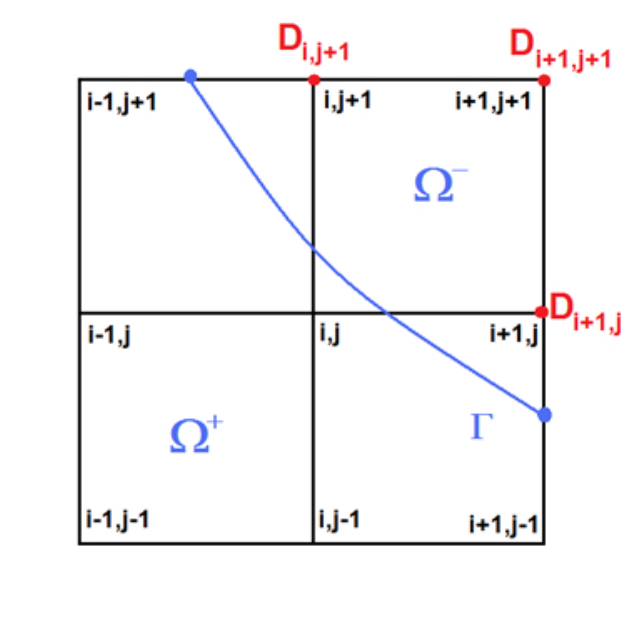}
 \end{center}
 \caption{The 9-point compact stencil next to the interface $\Gamma$.}
 \label{fig:D}
\end{figure}

To understand how the correction terms affect the discretization, let us
consider the situation depicted in figure~\ref{fig:D}. In this case, the node
$(i,j)$ lies in $\Omega^+$ while the nodes $(i+1,j)$, $(i+1,j+1)$, and
$(i,j+1)$ are in $\Omega^-$. Hence, to be able to use Eq.~\eqref{eq:9p}, we
need to compute $D_{i+1,j}$, $D_{i+1,j+1}$, and $D_{i,j+1}$.

After having solved for $D$ where necessary (see \S~\ref{sub:omega} and
\S~\ref{sub:solution}), we modify Eq.~\eqref{eq:9p} and write
\begin{equation}\label{eq:9p-C}
  L^5 u_{i,j} + \dfrac{1}{12}\p{h_x^2 + h_y^2}\hat{\partial}_{xx}\hat{\partial}_{yy}u_{i,j} =
  f_{i,j} + \dfrac{1}{12}\p{h_x^2\,(f_{xx})_{i,j} + h_y^2\,(f_{yy})_{i,j}} +
  C_{i,j},
\end{equation}
which differs from Eq.~\eqref{eq:9p} by the terms $C_{i,j}$ on the RHS only.
Here the $C_{i,j}$ are the CFM correction terms needed to complete the stencil
across the discontinuity at $\Gamma\/$. In the particular case illustrated in
figure~\ref{fig:D}, we have
\begin{equation}\label{eq:9p-D}
 \begin{split}
  C_{i,j} & = \left[\dfrac{1}{6}\dfrac{\p{h_x^2 + h_y^2}}{\p{h_xh_y}^2} -
             \dfrac{1}{h_y^2}\right]D_{i+1,j} +
             \left[\dfrac{1}{6}\dfrac{\p{h_x^2 + h_y^2}}{\p{h_xh_y}^2} -
             \dfrac{1}{h_x^2}\right]D_{i,j+1} \\
  & - \dfrac{1}{12}\dfrac{\p{h_x^2 + h_y^2}}{\p{h_xh_y}^2}D_{i+1,j+1}\/.
 \end{split}
\end{equation}
Similar formulas apply for the other possible arrangements of the Poisson's
equation stencil relative to the interface $\Gamma\/$. 

\subsection{Definition of $\Omega_{\Gamma}^{i,j}\/$.} \label{sub:omega}
There is some freedom on how to define $\Omega_{\Gamma}^{i,j}$. The basic requirements are
\begin{enumerate}[(i)]
 \item $\Omega_{\Gamma}^{i,j}$ should be a rectangle.
 \item the edges of $\Omega_{\Gamma}^{i,j}$ should be parallel to
 the grid lines.
 \item $\Omega_{\Gamma}^{i,j}$ should be small, since the problem's condition
 number increases exponentially with the distance from $\Gamma\/$ --- see
 remarks \ref{rem:4:3} and \ref{rem:4:4}.
 \item $\Omega_{\Gamma}^{i,j}$ should contain all the nodes where $D$ is needed.
 For the example, in figure~\ref{fig:D} we need to know $D_{i+1,j+1}$, $D_{i+1,j}$,
 and $D_{i,j+1}$. Hence, in this case, $\Omega_{\Gamma}^{i,j}$ should include the
 nodes $(i+1,j+1)$, $(i+1,j)$, and $(i,j+1)$.
 \item $\Omega_{\Gamma}^{i,j}$ should contain a segment of $\Gamma$, with a
 length that is as large as possible --- \ie\/~comparable to the length
 of the diagonal of $\Omega_{\Gamma}^{i,j}$. This follows from the calculation
 in remark~\ref{rem:4:3}, which indicates that the wavelength of the
 perturbations (along $\Gamma\/$) introduced by the discretization should be
 as long as possible. This should then minimize the condition number for the
 local problem --- see remark~\ref{rem:4:4}.
\end{enumerate}

Requirements (i) and (ii) are needed for algorithmic convenience only, and do not
arise from any particular argument in \S~\ref{sec:general}. Thus, in principle,
this convenience could be traded for improvements in other areas --- for example,
for better condition numbers for the local problems, or for additional flexibility
in dealing with complex geometries. However, for simplicity, in this paper we enforce
(i) and (ii). As explained earlier (see remark~\ref{rem:5:1}), we solve Eq.~\eqref{eq:D}
in a least squares sense. Hence integrations over $\Omega_{\Gamma}^{i,j}$ are required.
It is thus useful to keep $\Omega_{\Gamma}^{i,j}$ as simple as possible.

A discussion of various aspects regarding the proper definition of $\Omega_{\Gamma}^{i,j}$
can be found in appendix~\ref{ap:omega}. For instance, the requirement in item~(ii) is
convenient only when an implicit representation of the interface is used. Furthermore,
although the definition of $\Omega_{\Gamma}^{i,j}$ presented here proved robust for all
the applications of the \fth\/ order accurate scheme (see \S~\ref{sec:results}), there are
specific geometrical arrangements of the interface for which (ii) results in extremely
elongated $\Omega_{\Gamma}^{i,j}$. These elongated geometries can have negative effects
on the accuracy of the scheme. We noticed this effect in the \snd\/ order accurate version
of the method described in appendix~\ref{ap:scheme2}. These issues are addressed by the
various algorithms (of increasing complexity) presented in appendix~\ref{ap:omega}.

\begin{figure}[htb!]
 \begin{center}
  \includegraphics[width=2.5in]{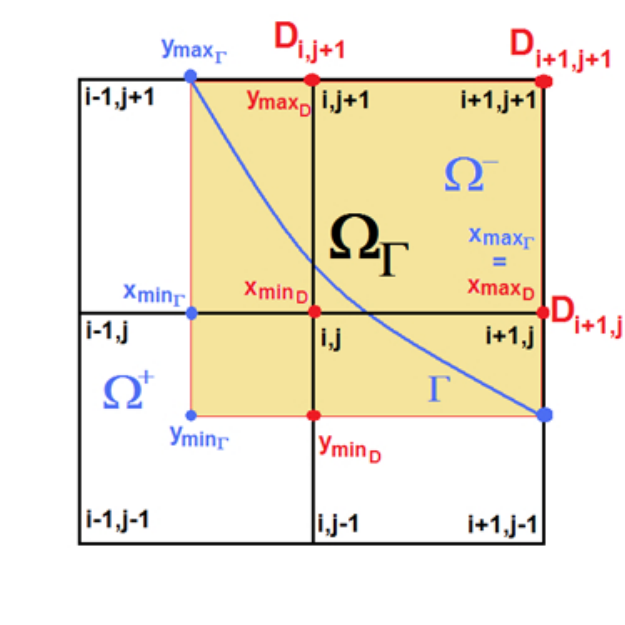}
 \end{center}
 \caption{The set $\Omega_{\Gamma}^{i,j}$ for the situation in
 figure~\ref{fig:D}.}
 \label{fig:omegag}
\end{figure}

With the points above in mind, here (for simplicity) we define $\Omega_{\Gamma}^{i,j}$ as the smallest rectangle that satisfies the requirements in (i), (ii), (iv), and (v) --- then (iii) follows automatically. Hence $\Omega_{\Gamma}^{i,j}$ can be constructed using the following three easy steps:
\begin{enumerate}[1.]
 \item Find the coordinates ($x_{\min_{\Gamma}}\/,\,x_{\max_{\Gamma}}$) and
 ($y_{\min_{\Gamma}}\/,\,y_{\max_{\Gamma}}$) of the smallest rectangle satisfying
 condition (ii), which completely encloses the section of the interface
 $\Gamma$ contained by the region covered by the 9-point stencil.
 \item Find the coordinates ($x_{\min_D}\/,\,x_{\max_D}$) and
 ($y_{\min_D}\/,\,y_{\max_D}$) of the smallest rectangle satisfying condition
 (ii), which completely encloses all the nodes at which $D$ needs to be known.
 \item Then  $\Omega_{\Gamma}^{i,j}$ is the smallest rectangle that encloses the
 two previous rectangles. Its edges are given by
  \begin{subequations}
   \begin{align}
    x_{\min} &= \min\p{x_{\min_{\Gamma}},x_{\min_D}},\\
    x_{\max} &= \max\p{x_{\max_{\Gamma}},x_{\max_D}},\\
    y_{\min} &= \min\p{y_{\min_{\Gamma}},y_{\min_D}},\\
    y_{\max} &= \max\p{y_{\max_{\Gamma}},y_{\max_D}}.
   \end{align}
  \end{subequations}
\end{enumerate}
Figure~\ref{fig:omegag} shows an example of $\Omega_{\Gamma}^{i,j}$ defined
using these specifications.

\begin{rmk}\label{rem:5:2}
Notice that for each node next to the interface we construct a domain
$\Omega_{\Gamma}^{i,j}$. When doing so, we allow the domains to overlap. For
example, the domain $\Omega_{\Gamma}^{i,j}$ shown in figure~\ref{fig:omegag} is
used to determine $C_{i,j}$. It should be clear that $\Omega_{\Gamma}^{i-1,j+1}$
(used to determine $C_{i-1,j+1}$), and  $\Omega_{\Gamma}^{i+1,j-1}$ (used to
determine $C_{i+1,j-1}$), each will overlap with $\Omega_{\Gamma}^{i,j}$.

The consequence of these overlaps is that different computed values for $D$
at the same node can (in fact, will) happen --- depending on which domain is
used to solve the local Cauchy problem. However, because we solve for $D$ ---
within each $\Omega_{\Gamma}^{i,j}$ --- to \fth\/ order accuracy, any differences
that arise from this multiple definition of $D$ lie within the order of
accuracy of the scheme. Since it is convenient to keep the computations
local, the values of $D\/$ resulting from the domain $\Omega_{\Gamma}^{i,j}$,
are used to evaluate the correction term $C_{i,j}$.\myremarkend
\end{rmk}

\begin{rmk}\label{rem:5:3}
While rare, cases where a single interface crosses the same stencil multiple times can occur.
In \S~\ref{sub:case2} we present such an example. A simple approach to deal with situations
like this is as follows: First associate each node where the correction function is needed to
a particular piece of interface crossing the stencil (say, the closest one). Then define one
$\Omega_{\Gamma}^{i,j}$ for each of the individual pieces of interface crossing the stencil.
 
For example, figure~\ref{fig:omegag-mult}(a) depicts a situation where the stencil is crossed
by two pieces of the same interface ($\Gamma_1$ and $\Gamma_2$), with $D$ needed at the nodes
$(i+1,j+1)$, $(i+1,j)$, $(i,j+1)$, and $(i-1,j-1)$. Then, first associate: (i)~$(i+1,j+1)$,
$(i+1,j)$, and $(i,j+1)$ to $\Gamma_1$, and (ii) $(i-1,j-1)$ to $\Gamma_2$. Second, define
\begin{enumerate}[1.]
\item
$\Omega_{\Gamma_1}^{i,j}$ is the smallest rectangle, parallel to the grid lines, that includes $\Gamma_1$ and the nodes $(i+1,j+1)$, $(i+1,j)$, and $(i,j+1)$.
\item
$\Omega_{\Gamma_2}^{i,j}$ is the smallest rectangle, parallel to the grid lines, that includes $\Gamma_2$ and the node $(i-1,j-1)$.
\end{enumerate} 
\begin{figure}[htb!]
 \begin{center}
  \includegraphics[width=2.5in]{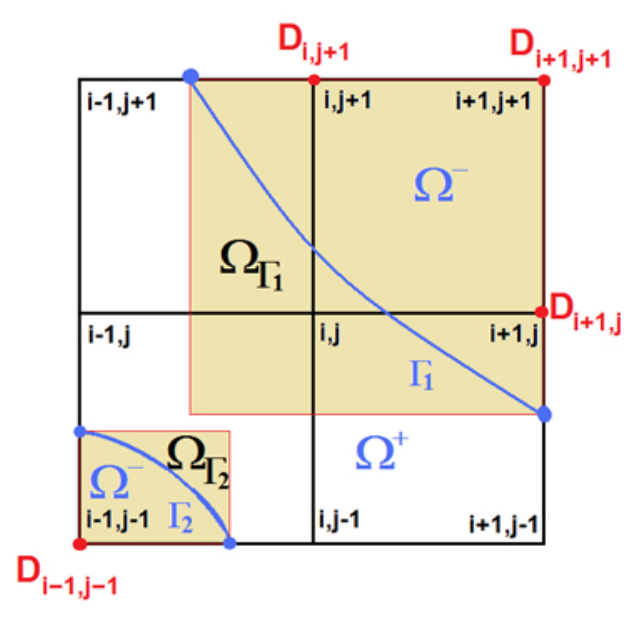}
  \includegraphics[width=2.5in]{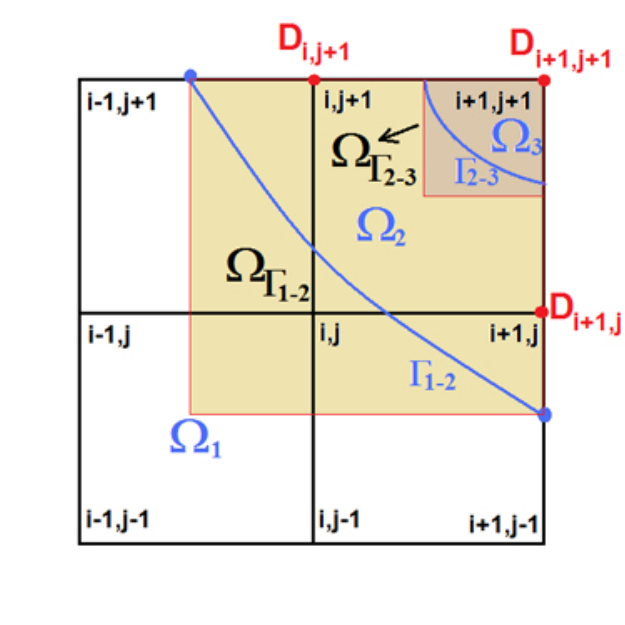}\\
  \text{(a)} \hspace{2.3in} \text{(b)}
 \end{center}
 \caption{Configuration where multiple $\Omega_{\Gamma_2}^{i,j}$ are needed in the same stencil. (a) Same interface crossing the stencil multiple times. (b) Distinct interfaces crossing the same stencil.}
 \label{fig:omegag-mult}
\end{figure}

After the multiple $\Omega_{\Gamma}^{i,j}$ are defined within a given stencil,
the local Cauchy problem is solved for each $\Omega_{\Gamma}^{i,j}$ separately.
For example, in the case shown in figure~\ref{fig:omegag-mult}(a), the solution
for $D$ inside $\Omega_{\Gamma_1}^{i,j}$ is done completely independent of the
solution for $D$ inside $\Omega_{\Gamma_2}^{i,j}$. The decoupling between multiple
crossings renders the CFM flexible and robust enough 
to handle complex geometries without any special algorithmic considerations.
\myremarkend
\end{rmk}

\begin{rmk}\label{rem:5:4}
When multiple distinct interfaces are involved, a single stencil can be crossed
by different interfaces -- \eg\/: see \S~\ref{sub:case4} and \S~\ref{sub:case5}.
This situation is similar to the one described in remark~\ref{rem:5:3}, but
with an additional complication: there may occur distinct domain regions that
are not separated by an interface, but rather by a third (or more) regions
between them. An example is shown in figure~\ref{fig:omegag-mult}(b), where
$\Gamma_{1-2}$ and $\Gamma_{2-3}$ are not part of the same interface. Here
$\Gamma_{1-2}$ is the interface between $\Omega_1$ and $\Omega_2$, while
$\Gamma_{2-3}$ is the interface between $\Omega_2$ and $\Omega_3$. There is no
interface $\Gamma_{1-3}$ separating $\Omega_1$ from $\Omega_3$, hence no jump
conditions between these regions are provided. Nonetheless,
$D_{1-3} = (u_{3} - u_{1})$ is needed at $(i+1,j+1)$.

Situations such as these can be easily handled by noticing that we can
distinguish between primary (\eg\/~$D_{1-2}$ and $D_{2-3}$) and secondary
correction functions, which can be written in terms of the primary functions
\linebreak (\eg\/~$D_{1-3} = D_{1-2} + D_{2-3}$) and need not be computed directly. Hence we
can proceed exactly as in remark~\ref{rem:5:3}, except that we have to make
sure that the intersections of the regions where the primary correction
functions are computed include the nodes where the secondary correction
functions are needed. For example, in the particular case in
figure~\ref{fig:omegag-mult}(b), we define
\begin{enumerate}[1.]
\item
$\Omega_{\Gamma_{1-2}}^{i,j}$ is the smallest rectangle, parallel to the grid lines, that includes $\Gamma_{1-2}$ and the nodes $(i+1,j+1)$, $(i+1,j)$, and $(i,j+1)$.
\item
$\Omega_{\Gamma_{2-3}}^{i,j}$ is the smallest rectangle, parallel to the grid lines, that includes $\Gamma_{2-3}$ and the node $(i+1,j+1)$.
\end{enumerate} 
\myremarkend
\end{rmk}
%
\subsection{Solution of the Local Cauchy Problem.} \label{sub:solution}
Since we use a \fth\/ order accurate discretization of the Poisson problem, we
need to find $D$ with \fth\/ order errors (or better) to keep the overall
accuracy of the scheme --- see \S~\ref{sub:erroranalysis}. Hence we
approximate $D$ using cubic Hermite splines (bicubic interpolants in
2D), which guarantees \fth\/ order accuracy --- see~\cite{nave:10}. Note also
that, even though the example scheme developed here is for 2D, this
representation can be easily extended to any number of dimensions.

Given a choice of basis functions,\footnote{The basis functions that we use
for the bicubic interpolation can be found in appendix~\ref{ap:bicubic}.}
we solve the local Cauchy problem defined in Eq.~\eqref{eq:D} in a least
squares sense, using a minimization procedure. Since we do not have boundary
conditions, but interface conditions, we must resort to a minimization
functional that is different from the standard one associated with the Poisson
equation. Thus we impose the Cauchy interface conditions by using a
penalization method. The functional to be minimized is then
\begin{equation}\label{eq:jp}
 \begin{split}
  J_P & = (\ell_c^{i,j})^3\int_{\Omega_{\Gamma}^{i,j}} \left[\nabla^2D\p{\vec{x}}
          - f_D\p{\vec{x}}\right]^2d\/V \\
      & + c_P\int_{\Gamma\cap\Omega_{\Gamma}^{i,j}}\left[D\p{\vec{x}}
        - a\p{\vec{x}}\right]^2dS\\
      & + c_P(\ell_c^{i,j})^2\int_{\Gamma\cap\Omega_{\Gamma}^{i,j}}
          \left[D_n\p{\vec{x}} - b\p{\vec{x}}\right]^2d\/S\/,
 \end{split}
\end{equation}
where $c_P > 0\/$ is the penalization coefficient used to enforce the interface
conditions, and $\ell_c^{i,j} > 0$ is a characteristic length associated with
$\Omega_{\Gamma}^{i,j}$ --- we used the shortest side length. Clearly $J_P\/$ is
a quadratic functional whose minimum (zero) occurs at the solution to
Eq.~\eqref{eq:D}.

In order to compute $D\/$ in the domain $\Omega_{\Gamma}^{i,j}$, its bicubic
representation is substituted into the formula above for $J_P\/$, with the
integrals approximated by Gaussian quadratures --- in this paper we used six
quadrature points for the 1D line integrals, and 36 points for the 2D area
integrals. The resulting discrete problem is then minimized. Because the
bicubic representation for $D\/$ involves 12 basis polynomials, the
minimization problem produces a $12\times12$ (self-adjoint) linear system.

\begin{rmk}\label{rem:5:5}
 We explored the option of enforcing the
 interface conditions using Lagrange multipliers. While this second approach
 yields good results, our experience shows that the penalization method is
 better.\myremarkend
\end{rmk}

\begin{rmk}\label{rem:5:6}
 The scaling using $\ell_c^{i,j}$ in Eq.~\eqref{eq:jp} is so that all the three
 terms in the definition of $J_p\/$ behave in the same fashion as the size of
 $\Omega_{\Gamma}^{i,j}$ changes with ($i\/,\,j$), or when the computational
 grid is refined.\footnote{The scaling also follows from dimensional
    consistency.}
 This follows because we expect that
 \begin{align*}
  \nabla^2 D - f = \mathcal{O}\p{\ell_c^2},\\
  D - a = \mathcal{O}\p{\ell_c^4},\\
  D_n - b = \mathcal{O}\p{\ell_c^3}.
 \end{align*}
 Hence each of the three terms in Eq.~\eqref{eq:jp} should be
 $\mathcal{O}\p{\ell_c^{9}}\/$.\myremarkend
\end{rmk}

\begin{rmk}\label{rem:5:7}
 Once all the terms in Eq.~\eqref{eq:jp}
 are guaranteed to scale the same way with the size of $\Omega_{\Gamma}^{i,j}$,
 the penalization coefficient $c_P\/$ should be selected so that the three terms
 have (roughly) the same size for the numerical solution (they will, of course,
 not vanish).
 In principle, $c_P\/$ could be determined from knowledge of the fourth order
 derivatives of the solution, which control the error in the numerical solution.
 This approach does not appear to be practical. A simpler method is based on the
 observation that $c_P\/$ should not depend on the grid size (at least to
 leading order, and we do not need better than this). Hence it can be determined
 empirically from a low resolution calculation.
 In the examples in this paper we found that $c_P \approx 50$ produced good
 results.\myremarkend
\end{rmk}

\begin{rmk}\label{rem:5:8}
 A more general version of $J_P\/$ in
 Eq.~\eqref{eq:jp} would involve different penalization coefficients for the
 two line integrals, as well as the possibility of these coefficients having a
 dependence on the position along $\Gamma\/$ of $\Omega_{\Gamma}^{i,j}\/$. These
 modifications could be useful in cases where the solution to the Poisson
 problem has large variations --- \eg\/~a very irregular interface
 $\Gamma\/$, or a complicated forcing $f\/$.\myremarkend
\end{rmk}

\subsection{Computational Cost.}
We can now infer something about the cost of the present scheme. To start with,
let us denote the number of nodes in the $x\/$ and $y\/$ directions by
\begin{align}
 N_x = \dfrac{1}{h_x}+1, && N_y = \dfrac{1}{h_y}+1\/,
\end{align}
assuming a 1 by 1 computational square. Hence, the total number of degrees of
freedom is $M=N_x\/N_y\/$. Furthermore, the number of nodes adjacent to the
interface is $\mathcal{O}(M^{1/2})\/$, since the interface is a 1D entity.

The discretization of Poisson's equation results in a $M\times M\/$ linear
system. Furthermore, the present method produces changes only on the RHS of
the equations. Thus, the basic cost of inverting the system is unchanged from
that of inverting the system resulting from a problem without an interface
$\Gamma\/$. Namely: it varies from $\mathcal{O}(M)$ to $\mathcal{O}(M^2)$,
depending on the solution method.

Let us now consider the computational cost added by the modifications to the
RHS. As presented above, for each node adjacent to the interface, we must
construct $\Omega_{\Gamma}^{i,j}$, compute the integrals that define the local
$12\times12$ linear system, and invert it. The cost associated with these tasks
is constant: it does not vary from node to node, and it does not change with
the size of the mesh. Consequently the resulting additional cost is a constant
times the number of nodes adjacent to the interface. Hence it scales as
$M^{1/2}$. Because of the (relatively large) coefficient of proportionality, for
small $M\/$ this additional cost can be comparable to the cost of inverting the
Poisson problem. Obviously, this extra cost becomes less significant as
$M\/$ increases.
%
\subsection{Interface Representation.}\label{sub:interface}
As far as the CFM is concerned, the framework needed to solve the local Cauchy
problems is entirely described above. However, there is an important issue that
deserves attention: the representation of the interface. This question is
independent of the CFM. Many approaches are possible, and the optimal choice is
geometry dependent. The discussion below is meant to shed some light on this
issue, and motivate the solution we have adopted.

In the present work, generally we proceed assuming that the interface is not
known exactly -- since this is what frequently happens. The only exceptions to
this are the examples in \S~\ref{sub:case4} and \S~\ref{sub:case5},
which involve two distinct (circular) interfaces touching at a
point. In the generic setting, in addition to a proper representation of the
interfaces, one needs to be able to identify the distinct interfaces, regions
in between, contact points, as well as distinguish between a single interface
crossing the same stencil multiple times and multiple distinct interfaces
crossing one stencil. While the CFM algorithm is capable of dealing with these
situations once they have been identified (\eg\/~see remarks~\ref{rem:5:3} and
\ref{rem:5:4}), the development of an algorithm with the capability to detect
such generic geometries is  beyond the scope of this paper, and a (hard) problem in
interface representation. For these reasons, in the examples in \S~\ref{sub:case4} and
\S~\ref{sub:case5} we use an explicit exact representation of the interface.

To guarantee the accuracy of the solution for $D\/$, the interface conditions
must be applied with the appropriate accuracy --- see
\S~\ref{sub:erroranalysis}. Since these conditions are imposed on the interface
$\Gamma\/$, the location of $\Gamma\/$ must be known with the same order of
accuracy desired for $D\/$. In the particular case of the \fth\/ order
implementation of the CFM algorithm in this paper, this means \fth\/ order
accuracy. For this reason, we adopted the gradient-augmented level set (GA-LS)
approach, as introduced in \cite{nave:10}. This method allows a simple and
completely local \fth\/ order accurate representation of the interface, using
Hermite cubics defined everywhere in the domain. The approach also allows the
computation of normal vectors in a straightforward and accurate fashion.

We point out that the GA-LS method is not the only option for an implicit
\fth\/ order representation of the interface. For example, a regular level set
method~\cite{osher:88}, combined with a high-order interpolation scheme, could
be used as well. Here we adopted the GA-LS approach because of the algorithmic
coherence that results from representing both the level set, and the correction
functions, using the same bicubic polynomial base.

\subsection{Error analysis.} \label{sub:erroranalysis}
A naive reading of the discretized system in Eq.~\eqref{eq:9p-C} suggests that,
in order to obtain a fourth order accurate solution $u\/$, we need to compute
the CFM correction terms $C_{i,j}$ with fourth order accuracy. Thus, from
Eq.~\eqref{eq:9p-D}, it would follow that we need to know the correction
function $D\/$ with sixth order accuracy! This is, however, \emph{not correct},
as explained below.

Since we need to compute the correction function $D\/$ only at grid-points an
$\mathcal{O}(h)\/$ distance away from $\Gamma\/$, it should be clear that
errors in the $D_{i,j}$ are equivalent to errors in $a\/$ and $b\/$ of the
same order. But errors in $a\/$ and $b\/$ produce errors of the same order in
$u\/$ --- see Eq.~\eqref{eq:poisson} and Eq.~\eqref{eq:jump}. Hence, if we
desire a fourth order accurate solution $u\/$, we need to compute the
correction terms $D_{i,j}$ with fourth order accuracy only. This argument is
confirmed by the convergence plots in figures~\ref{fig:case1-convergence},
\ref{fig:case2-convergence}, and \ref{fig:case3-convergence}.

\subsection{Computation of gradients.}\label{sub:gradient}
Some applications require not only the solution to the Poisson problem, but
also its gradient. Hence, in \S~\ref{sec:results}, we
include plots characterizing the behavior of the errors in the gradients
of the solutions. A key question is then: how are these gradients computed?

To compute the gradients near the interface, the correction function can be
used to extend the solution across the interface, so that a standard stencil
can be used. However, for this to work, it is important to discretize the
gradient operator using the same nodes that are part of the 9-point stencil
--- so that the same correction functions obtained while solving the Poisson
equation can be used. Hence we discretize the gradient operator with a
procedure similar to the one used to obtain the 9-point stencil. Specifically,
we use the following \fth\/ order accurate discretization:
\begin{align}
 \partial_x u_{i,j} &= \hat{\partial}_{x} u_{i,j} 
                     + \dfrac{h_x^2}{6}\left[ \hat{\partial}_{xx}\hat{\partial}_y u_{i,j} 
                                            - (f_x)_{i,j}\right],\\
 \partial_y u_{i,j} &= \hat{\partial}_{y} u_{i,j} 
                     + \dfrac{h_y^2}{6}\left[ \hat{\partial}_{yy}\hat{\partial}_y u_{i,j} 
                                            - (f_y)_{i,j}\right],
\end{align}
where
\begin{align}
 \hat{\partial}_x u_{i,j} &= \dfrac{u_{i+1,j} - u_{i-1,j}}{2h_x},\label{eq:partialx}\\
 \hat{\partial}_y u_{i,j} &= \dfrac{u_{i,j+1} - u_{i,j-1}}{2h_y},\label{eq:partialy}
\end{align}
and $\hat{\partial}_{xx}$ and $\hat{\partial}_{yy}$ are defined by
\eqref{eq:partialxx} and \eqref{eq:partialyy}, respectively. The terms
$(f_x)_{i,j}\/$ and  $(f_y)_{i,j}\/$ may be given analytically (if known),
or computed using appropriate second order accurate discretizations.

This discretization is \fth\/ order accurate. However, since the error in
the correction function is (generally) not smooth, the resulting gradient
will be less than \fth\/ order accurate (worse case scenario is \trd\/ order
accurate) next to the interface.

\section{Results.} \label{sec:results}
%
\subsection{General Comments.}
In this section we present five examples of computations in 2D using the algorithm introduced in \S~\ref{sec:scheme}. We solve the Poisson problem in the unit square $[0,1] \times [0,1]\/$ for five different configurations. Each example is defined below in terms of the problem parameters (source term $f\/$, and jump conditions across $\Gamma\/$), the representation of the interface(s) --- either exact or using a level set function, and the exact solution (needed to evaluate the errors in the convergence plots). Notice that
\begin{enumerate}
\item As explained in \S~\ref{sub:interface}, in examples 1 through 3 we represent the interface(s) using the GA-LS method. Hence, the interface is defined by a level set function $\phi\/$, with gradient $\vec{\nabla}\phi = \p{\phi_x,\phi_y}\/$ --- both of which are carried within the GA-LS framework~\cite{nave:10}.
\item Below the level set is described via an analytic formula. In examples 1 through 3 this formula is converted into the GA-LS representation for the level set, before it is fed into the code. Only this representation is used for the actual computations. This is done so as to test the code's performance under generic conditions -- where the interface $\Gamma\/$ would be known via a level set representation only.
\item Within the GA-LS framework we can, easily and accurately, compute the vectors normal to the interface --- anywhere in the domain. Hence, it is convenient to write the jump in the normal derivative, $\left[u_n\right]_{\Gamma}\/$, in terms of the jump in the gradient of $u\/$ dotted with the normal to the interface $\hat{n} = (n_x\/,\,n_y)\/$.
\end{enumerate}
The last two examples involve touching circular interfaces and were devised to demonstrate the robustness of the CFM in the presence of interfaces that are very close together. In these last two examples, for the reasons discussed in \S~\ref{sub:interface}, we decided to use an exact representation of the circular interfaces.
%
\subsection{Example 1.}\label{sub:case1}
\begin{itemize}
 \item Problem parameters:
  \begin{align*}
   f^+\p{x\/,\,y} & = - 2\/\pi^2\sin(\pi\/x)\sin(\pi\/y)\/,\\
   f^-\p{x\/,\,y} & = - 2\/\pi^2\sin(\pi\/x)\sin(\pi\/y)\/,
  \end{align*}
  \begin{align*}
   [u]_{\Gamma}              & = \sin(\pi\/x)\exp(\pi\/y)\/,\\
   \left[u_n\right]_{\Gamma} & = \pi\,\left[\cos(\pi\/x)\exp(\pi\/y)\,n_x
                              + \sin(\pi\/x)\exp(\pi\/y)\,n_y\right]\/.
  \end{align*}
 %
 \item 
 Level set defining the interface: \hfill
 $\phi\p{x\/,\,y}  = (x-x_0)^2 + (y-y_0)^2 - r_0^2\/$,
  where $x_0 = 0.5\/$, $y_0 = 0.5\/$, and $r_0 = 0.1\/$.
 %
 \item Exact solution:
  \begin{align*}
   u^+\p{x\/,\,y} & = \sin(\pi\/x)\sin(\pi\/y)\/,\\
   u^-\p{x\/,\,y} & = \sin(\pi\/x)[\sin(\pi\/y) - \exp(\pi\/y)]\/.
  \end{align*}
\end{itemize}

\begin{figure}[htb!]
 \begin{center}
  \includegraphics[width=3.5in]{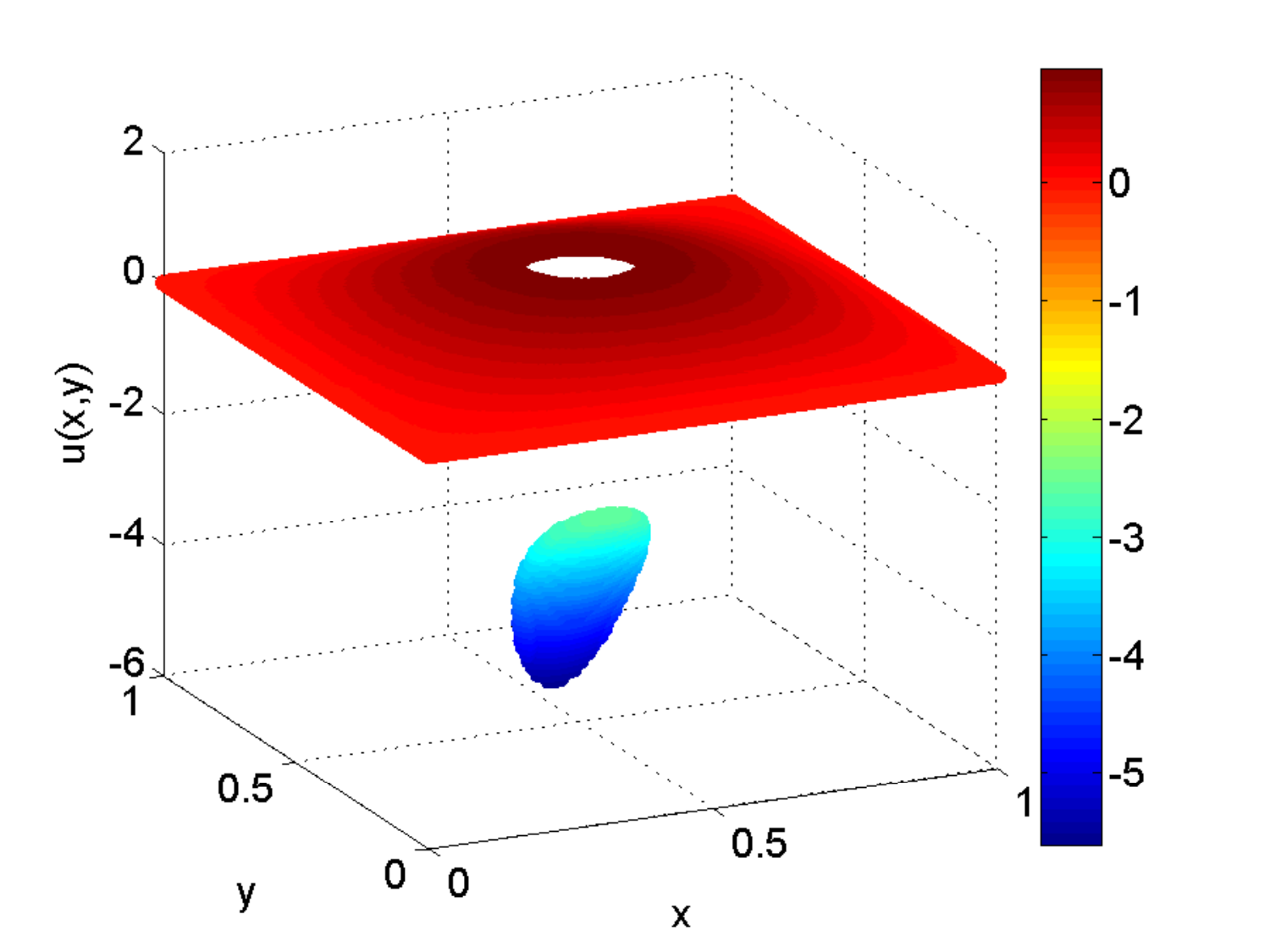}
 \end{center}
 \caption{Example 1 - numerical solution with $193 \times 193\/$ nodes.}
 \label{fig:case1-numerical}
\end{figure}

\begin{figure}[htb!]
 \begin{center}
  \includegraphics[width=2.6in]{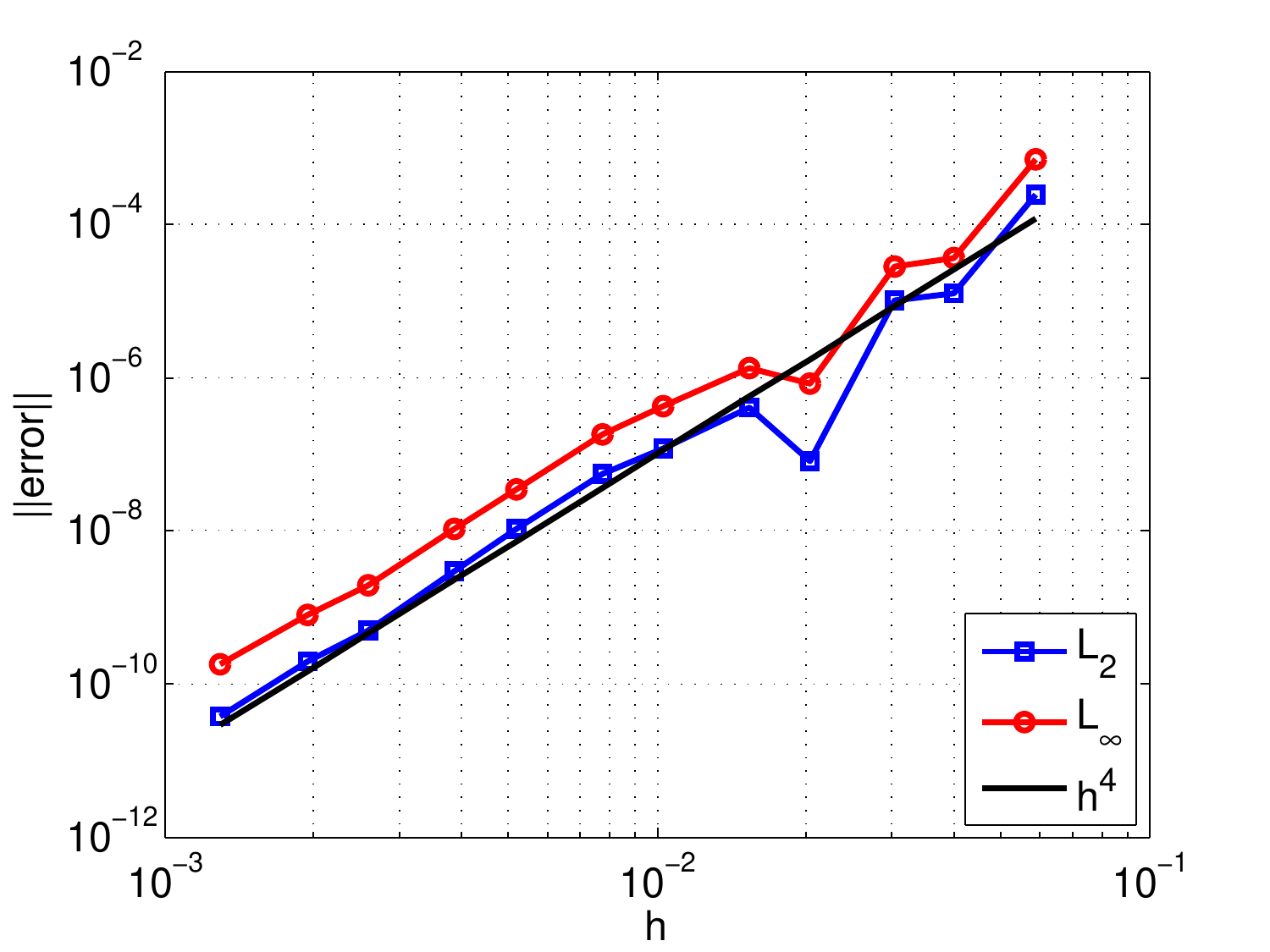}
  \includegraphics[width=2.6in]{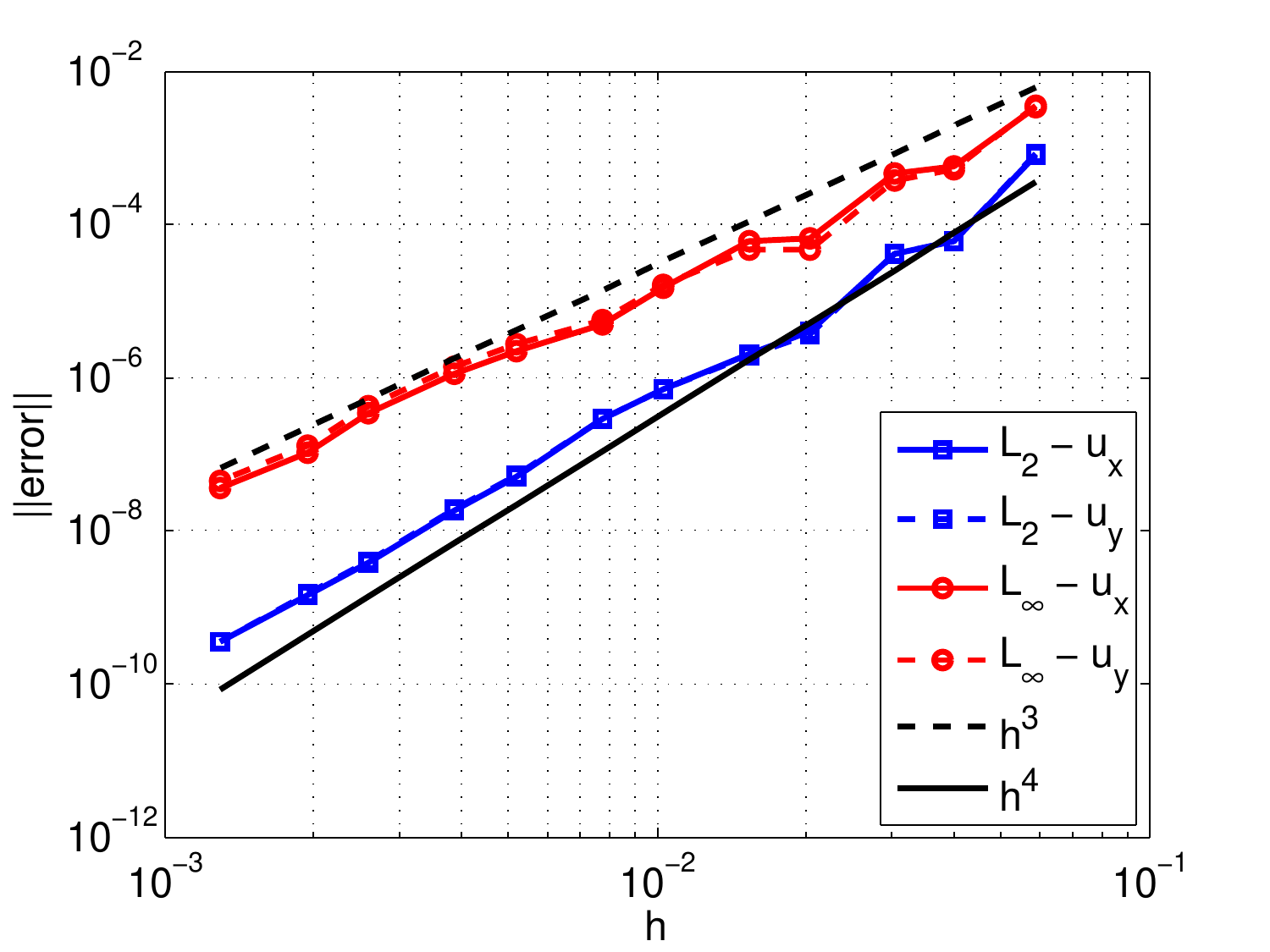}\\
  (a) Solution. \hspace*{1.6in} (b) Gradient.\hfill
 \end{center}
 \caption{Example 1 - Error behavior of the solution and its gradient in the $L_2\/$ and $L_\infty\/$ norms.}
 \label{fig:case1-convergence}
\end{figure}

Figure~\ref{fig:case1-numerical} shows the numerical solution with a fine grid
($193 \times 193$ nodes). The discontinuity is captured very sharply, and it
causes no oscillations in the solution. In addition, figure~\ref{fig:case1-convergence}
shows the behavior of the error of the solution and its gradient in the $L_2$ and
$L_{\infty}$ norms. As expected, the solution presents \fth\/ order convergence 
as the grid is refined. Moreover, the gradient converges to \trd\/ order in the
$L_{\infty}$ norm and to \fth\/ order in the $L_2$ norm, which is a reflection of
the fact that the error in the solution is not smooth in a narrow region close to
the interface only.

%
\subsection{Example 2.}\label{sub:case2}
\begin{itemize}
 \item Problem parameters:
  \begin{align*}
   f^+\p{x\/,\,y} & = 0\/,\\
   f^-\p{x\/,\,y} & = 0\/,
  \end{align*}
  \begin{align*}
   [u]_{\Gamma}              & = -\exp(x)\cos(y)\/,\\
   \left[u_n\right]_{\Gamma} & = -\exp(x)\cos(y)\,n_x + \exp(x)\sin(y)\,n_y\/.
  \end{align*}
 %
\item 
 Level set defining the interface: \hfill
 $\phi\p{x\/,\,y}  = (x-x_0)^2 + (y-x_0)^2 - r^2(\theta)\/$,
 where $r(\theta) = r_0 + \epsilon\,\sin(5\,\theta)\/$,
       $\theta\p{x\/,\,y} = \arctan\p{\dfrac{y-y_0}{x-x_0}}\/$,
       $x_0 = 0.5\/$, $y_0 = 0.5\/$, $r = 0.25\/$, and $\epsilon = 0.05\/$.
 \item Exact solution:
  \begin{align*}
   u^+\p{x\/,\,y} & = 0,\\
   u^-\p{x\/,\,y} & = \exp(x)\cos(y)\/.
  \end{align*}
\end{itemize}

\begin{figure}[htb!]
 \begin{center}
  \includegraphics[width=3.5in]{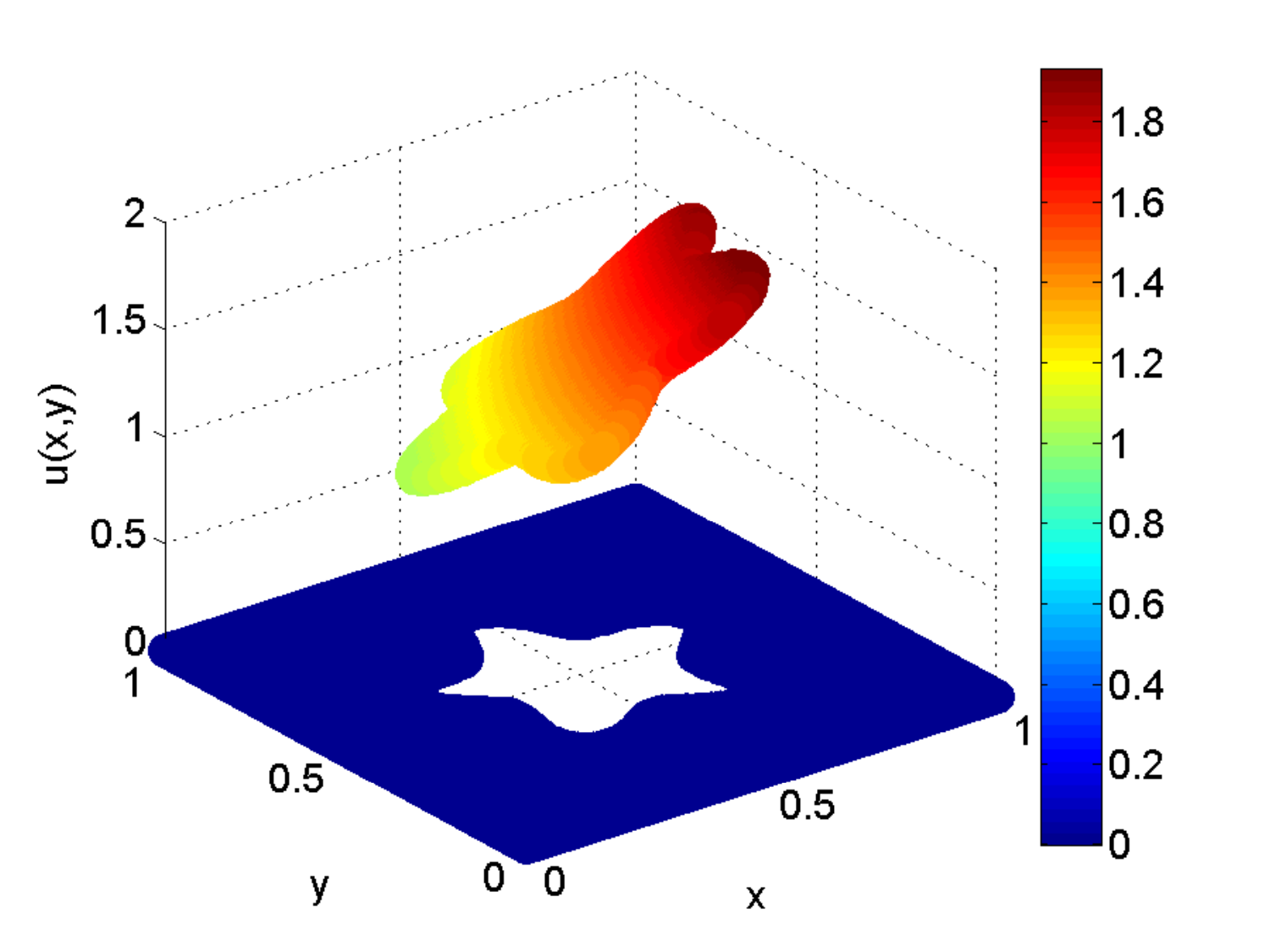}
 \end{center}
 \caption{Example 2 - numerical solution with $193 \times 193\/$ nodes.}
 \label{fig:case2-numerical}
\end{figure}

\begin{figure}[htb!]
 \begin{center}
  \includegraphics[width=2.6in]{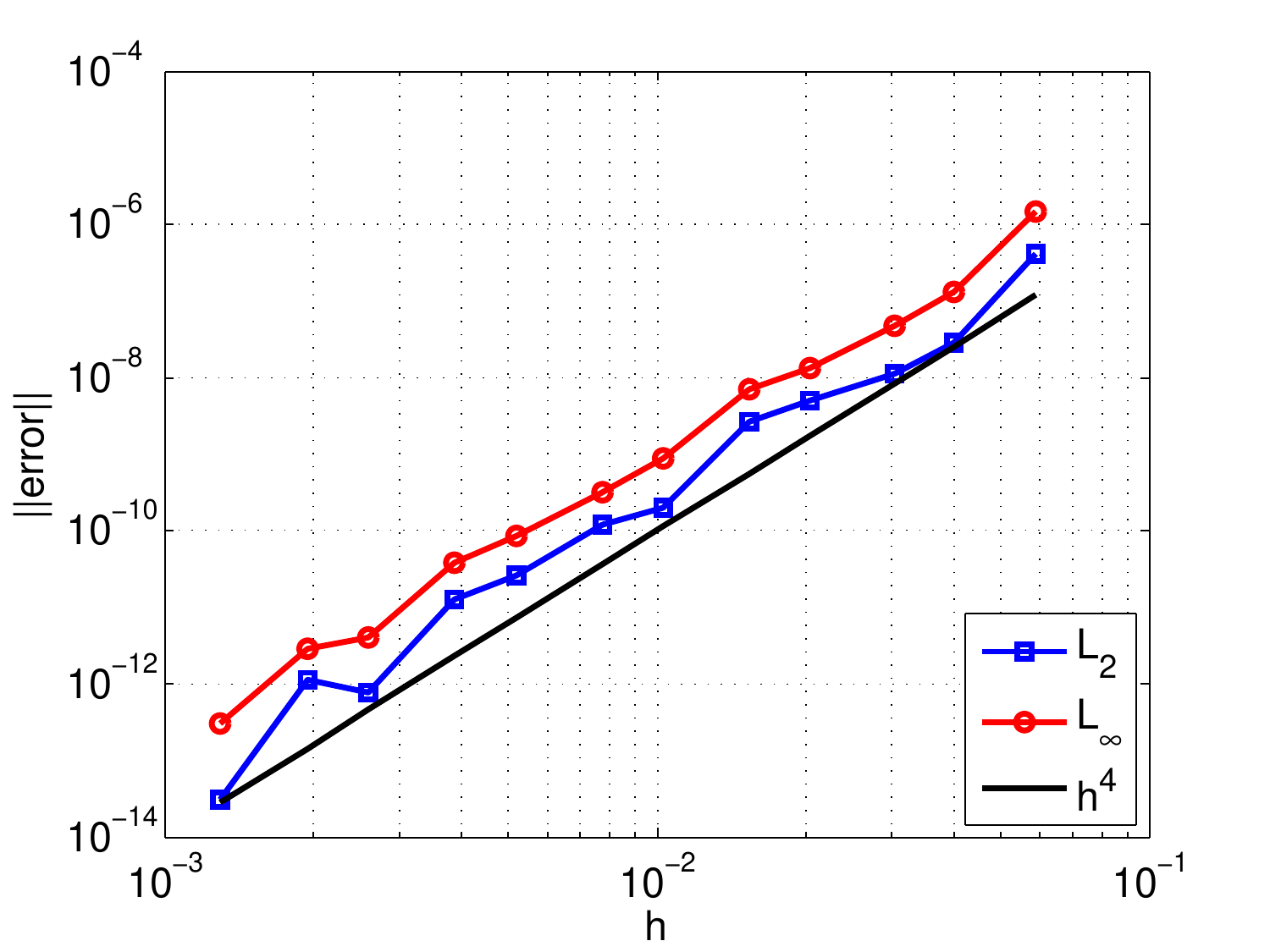}
  \includegraphics[width=2.6in]{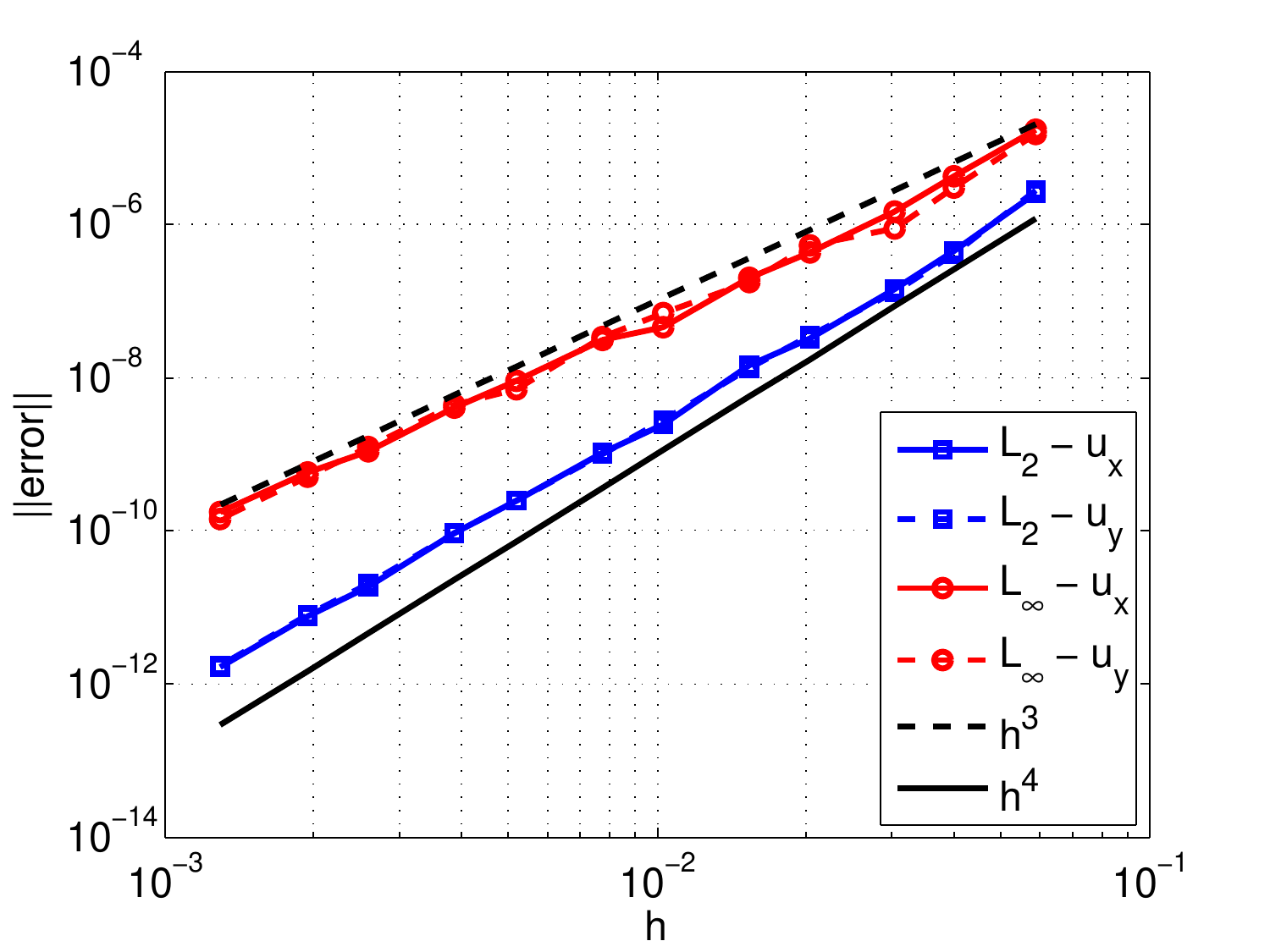}\\
  (a) Solution. \hspace*{1.6in} (b) Gradient.\hfill
 \end{center}
 \caption{Example 2 - Error behavior of the solution and its gradient in the $L_2\/$ and $L_\infty\/$ norms.}
 \label{fig:case2-convergence}
\end{figure}

Figure~\ref{fig:case2-numerical} shows the numerical solution with a fine grid
($193 \times 193$ nodes). Once again, the overall quality of the solution is
very satisfactory. Figure~\ref{fig:case2-convergence} shows the behavior of the
error of the solution and its gradient in the $L_2$ and $L_{\infty}$ norms.
Again, the solution converges to \fth\/ order, while the gradient converges
to \trd\/ order in the $L_{\infty}$ norm and close to \fth\/ order in the $L_2$ norm.
However, unlike what happens in example 1, small wiggles are observed in the error
plots. This behavior can be explained in terms of the construction of the sets
$\Omega_{\Gamma}^{i,j}$ --- see \S~\ref{sec:scheme}. The approach used to
construct $\Omega_{\Gamma}^{i,j}$ is highly dependent on the way in which the
grid points are placed relative to the interface. Thus, as the grid is refined,
the arrangement of the $\Omega_{\Gamma}^{i,j}$ can vary quite a lot --- specially
for a ``complicated'' interface such as the one in this example. What this means
is that, while one can guarantee that the correction function $D\/$ is obtained
with \fth\/ order precision, the proportionality coefficient is not
constant --- it may vary a little from grid to grid. This variation is
responsible for the small oscillations observed in the convergence plot.
Nevertheless, despite these oscillations, the overall convergence is clearly
\fth\/ order.
%

\subsection{Example 3.}\label{sub:case3}
\begin{itemize}
 \item Problem parameters:
  \begin{align*}
   f^+\p{x\/,\,y} & = \exp(x)\left[2 + y^2 + 2\/\sin(y)
                    + 4\/x\sin(y)\right]\/,\\
   f^-\p{x\/,\,y} & = 40\/,
  \end{align*}
  \begin{align*}
   [u]_{\Gamma}             & = \exp(x)\left[x^2\sin(y) + y^2\right]
                              - 10\/\p{x^2+y^2}\/,\\
   \left[u_n\right]_{\Gamma} & = \left\{\exp(x)\left[\left(x^2+2\/x\right)
                                \sin(y) + y^2\right] - 20\/x\right\}\,n_x\\
                            & + \left\{\exp(x)\left[x^2\cos(y) + 2\/y\right]
                              - 20\/y\right\}\,n_y\/.
  \end{align*}
 %
 \item 
 Level set defining the interface: \\
 $\phi\p{x\/,\,y} = \left[(x-x_1)^2 + (y-y_1)^2 - r_1^2\right]
                    \left[(x-x_2)^2 + (y-y_2)^2 - r_2^2\right]\/$,
  where $x_1 = y_1 = 0.25\/$, $r_1 = 0.15\/$,
        $x_2 = y_2 = 0.75\/$, and $r_2 = 0.1\/$.
 \item Exact solution:
  \begin{align*}
   u^+\p{x\/,\,y} & = \exp(x)\left[x^2\sin(y) + y^2\right]\/,\\
   u^-\p{x\/,\,y} & = 10\p{x^2+y^2}\/.
  \end{align*}
\end{itemize}

\begin{figure}[htb!]
 \begin{center}
  \includegraphics[width=3.5in]{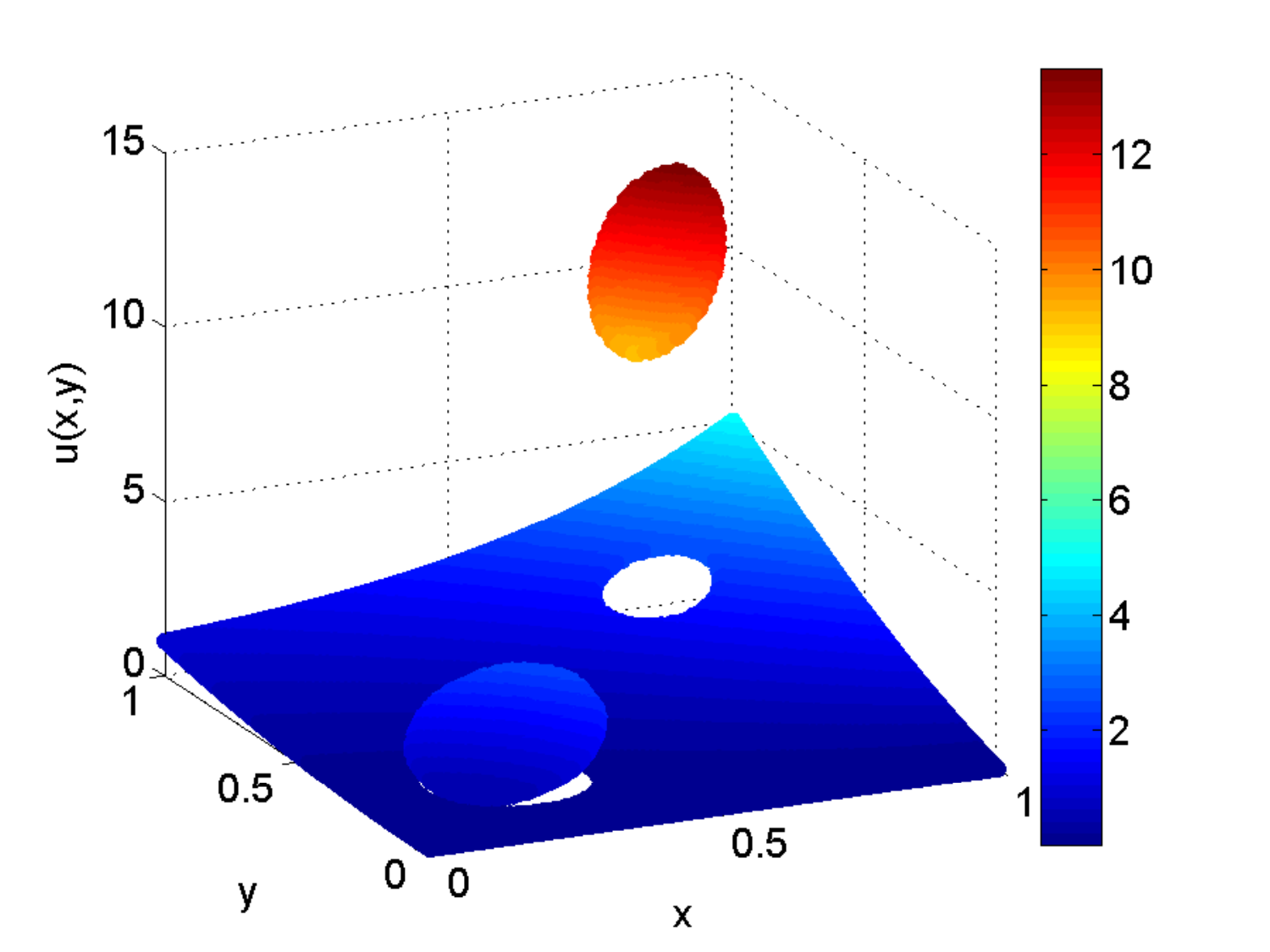}
 \end{center}
 \caption{Example 3 - numerical solution with $193 \times 193\/$ nodes.}
 \label{fig:case3-numerical}
\end{figure}

\begin{figure}[htb!]
 \begin{center}
  \includegraphics[width=2.6in]{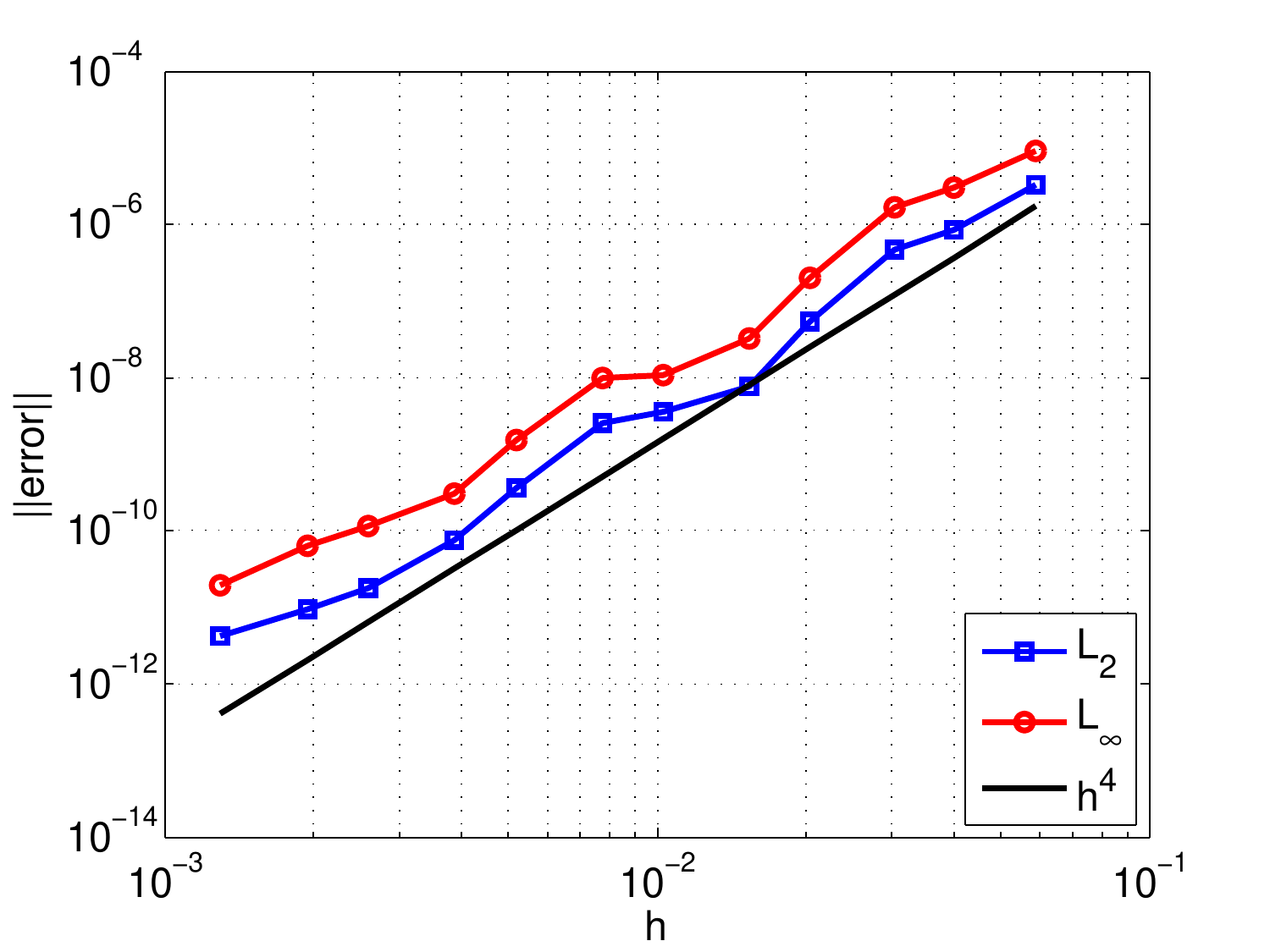}
  \includegraphics[width=2.6in]{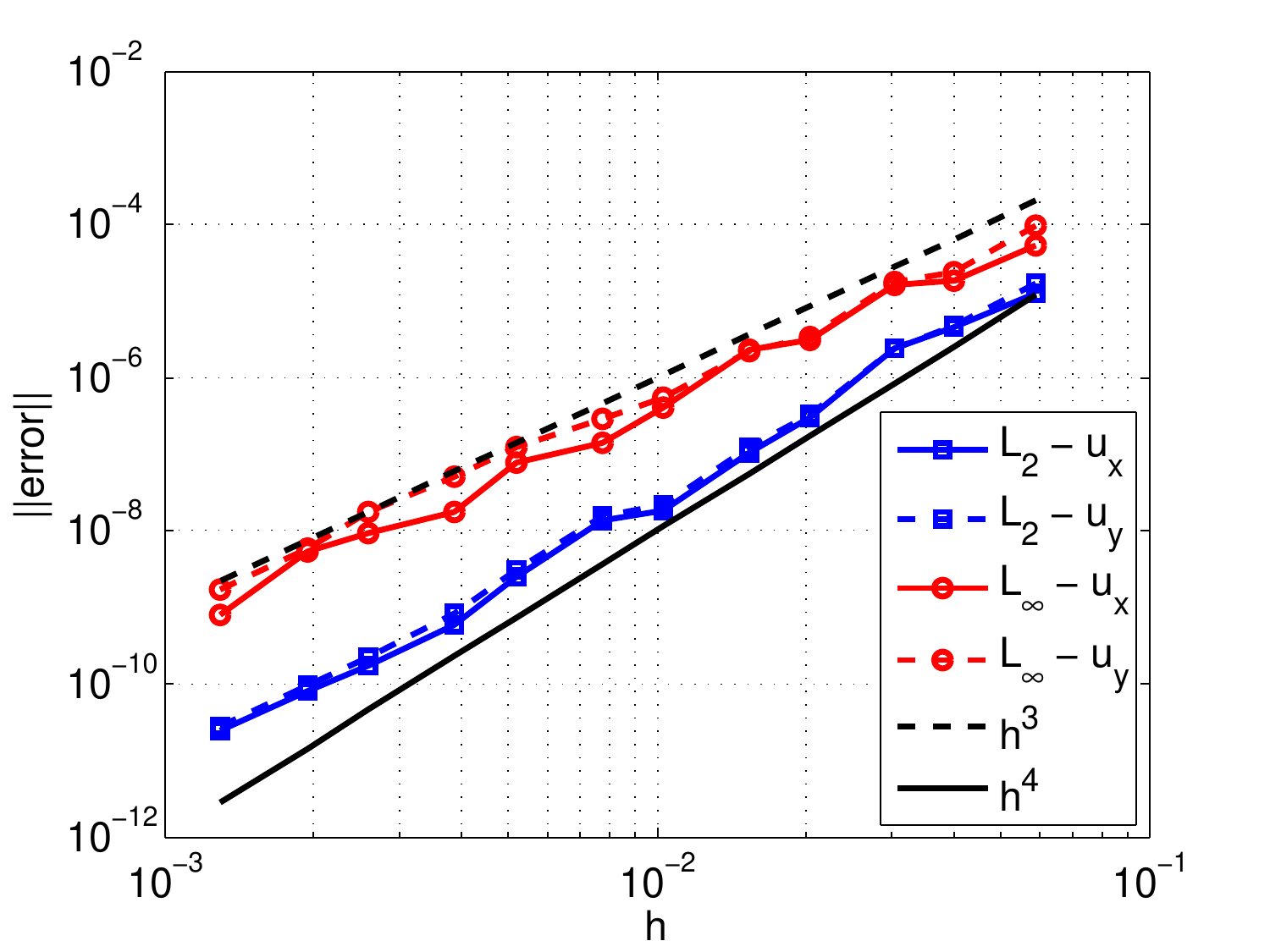}\\
  (a) Solution. \hspace*{1.6in} (b) Gradient.\hfill
 \end{center}
 \caption{Example 3 - Error behavior of the solution and its gradient in the $L_2\/$ and $L_\infty\/$ norms.}
 \label{fig:case3-convergence}
\end{figure}

Figure~\ref{fig:case3-numerical} shows the numerical solution with a fine grid
($193 \times 193\/$ nodes). In this example, there are two circular interfaces in
the solution domain. The two regions inside the circles make $\Omega^-\/$,
while the remainder of the domain is $\Omega^+\/$. This example shows that the
method is general enough to deal with multiple interfaces, keeping the same
quality in the solution. Figure~\ref{fig:case3-convergence} shows that the solution
converges to \fth\/ order in both $L_{\infty}$ and $L_2$ norms, while the gradient
converges to \trd\/ order in the $L_{\infty}$ norm and close to \fth\/ order in the $L_2$
norm.

\subsection{Example 4.}\label{sub:case4}
\begin{itemize}
 \item Problem parameters:
  \begin{align*}
   f_1\p{x\/,\,y} & = - 2\/\pi^2\sin(\pi\/x)\sin(\pi\/y)\/,\\
   f_2\p{x\/,\,y} & = \exp(x)\left[2 + y^2 + 2\/\sin(y)
                    + 4\/x\sin(y)\right]\/,\\
   f_3\p{x\/,\,y} & = - 2\/\pi^2\sin(\pi\/x)\sin(\pi\/y)\/,
  \end{align*}
  \begin{align*}
   [u]_{\Gamma_{1-2}}       & = \exp(x)\left[x^2\sin(y) + y^2\right] - \sin(\pi\/x)\sin(\pi\/y) - 5\/,\\
   \left[u_n\right]_{\Gamma_{1-2}} & = \{\exp(x)\left[(x^2 + 2x)\sin(y) + y^2\right] - \pi\cos(\pi\/x)\sin(\pi\/y)\}n_x\\
                                   & + \{\exp(x)\left[x^2\cos(y) + 2y\right] - \pi\sin(\pi\/x)\cos(\pi\/y)\}n_y\/,
  \end{align*}
  \begin{align*}
   [u]_{\Gamma_{2-3}}       & = \sin(\pi\/x)[\sin(\pi\/y) - \exp(\pi\/y)] - \exp(x)\left[x^2\sin(y) + y^2\right]\/,\\
   \left[u_n\right]_{\Gamma_{2-3}} & = \{\pi\cos(\pi\/x)[\sin(\pi\/y) - \exp(\pi\/y)] - \exp(x)\left[(x^2 + 2x)\sin(y) + y^2\right]\}n_x\\
                                   & + \{\pi\sin(\pi\/x)[\cos(\pi\/y) - \exp(\pi\/y)] - \exp(x)\left[x^2\cos(y) + 2y\right]\}n_y\/.
  \end{align*}
 %
 \item Interface (exact representation):
  \begin{itemize}
   \item Region 1: inside of the big circle.
   \item Region 2: outer region.
   \item Region 3: inside of the small circle.
  \end{itemize}
  \begin{itemize}
   \item Interface 1--2 (Big circle):
    \begin{align*}
     r_B &= 0.3,\\
     x_{0_B} &= 0.5,\\
     y_{0_B}  &= 0.5,
    \end{align*}
 
   \item Interface 2--3 (Small circle):
    \begin{align*}
     r_S &= 0.3,\\
     x_{0_S} &= x_{0_B} + r_B\cos(\pi\//e^2) - r_S\cos\p{\pi\/(1/e^2 + 1)},\\
     y_{0_S} &= y_{0_B} + r_B\sin(\pi\//e^2) - r_S\sin\p{\pi\/(1/e^2 + 1)}.
    \end{align*}
  \end{itemize}
 %
 \item Exact solution:
  \begin{align*}
   u_1\p{x\/,\,y} & = \sin(\pi\/x)\sin(\pi\/y) + 5\/,\\
   u_2\p{x\/,\,y} & = \exp(x)\left[x^2\sin(y) + y^2\right]\/,\\
   u_3\p{x\/,\,y} & = \sin(\pi\/x)[\sin(\pi\/y) - \exp(\pi\/y)]\/.
  \end{align*}
\end{itemize}

\begin{figure}[htb!]
 \begin{center}
  \includegraphics[width=3.5in]{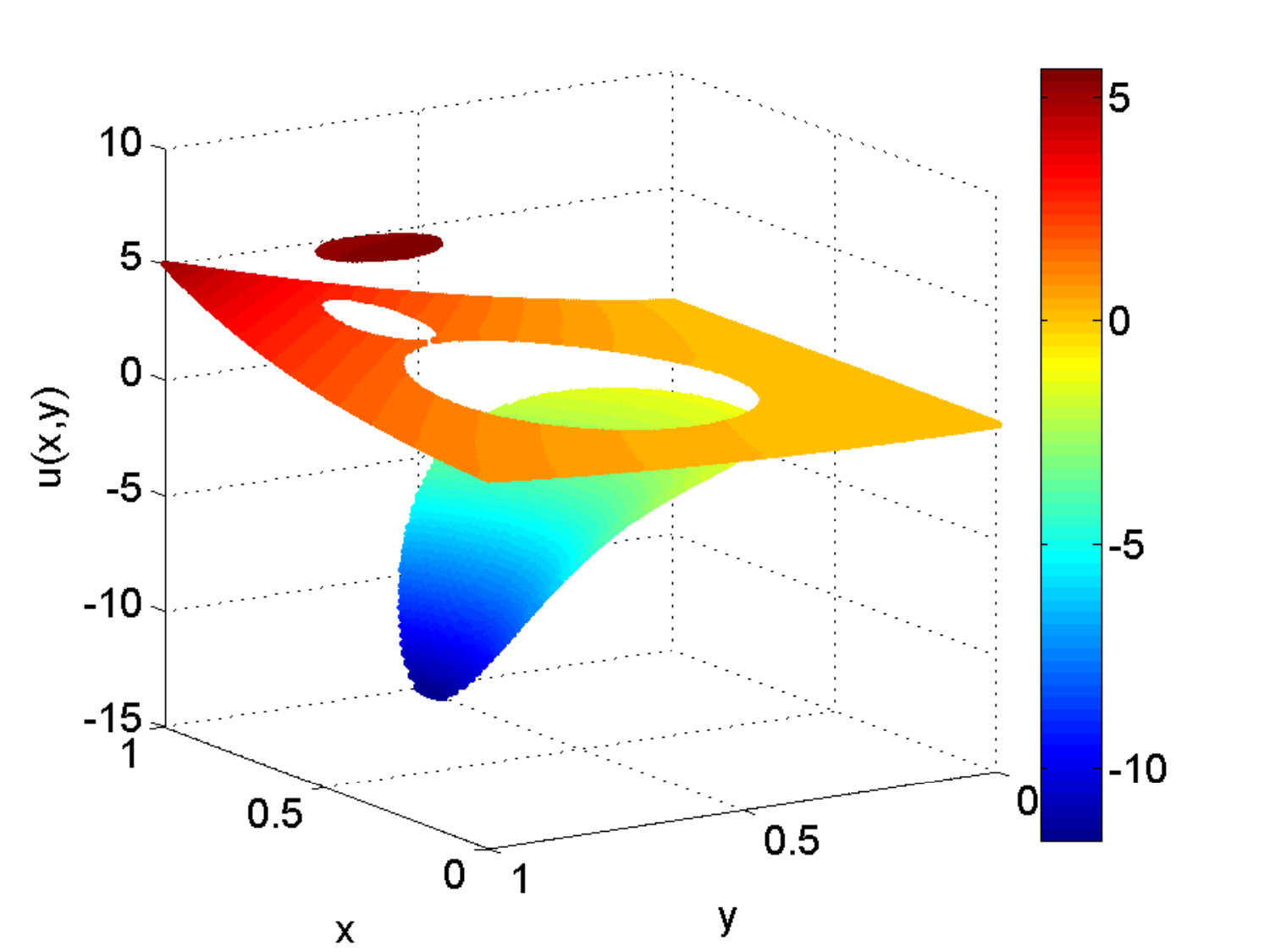}
 \end{center}
 \caption{Example 4 - numerical solution with $193 \times 193\/$ nodes.}
 \label{fig:case4-numerical}
\end{figure}

\begin{figure}[htb!]
 \begin{center}
  \includegraphics[width=2.6in]{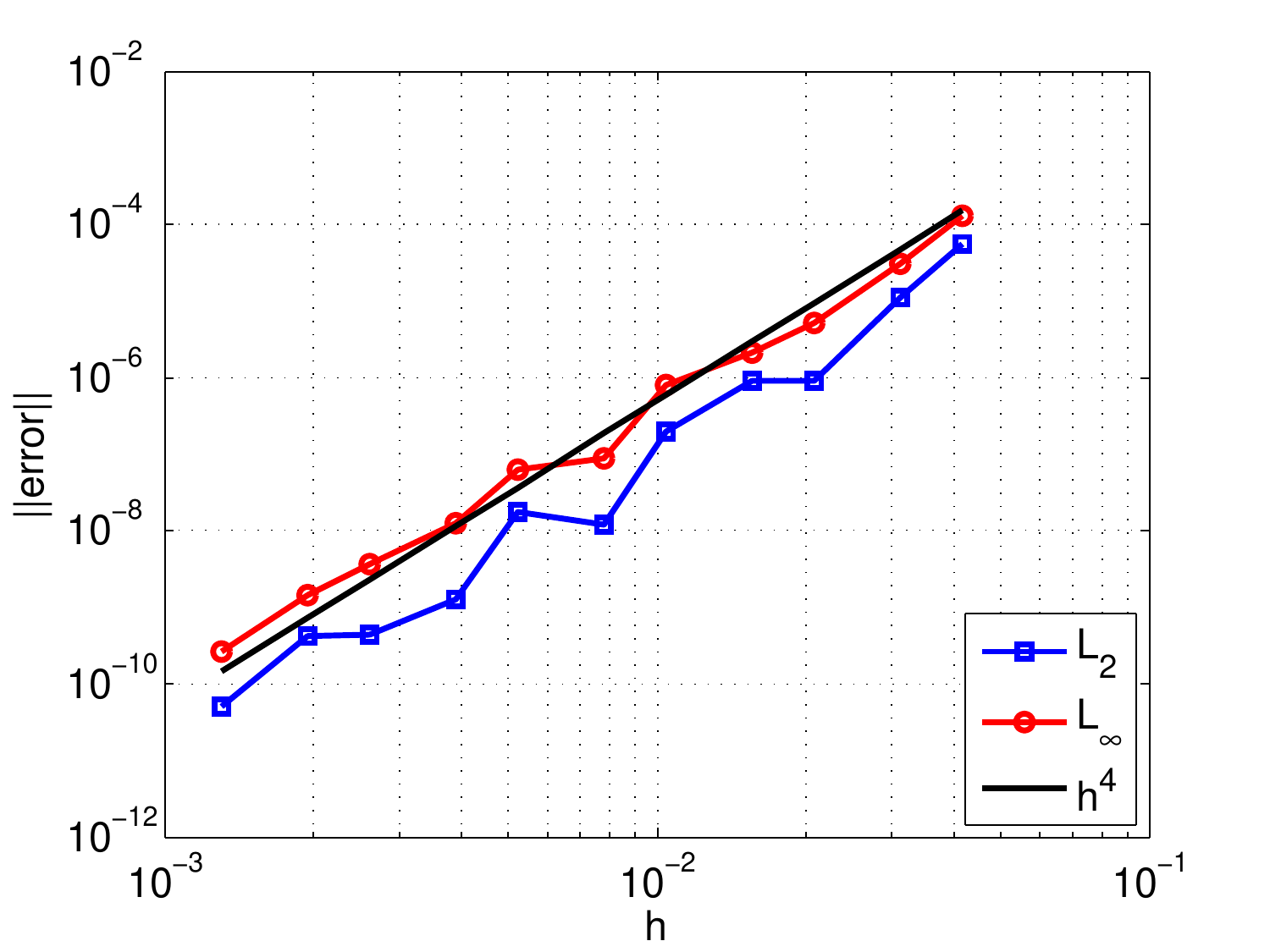}
  \includegraphics[width=2.6in]{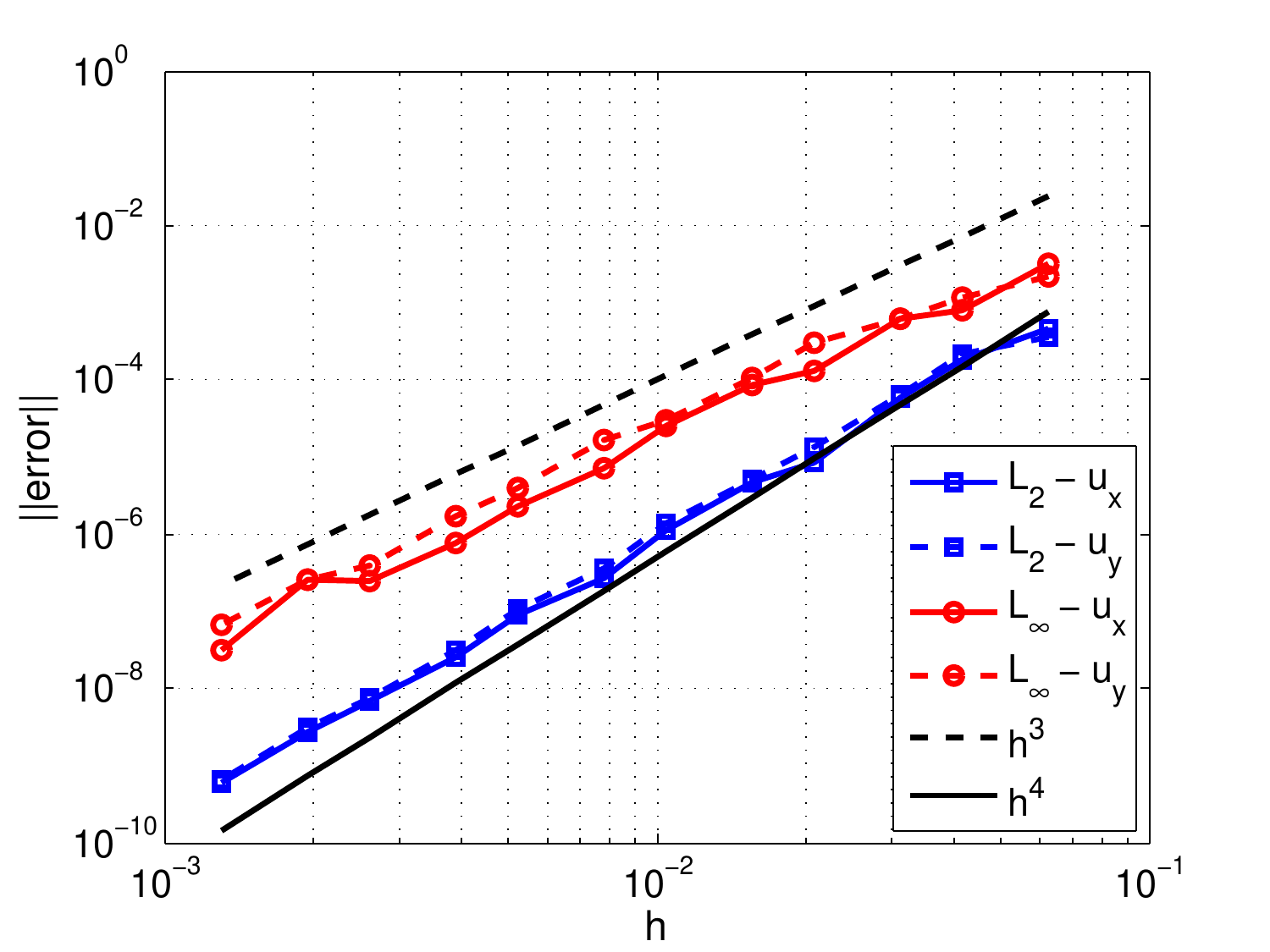}\\
  (a) Solution. \hspace*{1.6in} (b) Gradient.\hfill
 \end{center}
 \caption{Example 4 - Error behavior of the solution and its gradient in the $L_2\/$ and $L_\infty\/$ norms.}
 \label{fig:case4-convergence}
\end{figure}

Figure~\ref{fig:case4-numerical} shows the numerical solution with a fine grid
($193 \times 193\/$ nodes). In this example the big circle is centered within
the square integration domain and the small circle is external to it, with
a common point of tangency. The point of contact is placed along the boundary
of the big circle at the polar angle $\theta = \pi/e^2$ --- use the center of
the big circle as the polar coordinates' origin. These choices guarantee that,
as the grid is refined, a wide variety of configurations involving two distinct
interfaces crossing the same stencil occurs in a neighborhood of the contact
point. In particular, the selection of the angle $\theta\/$ is so that no
special alignments of the grid with the local geometry near the contact point
can happen. Figure~\ref{fig:case4-convergence} shows the behavior of the error
as $h \to 0\/$ in the $L_2$ and $L_{\infty}$ norms. Once again we observe \fth\/
order convergence (with small superimposed oscillations) for the solution. Moreover,
the gradient converges to \trd\/ order in the $L_{\infty}$ norm and close to \fth\/ order
in the $L_2$ norm. This example shows that the CFM is robust even in situations
where distinct interfaces can get arbitrarily close (tangent at a point).
%
%
%
\subsection{Example 5.}\label{sub:case5}
\begin{itemize}
 \item Problem parameters:
  \begin{align*}
   f_1\p{x\/,\,y} & = - 2\/\pi^2\sin(\pi\/x)\sin(\pi\/y)\/,\\
   f_2\p{x\/,\,y} & = - 2\/\pi^2\sin(\pi\/x)\sin(\pi\/y)\/,\\
   f_3\p{x\/,\,y} & = \exp(x)\left[2 + y^2 + 2\/\sin(y)
                    + 4\/x\sin(y)\right]\/,
  \end{align*}
  \begin{align*}
   [u]_{\Gamma_{1-2}}       & = -\left[\sin(\pi\/x)\exp(\pi\/y) + 5\right]\/,\\
   \left[u_n\right]_{\Gamma_{1-2}} & = -\pi\left[\cos(\pi\/x)\exp(\pi\/y)n_x + \sin(\pi\/x)\exp(\pi\/y)n_y\right]\/,
  \end{align*}
  \begin{align*}
   [u]_{\Gamma_{2-3}}       & = \exp(x)\left[x^2\sin(y) + y^2\right] - \sin(\pi\/x)[\sin(\pi\/y) - \exp(\pi\/y)]\/,\\
   \left[u_n\right]_{\Gamma_{2-3}} & = \{\exp(x)\left[(x^2 + 2x)\sin(y) + y^2\right] - \pi\cos(\pi\/x)[\sin(\pi\/y) - \exp(\pi\/y)]\}n_x\\
  & + \{\exp(x)\left[x^2\cos(y) + 2y\right] - \pi\sin(\pi\/x)[\cos(\pi\/y) - \exp(\pi\/y)]\}n_y.\\
  \end{align*}
 %
 \item Interface (exact representation):
  \begin{itemize}
   \item Region 1: inside small circle.
   \item Region 2: region between circles.
   \item Region 3: outer region.
  \end{itemize}
  \begin{itemize}
   \item Interface 2--3 (Big circle):
    \begin{align*}
     r_B &= 0.3,\\
     x_{0_B} &= 0.5,\\
     y_{0_B}  &= 0.5,
    \end{align*}
 
   \item Interface 1--2 (Small circle):
    \begin{align*}
     r_S &= 0.3,\\
     x_{0_S} &= x_{0_B} + (r_B - r_S)\cos(\pi\//e^2),\\
     y_{0_S} &= y_{0_B} + (r_B - r_S)\sin(\pi\//e^2).
    \end{align*}
  \end{itemize}
 %
 \item Exact solution:
  \begin{align*}
   u_1\p{x\/,\,y} & = \sin(\pi\/x)\sin(\pi\/y) + 5\/,\\
   u_2\p{x\/,\,y} & = \sin(\pi\/x)[\sin(\pi\/y) - \exp(\pi\/y)]\/,\\
   u_3\p{x\/,\,y} & = \exp(x)\left[x^2\sin(y) + y^2\right]\/.
  \end{align*}
\end{itemize}

\begin{figure}[htb!]
 \begin{center}
  \includegraphics[width=3.5in]{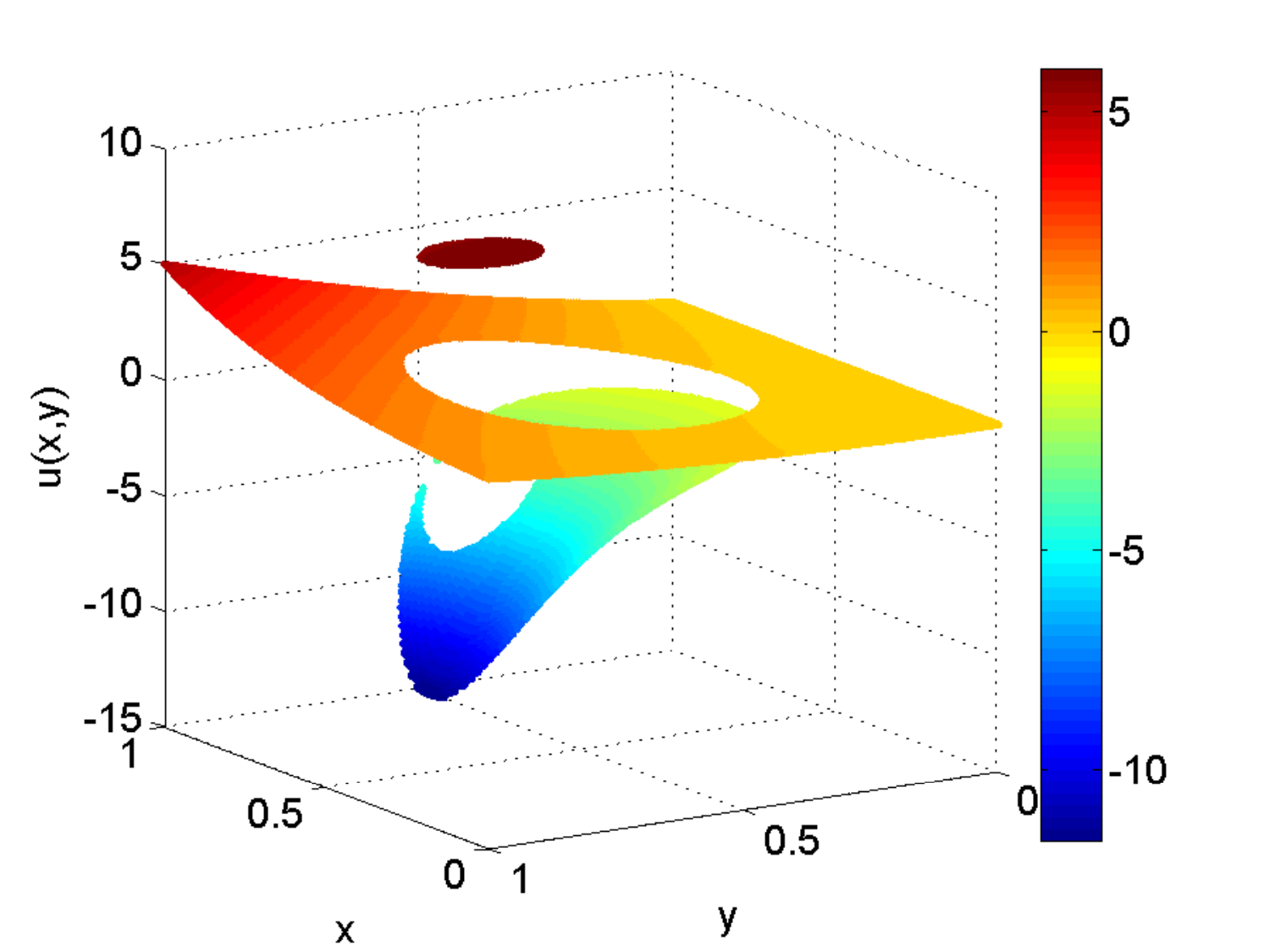}
 \end{center}
 \caption{Example 5 - numerical solution with $193 \times 193\/$ nodes.}
 \label{fig:case5-numerical}
\end{figure}

This example complements the example in~\S~\ref{sub:case4}, the sole
difference being that here the small circle is inside the big circle. The
results are entirely similar. Figure~\ref{fig:case5-numerical} shows the
numerical solution with a fine grid ($193 \times 193\/$ nodes), while
figure~\ref{fig:case4-convergence} shows the behavior of the error in the
$L_2$ and $L_{\infty}$ norms.
\clearpage
%
\begin{figure}[htb!]
 \begin{center}
  \includegraphics[width=2.6in]{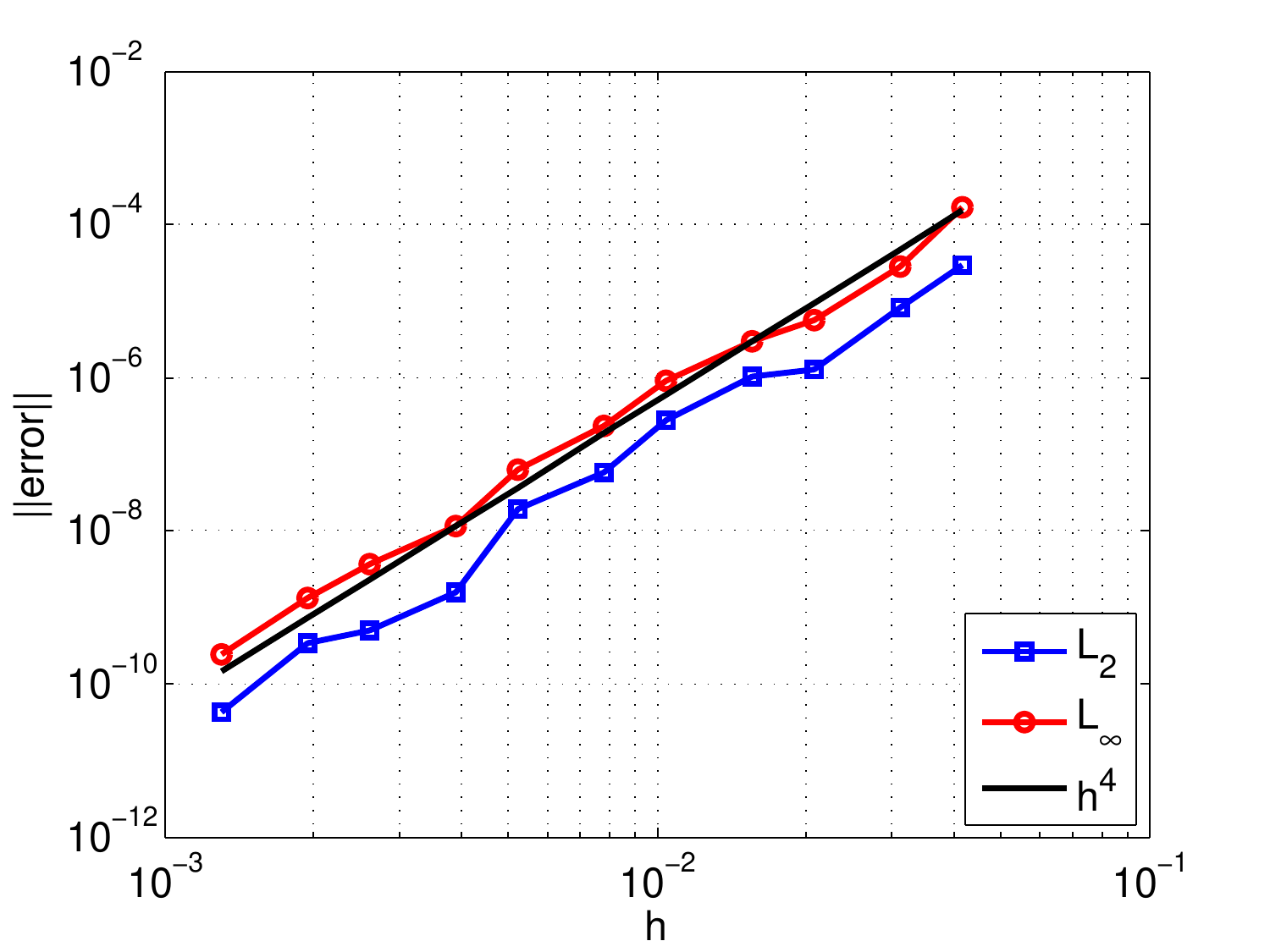}
  \includegraphics[width=2.6in]{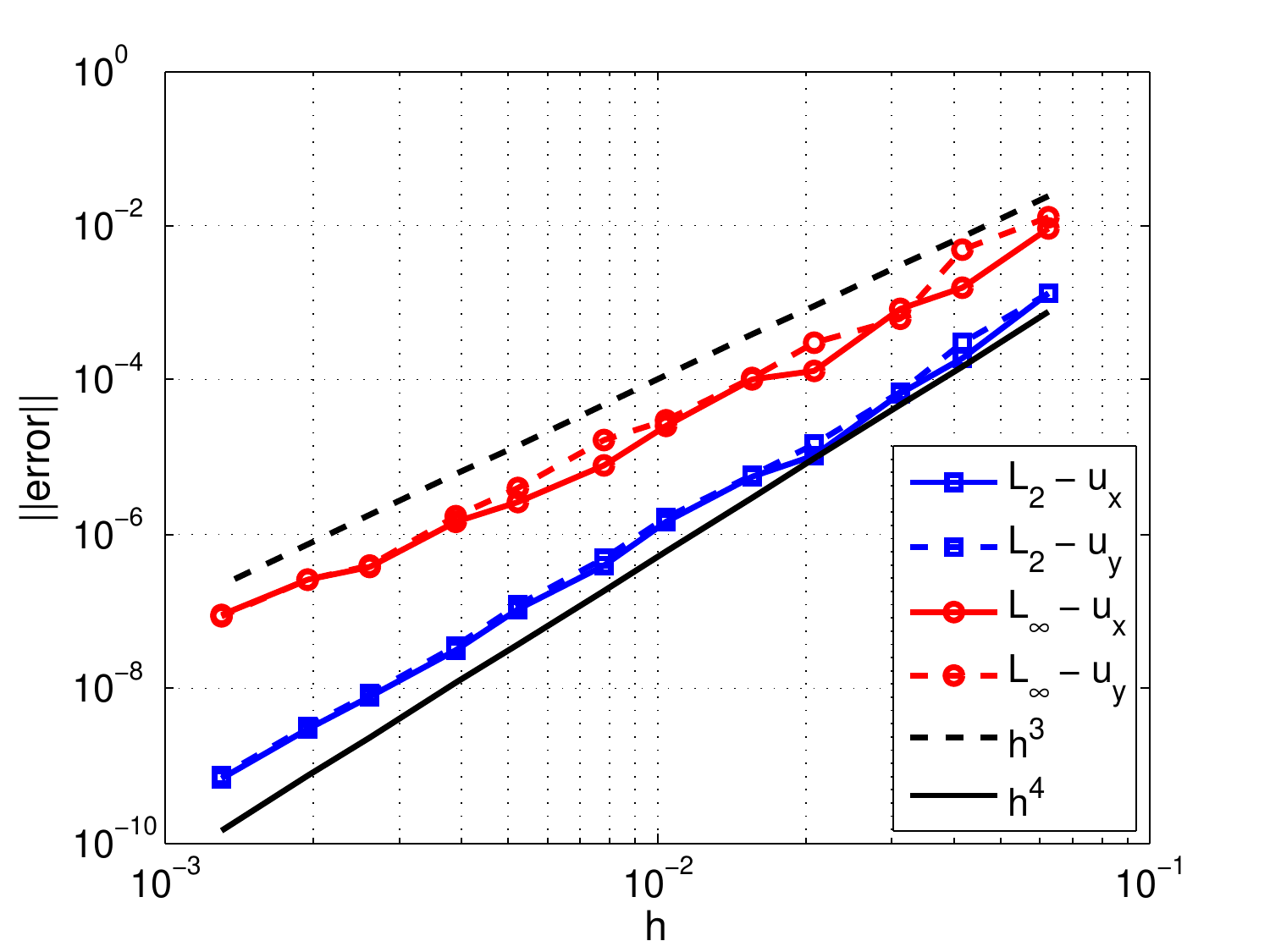}\\
  (a) Solution. \hspace*{1.6in} (b) Gradient.\hfill
 \end{center}
 \caption{Example 5 - Error behavior of the solution and its gradient in the $L_2\/$ and $L_\infty\/$ norms.}
 \label{fig:case5-convergence}
\end{figure}
\section{Conclusions.} \label{sec:conclusion}
In this paper we have introduced the Correction Function Method (CFM), which can,
in principle, be used to obtain (arbitrary) high-order accurate solutions to constant
coefficients Poisson problems with interface jump conditions. This method is based on
extending the correction terms idea of the Ghost Fluid Method (GFM) to that of a
correction function, defined in a (narrow) band enclosing the interface. This function
is the solution of a PDE problem, which can be solved (at least in principle) to any
desired order of accuracy. Furthermore, like the GFM, the CFM allows the use
of standard Poisson solvers. This feature follows from the fact that the interface jump
conditions modify (via the correction function) only the right-hand-side of the
discretized linear system of equations used in standard linear solvers for the Poisson
equation.

As an example application, the new method was used to create a \fth\/ order
accurate scheme to solve the constant coefficients Poisson equation with interface
jump conditions in 2D. In this scheme, the domain of definition of the correction
function is split into many grid size rectangular patches. In each patch the function
is represented in terms of a bicubic (with 12 free parameters), and the solution
is obtained by minimizing an appropriate (discretized) quadratic functional.
The correction function is thus pre-computed, and then is used to modify (in the standard way
of the GFM) the right hand side of the Poisson linear system, incorporating
the jump conditions into the Poisson solver. We used the standard \fth\/ order
accurate 9-point stencil discretization of the Laplace operator, to thus
obtain a \fth\/ order accurate method.

Examples were computed, showing the developed scheme to be robust, accurate,
and able to capture discontinuities sharply, without creating spurious
oscillations. Furthermore, the scheme is cost effective. First, because it
allows the use of standard ``black-box'' Poisson solvers, which are normally
tuned to be extremely efficient. Second, because the additional costs of
solving for the correction function scale linearly with the mesh spacing,
which means that they become relatively small for large systems.

Finally, we point out that the present method cannot be applied to all the
interesting situations where a Poisson problem must be solved with jump
discontinuities across an interface. Let us begin by displaying a very
general Poisson problem with jump discontinuities at an interface.
Specifically, consider
\begin{subequations}\label{eq:poissonGen}
 \begin{align}
   \vec{\nabla\/}\cdot\/\p{\beta^+\p{\vec{x}}\vec{\nabla\/}u^+\p{\vec{x}}}
      & = f^+\p{\vec{x}} & \mathrm{for}\;\; \vec{x} & \in \Omega^+\/, \\
   \vec{\nabla\/}\cdot\/\p{\beta^-\p{\vec{x}}\vec{\nabla\/}u^-\p{\vec{x}}}
      & = f^-\p{\vec{x}} & \mathrm{for}\;\; \vec{x} & \in \Omega^-\/, \\
   \left[\alpha\,u\right]_{\Gamma}
      & = a\p{\vec{x}}   & \mathrm{for}\;\; \vec{x} &\in \Gamma\/,
      \label{Gen:a}\\
   \left[(\gamma\,u)_n\right]_{\Gamma} + \left[\eta\,u\right]_{\Gamma}
      & = b\p{\vec{x}}   & \mathrm{for}\;\; \vec{x} &\in \Gamma\/,
      \label{Gen:b}\\
   u\p{\vec{x}}
      & = g\p{\vec{x}}   & \mathrm{for}\;\; \vec{x} &\in \partial\Omega\/,
      \label{Gen:BC}
 \end{align}
\end{subequations}
where we use the notation in \S~\ref{sec:problem} and
\begin{itemize}
 \item[\ref{sec:conclusion}.1\;]
 The brackets indicate jumps across the interface, for example:\\
    $\left[\alpha\,u\right]_{\Gamma} =
    (\alpha^+\,u^+)\p{\vec{x}}-(\alpha^-\,u^-)\p{\vec{x}}\quad$
    for $\;\; \vec{x} \in \Gamma\/$.
 \item[\ref{sec:conclusion}.2\;]
 The subscript $n\/$ indicates the derivative in the direction of $\hat{n}\/$,
 the unit normal to the interface $\Gamma\/$ pointing towards $\Omega^+\/$.
 \item[\ref{sec:conclusion}.3\;]
 $\beta^+> 0\/$ and $f^+\/$ are smooth functions of $\vec{x}$, defined in
    the union of $\Omega^+\/$ and some finite width band enclosing the
    interface $\Gamma\/$.
 \item[\ref{sec:conclusion}.4\;]
 $\beta^-> 0\/$ and $f^-\/$ are smooth functions of $\vec{x}$, defined in
    the union of $\Omega^-\/$ and some finite width band enclosing the
    interface $\Gamma\/$.
 \item[\ref{sec:conclusion}.5\;]
 $\alpha^{\pm} > 0\/$, $\gamma^{\pm} > 0\/$, and $\eta^{\pm}\/$ are smooth
 functions, defined on some finite width band enclosing the interface
 $\Gamma\/$.
 \item[\ref{sec:conclusion}.6\;]
 The Dirichlet boundary conditions in (\ref{Gen:BC}) could be replaced any
    other standard set of boundary conditions on $\partial\,\Omega\/$.
 \item[\ref{sec:conclusion}.7\;]
 As usual, the degree of smoothness of the various data functions involved
    determines how high an order an algorithm can be obtained.
\end{itemize}
Assume now that $\gamma^{\pm} = \alpha^{\pm}\/$,
$\eta^{\pm} = c\,\alpha^{\pm}\/$ --- where $c\/$ is a constant, and that
\begin{equation}\label{Gen:Assum}
 \left.
 \begin{array}{lclcl}
   \vec{A} & = & {\displaystyle
      \frac{1}{\beta^+}\,\vec{\nabla}\/\beta^+}
           & = & {\dfrac{1}{\beta^-}\,\vec{\nabla}\/\beta^-}\/,
    \\ \rule{0mm}{4.5ex}
    B      & = & {\dfrac{\alpha^+}{\beta^+}\,\vec{\nabla}\cdot \left(\beta^+
      \,\vec{\nabla}\,\left( \dfrac{1}{\alpha^+} \right) \right)}
           & = & {\dfrac{\alpha^-}{\beta^-}\,\vec{\nabla}\cdot \left(\beta^-
      \,\vec{\nabla}\,\left( \dfrac{1}{\alpha^-} \right) \right)}\/,
 \end{array}
 \right\}
\end{equation}
applies in the band enclosing $\Gamma\/$ where all the functions are
defined.\footnote{Note that (\ref{Gen:Assum}) implies that $\beta^+\/$
   is a multiple of  $\beta^-\/$.}
In this case the methods introduced in this paper can be used to deal with
the problem in (\ref{eq:poissonGen}) with minimal alterations. The main ideas
carry through, as follows
\begin{itemize}
 \item[\ref{sec:conclusion}.a]
 We assume that both $u^+\/$ and $u^-\/$ can be extended across $\Gamma\/$,
 so that they are defined in some band enclosing the interface.
 \item[\ref{sec:conclusion}.b]
 We define the correction function, in the band enclosing $\Gamma\/$ where
 the $\alpha^{\pm}\/$ and the $u^{\pm}\/$ exist, by
 $D = \alpha^+(\vec{x})\,u^+\p{\vec{x}}-\alpha^-(\vec{x})\,u^-\p{\vec{x}}\/$.
 \item[\ref{sec:conclusion}.c]
 We notice that the correction function satisfies the elliptic Cauchy
 problem
 \begin{equation}
    \nabla^2\,D + \vec{A}\cdot\vec{\nabla}\/D + B\,D =
    \frac{\alpha^+}{\beta^+}\,f^+ - \frac{\alpha^-}{\beta^-}\,f^-\/,
 \end{equation}
 with $D = a\/$ and $D_n = b - c\,a\/$ at the interface $\Gamma\/$.
 \item[\ref{sec:conclusion}.d]
 We notice that the GFM correction terms can be written with knowledge of
 $D\/$.
\end{itemize}
Unfortunately, the conditions in (\ref{Gen:Assum}) exclude some interesting physical
phenomena. In particular, in two-phase flows the case where $\alpha^{\pm}\/$ and $\gamma^{\pm}\/$
are constants (but distinct) and $\eta^{\pm} = 0\/$ arises. We are currently
investigating ways to circumvent these limitations, so as to extend our method
to problems involving a wider range of physical phenomena.


\appendix
\renewcommand{\thesection}{\Alph{section}}
\renewcommand{\thermk}{\Alph{section}.\arabic{rmk}}
\renewcommand{\thefigure}{\Alph{section}.\arabic{figure}}
\setcounter{rmk}{0}
\setcounter{figure}{0}
\section{Bicubic interpolation.} \label{ap:bicubic}
Bicubic interpolation is similar to bilinear interpolation, and can also be
used to represent a function in a rectangular domain. However, whereas bilinear
interpolation requires one piece of information per vertex of the rectangular
domain, bicubic interpolation requires 4 pieces of information: function value,
function gradient, and first mixed derivative (\ie\/~$f_{x\/y}\/$). For
completeness, the relevant formulas for bicubic interpolation are presented
below.

We use the classical multi-index notation, as in Ref.~\cite{nave:10}. Thus, we
represent the 4 vertices of the domain using the vector index
$\vec{v} \in \{0,1\}^2\/$. Namely, the 4 vertices are
$\vec{x}_{\vec{v}} = (x_1^0 + v_1\,\Delta x_1\/,\, x_2^0 + v_2\,\Delta x_2)\/$,
where $(x_1^0\/,\,x_2^0)\/$ are the coordinates of the left-bottom vertex and
$\Delta x_i\/$ is the length of the domain in the $x_i\/$ direction.
Furthermore, given a scalar function $\phi\/$, the 4 pieces of information
needed per vertex are given by
\begin{equation}
 \phi_{\vec{\alpha}}^{\vec{v}} = \partial^{\vec{\alpha}}\phi\p{\vec{x}_{\vec{v}}}\/,
\end{equation}
where both $\vec{v}\/,\,\vec{\alpha} \in \{0,1\}^2$ and
\begin{align}
 \partial^{\vec{\alpha}} = \partial_1^{\alpha_1}\,\partial_2^{\alpha_2}\/, &&
 \partial_i^{\alpha_i} = \p{\Delta x_i}^{\alpha_i}\,
 \dfrac{\partial^{\alpha_i}}{\partial x_i^{\alpha_i}}\/.
\end{align}
Then the 16 polynomials that constitute the standard basis for the bicubic
interpolation can be written in the compact form
\begin{equation}\label{eq:tensor}
 W_{\vec{\alpha}}^{\vec{v}} = \prod_{i=1}^2w_{\alpha_i}^{v_i}\p{\bar{x}_i}\/,
\end{equation}
where $\bar{x}_i = \tfrac{x_i - x_i^0}{\Delta x_i}\/$, and $w_{\alpha}^{v}\/$ is
the cubic polynomial
\begin{equation}\label{eq:basis}
  w_{\alpha}^v(x) = \left\{
  \begin{array}{lll}
   f(x)    & \mbox{for}\; v = 0 & \mbox{and}\; \alpha = 0\/,\\
   f(1-x)  & \mbox{for}\; v = 1 & \mbox{and}\; \alpha = 0\/,\\
   g(x)    & \mbox{for}\; v = 0 & \mbox{and}\; \alpha = 1\/,\\
   -g(1-x) & \mbox{for}\; v = 1 & \mbox{and}\; \alpha = 1\/,
  \end{array}\right.
\end{equation}
where $f(x) = 1 - 3\,x^2 + 2\,x^3\/$ and $g(x) = x\,(1-x)^2\/$.

Finally, the bicubic interpolation of a scalar function $\phi\/$ is given by
the following linear combination of the basis functions:
\begin{equation}\label{eq:int}
 \mathcal{H}\p{\vec{x}} = \sum_{\vec{v}\/,\,\vec{\alpha}\,\in\,\{0,1\}^2}
 W_{\vec{\alpha}}^{\vec{v}}\,\phi_{\vec{\alpha}}^{\vec{v}}
\end{equation}
As defined above (standard bicubic interpolation), 16 parameters are needed to
determine the bicubic. However, in Ref.~\cite{nave:10} a method (``cell-based
approach'') is introduced, that reduces the number of degrees of freedom to 12,
without compromising accuracy. This method uses information from the first
derivatives to obtain approximate formulae for the mixed derivatives. In the
present work, we adopt this cell-based approach.

\setcounter{rmk}{0}
\setcounter{figure}{0}
%
\section{Issues affecting the construction of $\Omega_{\Gamma}^{i,j}\/$.}
\label{ap:omega}
%
%
\subsection{Overview.} \label{sub:overview}
As discussed in \S~\ref{sec:scheme}, the CFM is based on local solutions to
the PDE~\eqref{eq:D} in sub-regions of $\Omega_{\Gamma}$ -- which we call
$\Omega_{\Gamma}^{i,j}\/$. However, there is a certain degree of arbitrariness
in how $\Omega_{\Gamma}^{i,j}\/$ is defined. Here we discuss several factors
that must be consider when solving equation~\eqref{eq:D}, and how they
influence the choice of $\Omega_{\Gamma}^{i,j}\/$. We also present four distinct
approaches to constructing $\Omega_{\Gamma}^{i,j}\/$, of increasing level of
robustness (and, unfortunately, complexity).

The only requirements on  $\Omega_{\Gamma}^{i,j}\/$ that the discussion in
\S~\ref{sec:general} imposes are
\begin{itemize}
 \item
 $\Omega_{\Gamma}^{i,j}$ should be small, since the local problems' condition
 numbers increase exponentially with distance from $\Gamma\/$ --- see
 remarks \ref{rem:4:3} and \ref{rem:4:4}.
 \item
 $\Omega_{\Gamma}^{i,j}$ should contain all the nodes where the correction
 function $D$ is needed.
\end{itemize}
In addition, practical algorithmic considerations further restrict the
definition of $\Omega_{\Gamma}^{i,j}\/$, as explained below.

First, we solve for $D$ in a weak fashion, by locally minimizing a discrete
version of the functional $J_P$ defined in equation~\eqref{eq:jp}. This
procedure involves integrations over $\Omega_{\Gamma}^{i,j}\/$. Thus, it is
useful if $\Omega_{\Gamma}^{i,j}\/$ has an elementary geometrical shape, so
that simple quadrature rules can be applied to evaluate the integrals. Second,
if $\Omega_{\Gamma}^{i,j}\/$ is a rectangle, we can use a (high-order)
bicubic (see appendix~\ref{ap:bicubic}) to represent $D$ in $\Omega_{\Gamma}^{i,j}\/$.
Hence we restrict $\Omega_{\Gamma}^{i,j}\/$ to be a rectangle\footnote{Clearly,
   other simple geometrical shapes, with other types approximations for $D\/$,
   should be possible --- though we have not investigated them.}.
Third, another consideration when constructing $\Omega_{\Gamma}^{i,j}\/$ is how
the interface is represented. In principle the solution to the PDE
in~\eqref{eq:D} depends on the information given along the interface only, and
it is independent of the underlying grid. Nevertheless, consider an interface
described implicitly by a level set function, known only in terms of its values
(and perhaps derivatives) at the grid points. It is then convenient if the
portions of the interface contained within $\Omega_{\Gamma}^{i,j}\/$ can be
easily described in terms of the level set function discretization --- \eg\/
in terms of a pre-determined set of grid cells, such as the cells that define
the discretization stencil. The approaches in \S~\ref{sub:naive} through
\S~\ref{sub:free} are based on this premise.

Although the prior paragraph's strategy makes for an easier implementation, it
also ties $\Omega_{\Gamma}^{i,j}\/$ to the underlying grid, whereas it should
depend only on the interface geometry. Hence it results in definitions for
$\Omega_{\Gamma}^{i,j}\/$ that cannot track the interface optimally. For this
reason, we developed the approach in \S~\ref{sub:node}, which allows
$\Omega_{\Gamma}^{i,j}\/$ to freely adapt to the local interface geometry,
regardless of the underlying grid. The idea is to first identify a piece of
the interface based on a pre-determined set of grid cells, and then use this
information to construct an optimal $\Omega_{\Gamma}^{i,j}\/$. This yields a
somewhat intricate, but very robust definition for $\Omega_{\Gamma}^{i,j}\/$.

Finally, note that explicit representations of the interface are not constrained
by the underlying grid. Moreover, information on the interface geometry is
readily available anywhere along the interface. Hence, in this case, an optimal
$\Omega_{\Gamma}^{i,j}\/$ can be constructed without the need to identify a piece
of the interface in terms of a pre-determined set of grid cells. This fact makes
the approach \S~\ref{sub:node} straightforward with explicit representations of
the interface. By contrast, the less robust approaches in \S~\ref{sub:naive}
through \S~\ref{sub:free} become more involved in this context, because they
requires the additional work of constraining the explicit representation to
the underlying grid.

Obviously, the algorithms presented here (\S~\ref{sub:naive} through
\S~\ref{sub:node}) represent only a few of the possible ways in which
$\Omega_{\Gamma}^{i,j}\/$ can be defined. Nevertheless, these approaches serve
as practical examples of how different factors must be balanced to design
robust schemes.

\subsection{Naive Grid--Aligned Stencil--Centered Approach.}\label{sub:naive}
\begin{figure}[htb!]
 \begin{center}
  \includegraphics[width=2.5in]{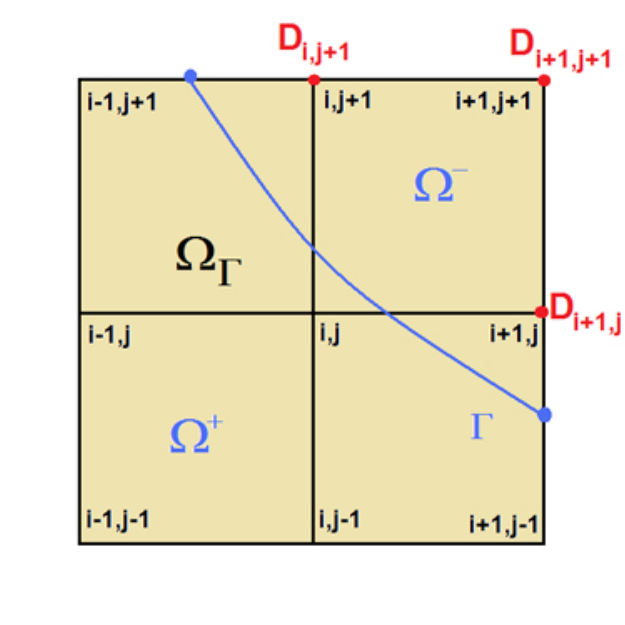}
  \includegraphics[width=2.5in]{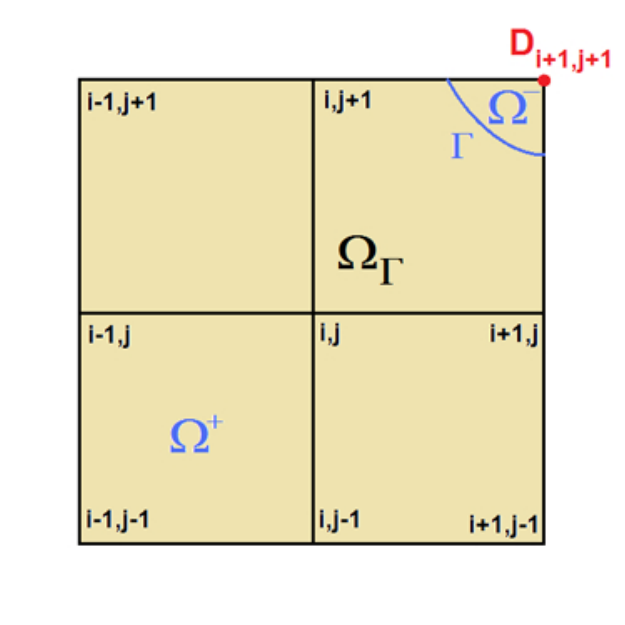}\\
  \text{(a) Well--posed.}\hspace{1.6in}\text{(b) Ill--posed.}
 \end{center}
 \caption{$\Omega_{\Gamma}^{i,j}$ as defined by the naive grid--aligned
          stencil--centered approach.}
 \label{fig:naive}
\end{figure}

In this approach, we fit $\Omega_{\Gamma}^{i,j}\/$ to the underlying grid by
defining it as the $2h_x\times 2h_y$ box that covers the 9-point stencil.
Figure~\ref{fig:naive} shows two examples.

This approach is very appealing because of its simplicity, but it has serious
flaws and \emph{we do not recommend it.} The reason is that the piece of the
interface contained within $\Omega_{\Gamma}^{i,j}\/$ can become arbitrarily small
-- see figure~\ref{fig:naive}(b). Then the arguments that make the local Cauchy
problem well posed no longer apply --- see remarks~\ref{rem:4:3}, \ref{rem:4:4},
and \ref{rem:4:6}. In essence, the biggest frequency encoded in the interface,
$k_{\max} \approx 1/\text{length}(\Gamma/\Omega_{\Gamma}^{i,j}\/)$, can become
arbitrarily large --- while the characteristic length of
$\Omega_{\Gamma}^{i,j}\/$ remains $\mathcal{O}(h)$. As a consequence, the
condition number for the local Cauchy problem can become arbitrarily large. We
describe this approach here merely as an example of the problems that can arise
from too simplistic a definition of $\Omega_{\Gamma}^{i,j}\/$.

\subsection{Compact Grid--Aligned Stencil--Centered Approach.}
\label{sub:compact}
\begin{figure}[htb!]
 \begin{center}
  \hspace*{0ex} \hfill
  \parbox{1.71in}{
  \includegraphics[width=1.7in]{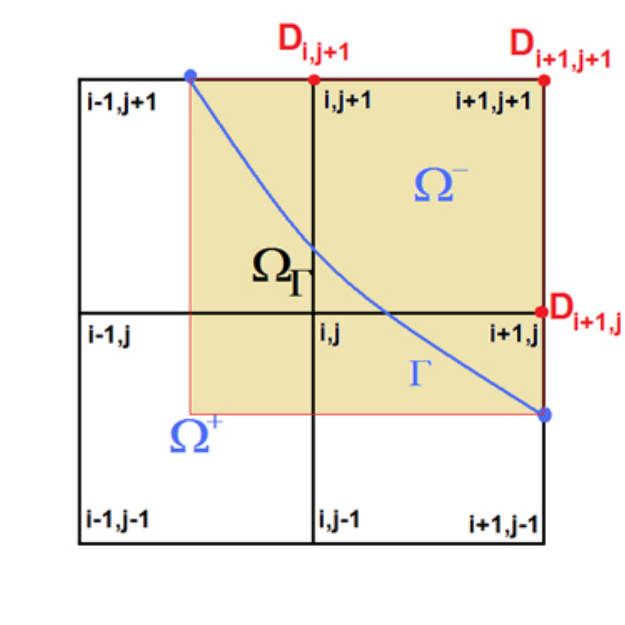}
  \\ \hspace*{0ex} \hfill  (a) Well balanced. \hfill \hspace*{0ex}}
  \hfill
  \parbox{1.71in}{
  \includegraphics[width=1.7in]{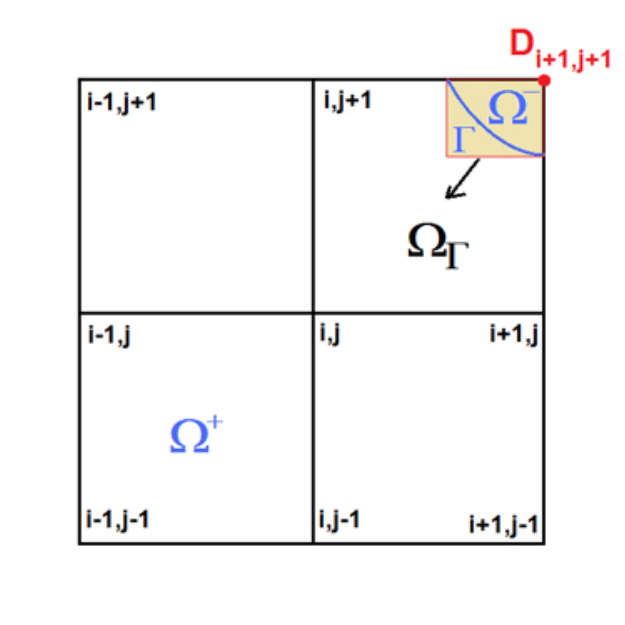}
  \\ \hspace*{0ex} \hfill  (b) Well balanced. \hfill \hspace*{0ex}}
  \hfill
  \parbox{1.71in}{
  \includegraphics[width=1.7in]{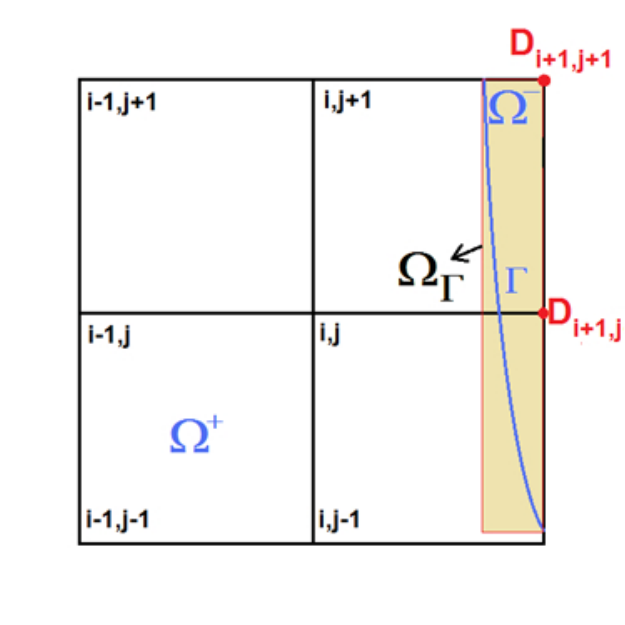}
  \\ \hspace*{0ex} \hfill  (c) Elongated. \hfill \hspace*{0ex}}
  \hfill \hspace*{0ex} \\  \hspace*{1.5ex}
 \end{center}
 \caption{$\Omega_{\Gamma}^{i,j}$ as defined by the compact grid--aligned
          stencil--centered approach.}
 \label{fig:compact}
\end{figure}

This is the approach described in detail in \S~\ref{sub:omega}. In summary:
$\Omega_{\Gamma}^{i,j}\/$ is defined as the smallest rectangle that
\begin{enumerate}[(i)]
 \item Is aligned with the grid.
 \item Includes the piece of the interface contained within the stencil.
 \item Includes all the nodes where $D$ is needed.
\end{enumerate}
Figure~\ref{fig:compact} shows three examples of this definition. As it
should be clear from this figure, a key consequence of (i--iii) is that the
piece of interface contained within $\Omega_{\Gamma}^{i,j}\/$ is always close to
its diagonal --- hence it is never too small relative to the size of
$\Omega_{\Gamma}^{i,j}\/$. Consequently, this approach is considerably more
robust than the one in \S~\ref{sub:naive}. In fact, we successfully employed it
for all the examples using the \fth\/ order accurate scheme ---
see \S~\ref{sec:results}.

Unfortunately, the requirements (i--ii) in this approach tie
$\Omega_{\Gamma}^{i,j}\/$ to the grid and the stencil. As mentioned earlier, these
constraints may lead to an $\Omega_{\Gamma}^{i,j}\/$ which is not the best fit to
the geometry of the interface. Figure~\ref{fig:compact}(c) depicts a situation
where this strategy yields relatively poor results. This happens when there is
an almost perfect alignment of the interface with the grid, 
which can result in an
excessively elongated $\Omega_{\Gamma}^{i,j}\/$ --- in the worse case scenario,
this set could reduce to a line. Although the local Cauchy problem remains well
conditioned,
the elongated sets can interfere with the process we use to solve the equation
for $D$ in~\eqref{eq:D}. 
Essentially, the representation of a
function by a bicubic (or a modified bilinear in the case of the \snd\/ accurate
scheme in appendix~\ref{ap:scheme2}) becomes an ill-defined problem as the aspect
ratio of the rectangle $\Omega_{\Gamma}^{i,j}\/$ vanishes (however, see the
next paragraph). In the authors' experience, when this sort of alignment
happens the solution remains valid and the errors are still relatively small
--- but the convergence rate may be affected if the bad alignment persists over
several grid refinements.

We note that we observed the difficulties described in the paragraph above with
the \snd\/ accurate scheme in appendix~\ref{ap:scheme2} only. We attribute this to
the fact that
the bicubics result in a much better enforcement of the PDE~\eqref{eq:D} than the
modified bilinears. The latter are mostly determined by the interface conditions
--- see remark~\ref{rem:C:1}. Hence, the PDE~\eqref{eq:D} provides a much stronger
control over the bicubic parameters, making this interpolation more robust in
elongated sets.

Note that this issue is much less severe than the problem affecting the
approach in \S~\ref{sub:naive}, and could be corrected by making the Cauchy
solver ``smarter'' when dealing with elongated sets. Instead we adopted the
simpler solution of having a definition for $\Omega_{\Gamma}^{i,j}\/$ that avoids
elongated sets. The most robust way to do this is to abandon the requirements
in (i--ii), and allow $\Omega_{\Gamma}^{i,j}\/$ to adapt to the local geometry
of the interface. This is the approach introduced in \S~\ref{sub:node}. A
simpler, compromise solution, is presented in \S~\ref{sub:free}.
%
\subsection{Free Stencil--Centered Approach.} \label{sub:free}
Here we present a compromise solution for avoiding elongated
$\Omega_{\Gamma}^{i,j}\/$, which abandons the requirement in (i), but not in (ii)
--- since (ii) is convenient when the interface is represented implicitly. In
this approach $\Omega_{\Gamma}^{i,j}\/$ is defined as the smallest rectangle that
\begin{enumerate}[(i*)]
 \item Is aligned with the grid rotated by an angle $\theta_r$, where
 $\theta_r = \theta_{\Gamma} - \pi/4$ and $\theta_{\Gamma}$ characterizes the
 interface alignment with respect to the grid (\eg\/ the polar angle for the
 tangent vector to the interface section at its mid-point inside the stencil).
 \item Includes the piece of the interface contained within the stencil.
 \item Includes all the nodes where $D$ is needed.
\end{enumerate}
Figure~\ref{fig:compact} shows two examples of this approach.
\begin{figure}[htb!]
 \begin{center}
  \includegraphics[width=2.1in]{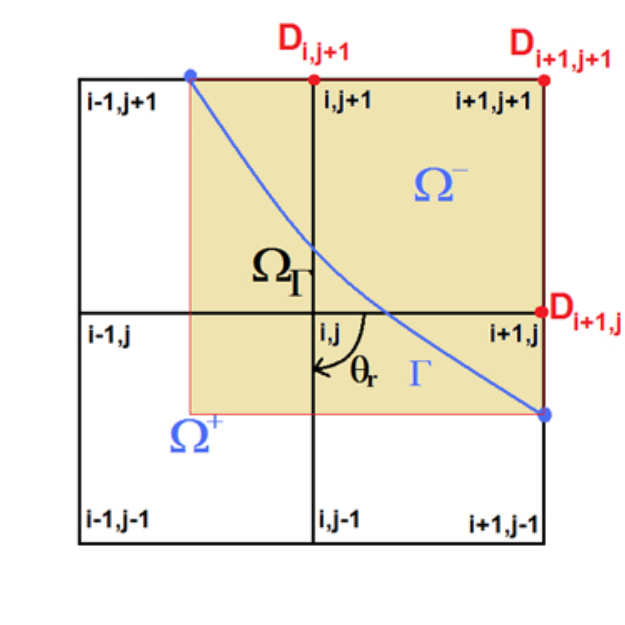}
  \includegraphics[height=2.1in]{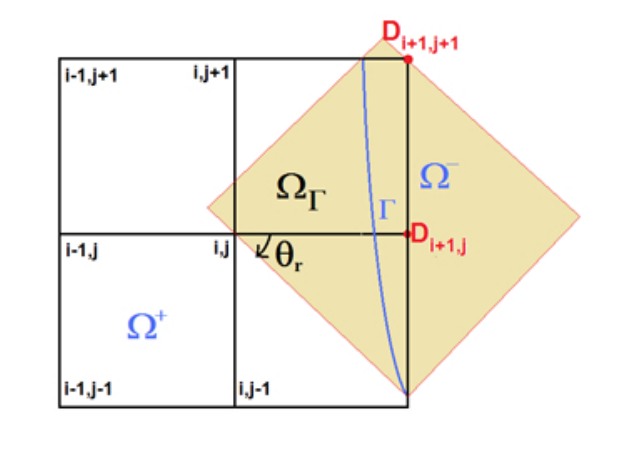}\\
  \text{(a)}\hspace{2.0in}\text{(b)\hspace{0.5in}}
 \end{center}
 \caption{$\Omega_{\Gamma}^{i,j}$ as defined by the free stencil--centered
           approach.}
 \label{fig:free}
\end{figure}

The implementation of the present approach is very similar to that of the one
in \S~\ref{sub:compact}. The only additional work is to compute $\theta_r$ and
to write the interface and points where $D$ is needed in the rotated frame of
reference. In both these approaches the diagonal of $\Omega_{\Gamma}^{i,j}$ is
very close to the piece of interface contained within the stencil, which
guarantees a  well conditioned local problem. However, here the addition of a
rotation keeps  $\Omega_{\Gamma}^{i,j}$ nearly square, and avoids elongated
geometries. The price paid for this is that the sets $\Omega_{\Gamma}^{i,j}$
created using \S~\ref{sub:free} can be a little larger than the ones from
\S~\ref{sub:compact} --- with both sets including the exact same piece of
interface. In such situations, the present approach results in a somewhat
larger condition number.
%
\subsection{Node-Centered Approach.} \label{sub:node}
Here we define $\Omega_{\Gamma}^{i,j}$ in a fashion that is completely independent
from the underlying grid and stencils. In fact, instead of associating
each $\Omega_{\Gamma}^{i,j}$ to a particular stencil, we define a different
$\Omega_{\Gamma}^{i,j}$ for each node where the correction is needed ---
hence the name node-centered, rather than stencil-centered. As a consequence,
whereas the prior strategies lead to multiple values of $D$ at the same node
(one value per stencil, see remark~\ref{rem:5:2}), here there is a unique
value of $D$ at each node.

\begin{figure}[htb!]
 \begin{center}
  \includegraphics[width=2.3in]{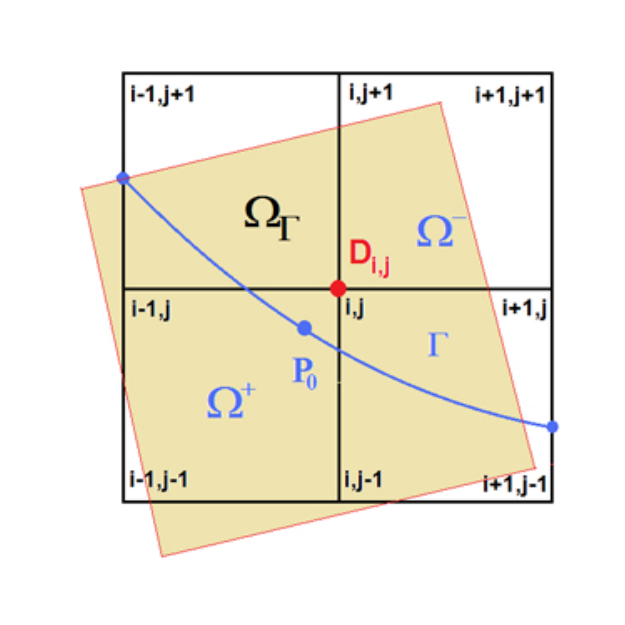}
  \includegraphics[width=2.3in]{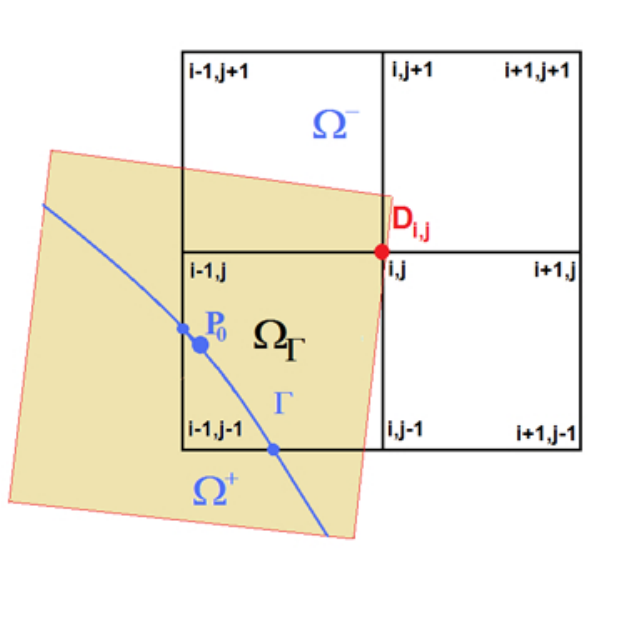}\\
  \text{(a)}\hspace{2.1in}\text{(b)\hspace{0.5in}}
 \end{center}
 \caption{$\Omega_{\Gamma}^{i,j}$ as defined by the node--centered approach.}
 \label{fig:node}
\end{figure}

In this approach, $\Omega_{\Gamma}^{i,j}$ is defined by the following steps:
\begin{enumerate}[1.]
 \item 
 Identify the interface in the 4 grid cells that surround a given node. This
 step can be skipped if the interface is represented explicitly.
 \item 
 Find the point, $P_0\/$, along the interface that is closest to the node. This
 point becomes the center of $\Omega_{\Gamma}^{i,j}$. There is no need to obtain
 $P_0$ very accurately. Small errors in $P_0$ result in small shifts in
 $\Omega_{\Gamma}^{i,j}$ only.
 \item 
 Compute $\hat{t}_0\/$, the vector tangent to the interface at $P_0\/$. This
 vector defines one of the diagonals of $\Omega_{\Gamma}^{i,j}\/$. The normal
 vector $\hat{n}_0$ defines the other diagonal. Again, high accuracy is not
 needed.
 \item 
 Then $\Omega_{\Gamma}^{i,j}\/$ is the square with side length
 $2\sqrt{h_x^2 + h_y^2}\/$, centered at $P_0$, and diagonals parallel to
 $\hat{t}_0\/$ and $\hat{n}_0\/$ --- $\Omega_{\Gamma}^{i,j}\/$ need not
 be aligned with the grid.
\end{enumerate}
Figure~\ref{fig:node} shows two examples of this approach.
%

Note that the piece of interface contained within $\Omega_{\Gamma}^{i,j}\/$, as
defined by steps 1--4 above, is not necessarily the same found in step 1. Hence,
after defining $\Omega_{\Gamma}^{i,j}$, we still need to identify the piece of
interface that lies within it. For an explicit representation of the interface,
this additional step is not particularly costly, but the same is not true for an
implicit representation.

This approach is very robust because it always creates a square
$\Omega_{\Gamma}^{i,j}\/$, with the interface within it close to one of the
diagonals --- as guaranteed by steps 2 and 3. Hence, the local Cauchy problem
is always well conditioned. Furthermore, making $\Omega_{\Gamma}^{i,j}\/$ square
(to avoid elongation) does not result in larger condition numbers as in
\S~\ref{sub:free} because the larger $\Omega_{\Gamma}^{i,j}\/$ contain an
equally larger piece of the interface.

Finally, we point out that the (small) oscillations observed in the convergence
plots shown in \S~\ref{sec:results} occur because in these calculations we use
the approach in \S~\ref{sub:compact} --- which produces sets
$\Omega_{\Gamma}^{i,j}$ that are not uniform in size, nor shape, along the
interface. Tests we have done show that these oscillations do not occur with
the node-centered approach here, for which all the $\Omega_{\Gamma}^{i,j}\/$ are
squares of the same size. Unfortunately, as pointed out earlier, the
node-centered approach is not well suited for calculations using an interface
represented implicitly.
\setcounter{rmk}{0}
\setcounter{figure}{0}
%
\section{\snd\/ Order Accurate Scheme in 2D.} \label{ap:scheme2}
%
\subsection{Overview.}\label{sub:scheme2:overview}
In \S~\ref{sec:scheme} we present a \fth\/ order accurate scheme to solve the
2D Poisson problem with interface jump conditions, based on the correction
function defined in \S~\ref{sec:general}. However, there are many situations
in which a \snd\/ order
version of the method could be of practical relevance. Hence, in this section
we use the general framework provided by the correction function method (CFM)
to develop a specific example of a \snd\/ order accurate scheme in 2D. The basic
approach is analogous to the \fth\/ version presented in \S~\ref{sec:scheme}.
A few key points are
\begin{enumerate}[(a)]
 \item 
 We discretize Poisson's equation using the standard 5-point stencil. This
 stencil is compact, which is an important requirement for a well conditioned
 problem for the correction function $D\/$ --- see
 remarks~\ref{rem:4:2}--\ref{rem:4:4}.
 \item 
 We approximate $D\/$ using modified bilinear interpolants (defined in
 \S~\ref{sub:bilinear}), each valid in a small neighborhood
 $\Omega_{\Gamma}^{i,j}\/$ of the interface. This guarantees local \snd\/ order
 accuracy with only 5 interpolation parameters. Each $\Omega_{\Gamma}^{i,j}\/$
 corresponds to a grid point at which the standard discretization of
 Poisson's equation involves a stencil that straddles the interface $\Gamma\/$.
 \item 
 The domains $\Omega_{\Gamma}^{i,j}\/$ are rectangular regions, each enclosing a
 portion of $\Gamma\/$, and all the nodes where $D\/$ is needed to complete the
 discretization of Poisson's equation at the ($i\/,\,j$)-th stencil. Each is
 a sub-domain of $\Omega_\Gamma\/$.
 \item 
 Starting from (b) and (c), we design a local solver that provides an
 approximation to $D\/$ inside each domain $\Omega_{\Gamma}^{i,j}\/$.
 \item 
 The interface $\Gamma\/$ is represented using the standard level set approach
 --- see~\cite{osher:88}. This guarantees a local \snd\/ order representation
 of the interface, as required to keep the overall accuracy of the scheme.
 \item 
 In each $\Omega_{\Gamma}^{i,j}\/$, we solve the PDE in \eqref{eq:D} in a least
 squares sense --- see remark~\ref{rem:5:1}. Namely, we seek the minimum of the
 positive quadratic integral quantity $J_P\/$ in~\eqref{eq:jp}, which
 vanishes at the solution: We substitute the modified bilinear approximation
 for $D\/$ into $J_P\/$, discretize the integrals using Gaussian quadratures,
 and minimize the resulting discrete $J_P\/$.
\end{enumerate}
%
\begin{rmk}\label{rem:C:1}
 Using standard bilinear interpolants to approximate $D\/$ in each
 $\Omega_{\Gamma}^{i,j}\/$ also yields \snd\/ order accuracy. However, the
 Laplacian of a standard bilinear interpolant vanishes. Thus, this basis cannot
 take full advantage of the fact that $D\/$ is the solution to the Cauchy
 problem in \eqref{eq:D}. The modified bilinears, on the other hand, incorporate
 the average of the Laplacian into the formulation --- see
 \S~\ref{sub:bilinear}. In \S~\ref{sub:results2} we include results with both
 the standard and the modified bilinears, to demonstrate the advantages of the
 latter.
\myremarkend
\end{rmk}

Below, in \S~\ref{sub:bilinear}--\ref{sub:solution2} we describe
the \snd\/ order accurate scheme, and in \S~\ref{sub:results2} we present
applications of the scheme to three test cases.
%
\subsection{Modified Bilinear.}\label{sub:bilinear}
The modified bilinear interpolant used here builds on the standard bilinear
polynomials. To define them, we use the multi-index notation --- see
appendix~\ref{ap:bicubic} and Ref.~\cite{nave:10}. Thus, we label
the 4 vertices of a rectangular cell using the vector index
$\vec{v} \in \{0,1\}^2\/$. Then the vertices are
$\vec{x}_{\vec{v}} = (x_1^0 + v_1\,\Delta x_1\/,\, x_2^0 + v_2\,\Delta x_2)\/$,
where $(x_1^0\/,\,x_2^0)\/$ are the coordinates of the left-bottom vertex and
$\Delta x_i\/$ is the length of the domain in the $x_i\/$ direction. The
standard bilinear interpolation basis is then given by the 4 polynomials
\begin{equation}
 W^{\vec{v}} = \prod_{i=1}^2w^{v_i}\p{\bar{x}_i}\/,
\end{equation}
where $\bar{x}_i = \tfrac{x_i - x_i^0}{\Delta x_i}\/$, and $w^{v}\/$ is the
linear polynomial
\begin{equation}
  w^v(x) = \left\{
  \begin{array}{lll}
   1-x  & \mbox{for}\; v = 0\/,\\
   x    & \mbox{for}\; v = 1\/.
  \end{array}\right.
\end{equation}
The standard bilinear interpolation of a scalar function $\phi\/$ is given (in
each cell) by
\begin{equation}
 \mathcal{H}_s\p{\vec{x}} = \sum_{\vec{v}\,\in\,\{0,1\}^2}
  W^{\vec{v}}\,\phi\p{\vec{x}_{\vec{v}}}.
\end{equation}
In the modified version, we add a quadratic term proportional to $x^2 + y^2$,
so that the Laplacian of the modified bilinear is no longer identically zero.
The coefficient of the quadratic term can be written in terms of the average
value of the Laplacian over the domain, $\overline{\nabla^2\,\phi}$. This yields
the following formula for the modified bilinear interpolant
\begin{equation}
 \mathcal{H}_m\p{\vec{x}} = \mathcal{H}_s\p{\vec{x}} -
 \dfrac{1}{4}\left[w^0(\bar{x})w^1(\bar{x})(\Delta x)^2 +
 w^0(\bar{y})w^1(\bar{y})(\Delta y)^2\right]\overline{\nabla^2\,\phi}\/.
\end{equation}
%
\subsection{Standard Stencil.}
We use the standard \snd\/ order accurate 5-point discretization of the
Poisson equation
\begin{equation}\label{eq:5p}
 L^5 u_{i,j} = f_{i,j}\/,
\end{equation}
where $L^5\/$ is defined in \eqref{eq:L5}. In the absence of discontinuities,
\eqref{eq:5p} provides a compact \snd\/ order accurate representation of the
Poisson equation. In the vicinity of the discontinuities at the interface
$\Gamma\/$, we define an appropriate domain $\Omega_{\Gamma}^{i,j}\/$, and use
it to compute the correction terms needed by \eqref{eq:5p} --- as
described below.

\begin{figure}[htb!]
 \begin{center}
  \includegraphics[width=2.5in]{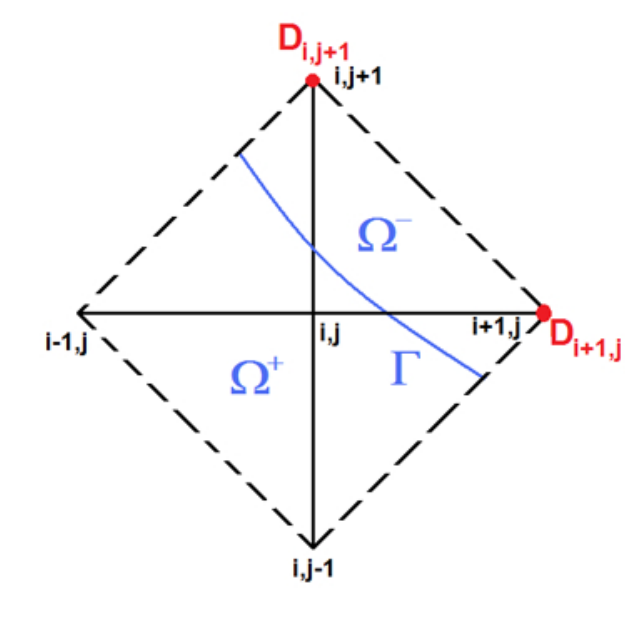}
 \end{center}
 \caption{The 5-point stencil next to the interface $\Gamma$. The dashed box
          shows a compact quadrangular region that contains the stencil.}
 \label{fig:D2}
\end{figure}

To understand how the correction terms affect the discretization,
consider the situation in figure~\ref{fig:D2}. In this case, the node
$(i,j)$ lies in $\Omega^+$ while the nodes $(i+1,j)$ and $(i,j+1)$ are in
$\Omega^-$. Hence, to be able to use equation~\eqref{eq:5p}, we need to
compute $D_{i+1,j}$ and $D_{i,j+1}$.
%

After having solved for $D$ where necessary (see \S~\ref{sub:omega2} and
\S~\ref{sub:solution2}), we modify equation~\eqref{eq:5p} and write
\begin{equation}\label{eq:5p-C}
  L^5 u_{i,j} = f_{i,j} + C_{i,j}\/,
\end{equation}
which differs from \eqref{eq:5p} by the terms $C_{i,j}$ on the RHS only. Here
the $C_{i,j}$ are the CFM correction terms needed to complete the stencil
across the discontinuity at $\Gamma\/$. In the particular case illustrated by
figure~\ref{fig:D2}, we have
\begin{equation}\label{eq:5p-D}
  C_{i,j} = - \dfrac{1}{h_x^2}D_{i+1,j} - \dfrac{1}{h_y^2}D_{i,j+1}\/.
\end{equation}
Similar formulas apply for all the other possible arrangements of the stencil
for the Poisson's equation, relative to the interface $\Gamma\/$. 
%
\subsection{Definition of $\Omega_{\Gamma}^{i,j}\/$.} \label{sub:omega2}
As discussed in appendix~\ref{ap:omega}, the construction of
$\Omega_{\Gamma}^{i,j}\/$ presented in \S~\ref{sub:omega} may lead to elongated
shapes for $\Omega_{\Gamma}^{i,j}\/$, which can cause accuracy loses. As
mentioned earlier (see the end of \S~\ref{sub:compact}) this is not a problem
for the \fth\/ order scheme, but it can be one for the  \snd\/ order one.
To resolve this issue we could have implemented the robust construction of
$\Omega_{\Gamma}^{i,j}\/$ given in \S~\ref{sub:node}. However the simpler
compromise version in \S~\ref{sub:free} proved sufficient.

\begin{figure}[htb!]
 \begin{center}
  \includegraphics[width=2.5in]{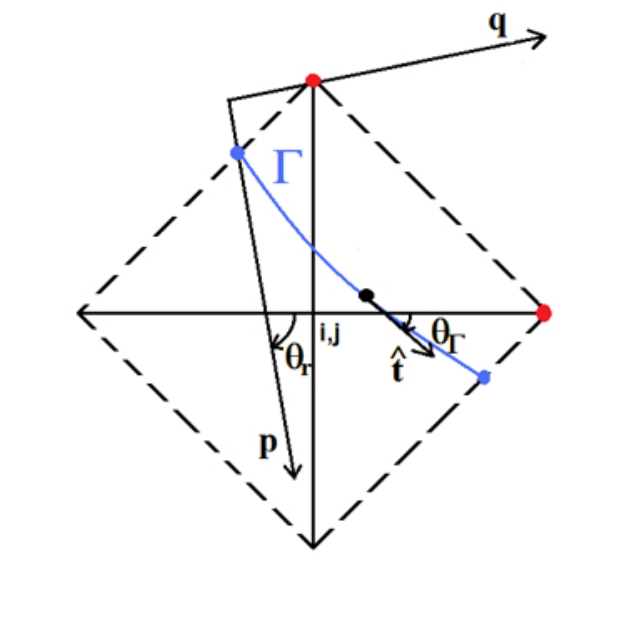}
  \includegraphics[width=2.5in]{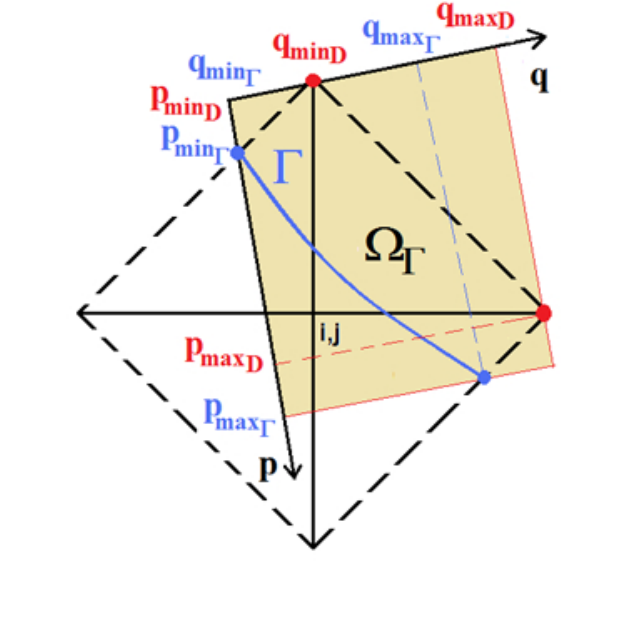}\\
  \text{(a) Step 1.}\hspace{1.5in} \text{(b) Steps 2--4.}
 \end{center}
 \caption{The set $\Omega_{\Gamma}^{i,j}$ for the situation in
          figure~\ref{fig:D2}.}
 \label{fig:omegag2}
\end{figure}

The approach in \S~\ref{sub:free} requires $\Omega_{\Gamma}^{i,j}\/$ to include
the piece of interface contained ``within the stencil.'' For the 9-point stencil,
this naturally means ``within the $2h_x\times 2h_y$ box aligned with the grid
that includes the nine points of the stencil.'' We could use the same meaning
for the 5-point stencil, but this would not take full advantage of the 5-point
stencil compactness. A better choice is to use the quadrilateral defined by the
stencil's four extreme points --- \ie\/ the dashed box in figure~\ref{fig:D2}.
This choice is used only to determine the piece of interface to be included
within $\Omega_{\Gamma}^{i,j}\/$. It does not affect the definition of the
interface alignment angle $\theta_L\/$ used by the \S~\ref{sub:free} approach.
Hence the following four steps define $\Omega_{\Gamma}^{i,j}\/$
\begin{enumerate}[1.]
 \item 
 Find the angle $\theta_{\Gamma}\/$ between the vector tangent to the interface
 at the mid-point within the stencil, and the $x$-axis. Introduce the
 coordinate system $p$-$q$, resulting from rotating $x$-$y$ by
 $\theta_r = \theta_{\Gamma} - \pi\//4$ --- see figure~\ref{fig:omegag2}(a).
 \item 
 Find the coordinates $(p_{\min_{\Gamma}}\/,\,p_{\max_{\Gamma}})\/$ and
 $(q_{\min_{\Gamma}}\/,\,q_{\max_{\Gamma}})$ characterizing the smallest $p$-$q$
 coordinate rectangle enclosing the section of the interface contained within
 the stencil --- see figure~\ref{fig:omegag2}(b). 
 \item 
 Find the coordinates $(p_{\min_D}\/,\,p_{\max_D})\/$ and
 $(q_{\min_D}\/,\,q_{\max_D})\/$ characterizing the smallest $p$-$q$ coordinate
 rectangle enclosing all the nodes at which $D\/$ is needed --- see
 figure~\ref{fig:omegag2}(b).
 \item 
 $\Omega_{\Gamma}^{i,j}\/$ is the smallest $p$-$q$ coordinate rectangle enclosing
 the two previous rectangles. Its edges are characterized by
  \begin{subequations}
   \begin{align}
    p_{\min} &= \min\p{p_{\min_{\Gamma}},p_{\min_D}},\\
    p_{\max} &= \max\p{p_{\max_{\Gamma}},p_{\max_D}},\\
    q_{\min} &= \min\p{q_{\min_{\Gamma}},q_{\min_D}},\\
    q_{\max} &= \max\p{q_{\max_{\Gamma}},q_{\max_D}}.
   \end{align}
  \end{subequations}
\end{enumerate}
Figure~\ref{fig:omegag2} shows an example of $\Omega_{\Gamma}^{i,j}$ defined in
this way.
%
%
%
\subsection{Solution of the Local Cauchy Problem.} \label{sub:solution2}
The remainder of the \snd\/ order accurate scheme follows the exact same
procedure used for the \fth\/ order accurate version described in detail in
\S~\ref{sec:scheme}. In particular, as explained in item~(f) of
\S~\ref{sub:scheme2:overview}, we solve the local Cauchy problem defined by
equation~\eqref{eq:D} in a least squares sense (using the same minimization
procedure described in \S~\ref{sub:solution}). The only differences are that:
(i) In the \snd\/ order accurate version of the method we use 4 Gaussian
quadrature points for the 1D line integrals, and 16 points for the 2D area
integrals. (ii) Because the modified bilinear representation for $D\/$ involves
5 basis polynomials, the minimization problem produces a $5\times5$
(self-adjoint) linear system --- instead of the  $12\times12$ system that
occurs for the \fth order algorithm.
%
\subsection{Results}\label{sub:results2}
Here we present three examples of solutions to the 2D Poisson's equation using
the algorithm described above. Each example is defined below in terms of the
problem parameters (source term $f\/$ and jump conditions across the
interface $\Gamma\/$), the representation of the interface, and the exact
solution (needed to evaluate the errors in the convergence plots). Note that
\begin{enumerate}
 \item 
 In all cases we represent the interface using a level set function $\phi\/$. 
 \item 
 Below the level set function is described via an analytic formula. This formula
 is converted into a discrete level set representation used by the code. Only
 this discrete representation is used for the actual computations.
 \item 
 The level set formulation allows us to compute the vectors normal to the
 interface to \snd\/ order accuracy with a combination of finite differences and
 standard bilinear interpolation. Hence, it is convenient to write the jump in
 the normal derivative, $\left[u_n\right]_{\Gamma}\/$, in terms of the jump in
 the gradient of $u\/$ dotted with the normal to the interface
 $\hat{n} = (n_x\/,\,n_y)\/$.
\end{enumerate}
We use\footnote{A subscript s is added to the example numbers of the \snd
   order method, to avoid confusion with the \fth order method examples.}
example~1$_s$ to test the \snd\/ order scheme in a generic problem with a
non-trivial interface geometry. In addition, in this example the correction
function $D\/$ has a non-zero Laplacian. Thus, we used this problem to compare
the performance of the Standard Bilinear (SB) and Modified Bilinear (MB) interpolations
as basis for the correction function. Finally, examples~2$_s$ and 3$_s$ here
correspond to examples~1 and 3 in~\cite{leveque:94}, respectively. Hence
they can be used to compare the performance of the present \snd\/ order scheme
with the Immersed Interface Method (IIM), in two distinct situations.
\subsubsection{Example 1$_s$.}\label{subsub:case1}
\begin{figure}[htb!]
 \begin{center}
  \includegraphics[width=3.5in]{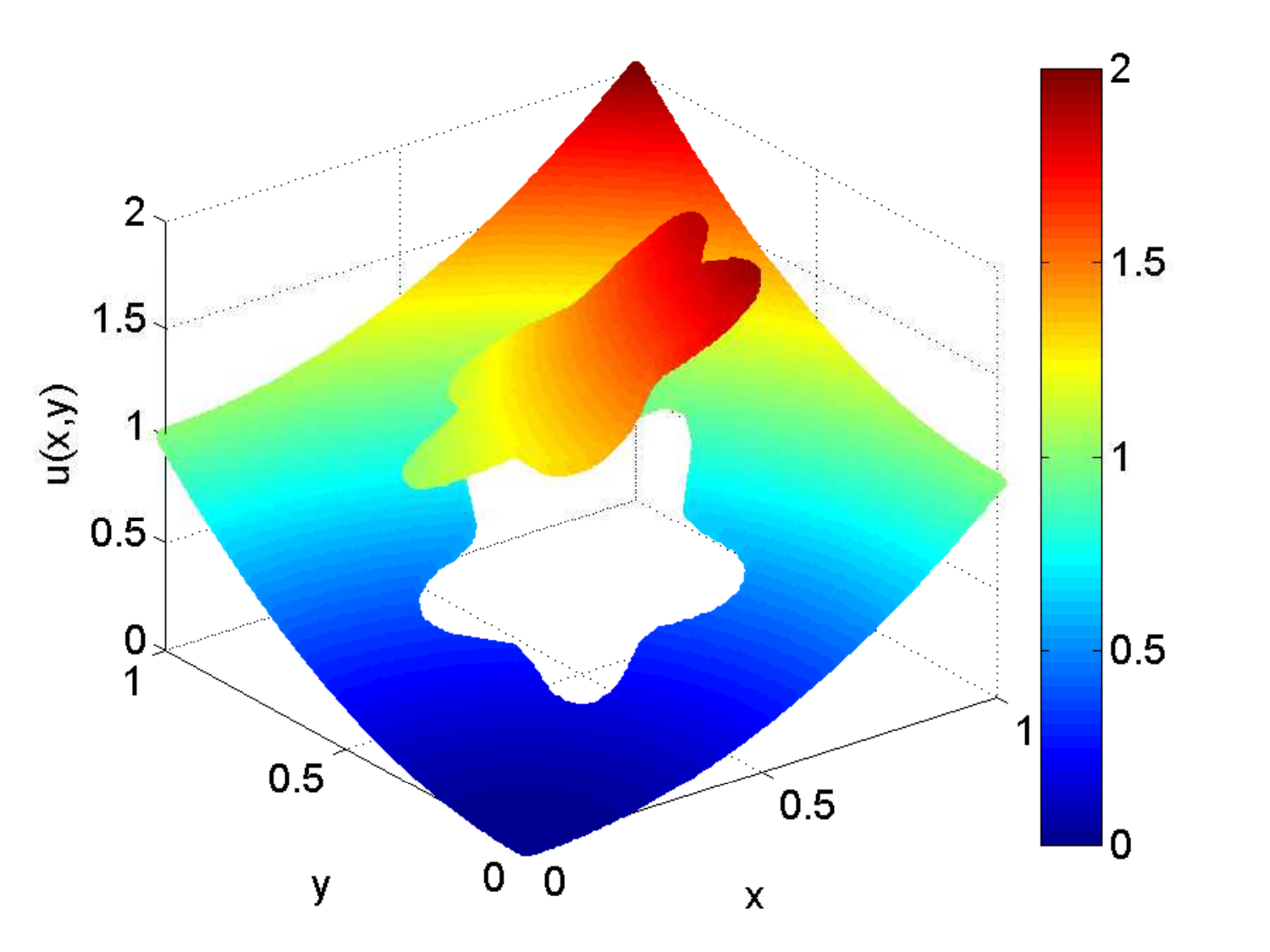}
 \end{center}
 \caption{Example 1$_s$. Numerical solution with $193 \times 193\/$ nodes,
          \snd order scheme, modified bilinear version.}
 \label{fig:case1-2-numerical}
\end{figure}

\begin{itemize}
 \item 
 Domain: $(x,y) \in [0,1] \times [0,1]$.
 \item 
 Problem parameters:
 \begin{align*}
   f^+\p{x\/,\,y}           & = 4\/,\\
   f^-\p{x\/,\,y}           & = 0\/,\\
   [u]_{\Gamma}              & = x^2 + y^2 - \exp(x)\cos(y)\/,\\
   \left[u_n\right]_{\Gamma} & = \left[2x - \exp(x)\cos(y)\right]\,n_x +
                                \left[2y + \exp(x)\sin(y)\right]\,n_y\/.
 \end{align*}
 \item 
 Level set defining the interface: \hfill
 $\phi\p{x\/,\,y} = \sqrt{(x-x_0) + (y-x_0)} - r(\theta)\/$,\\
 where $r(\theta) = r_0 + \epsilon\,\sin(5\,\theta)\/$,
       $\theta\p{x\/,\,y} = \arctan\p{\dfrac{y-y_0}{x-x_0}}\/$,
       $x_0 = 0.5\/$, $y_0 = 0.5\/$, $r = 0.25\/$, and $\epsilon = 0.05\/$.
 \item 
 Exact solution:
 \begin{align*}
   u^+\p{x\/,\,y} & = x^2 + y^2\/,\\
   u^-\p{x\/,\,y} & = \exp(x)\cos(y)\/.
 \end{align*}
\end{itemize}
Figure~\ref{fig:case1-2-numerical} shows the numerical solution with a fine grid
($193 \times 193$ nodes). The non-trivial contour of the interface is accurately
represented and the discontinuity is captured very sharply. For comparison we
solved this problem using both the SB and MB interpolations to represent the
correction function. Although figure~\ref{fig:case1-2-numerical} shows the
solution obtained with the modified bilinear version, both versions produce
small errors and are visually indistinguishable.

\begin{figure}[htb!]
 \begin{center}
  \includegraphics[width=2.6in]{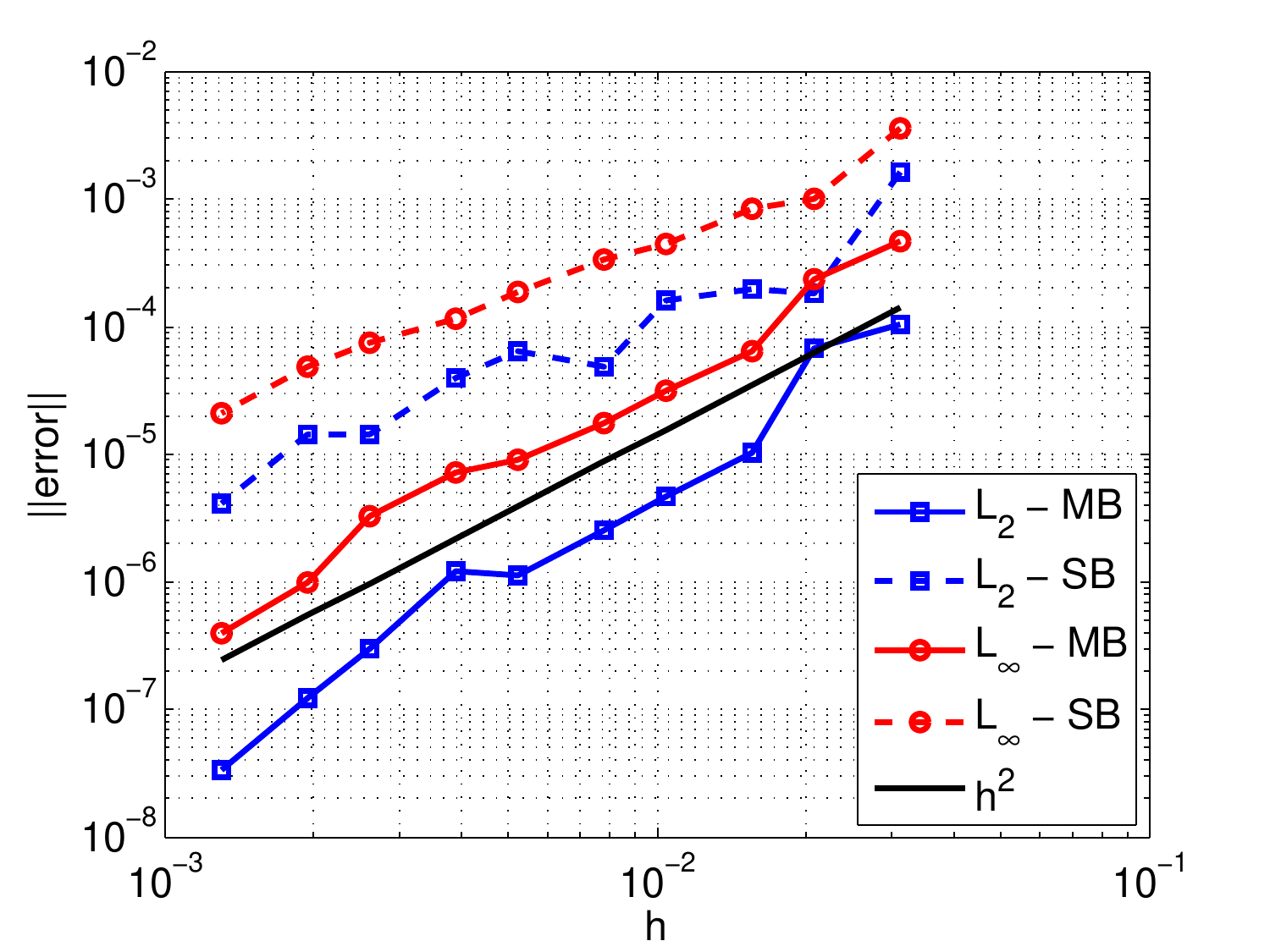}
  \includegraphics[width=2.6in]{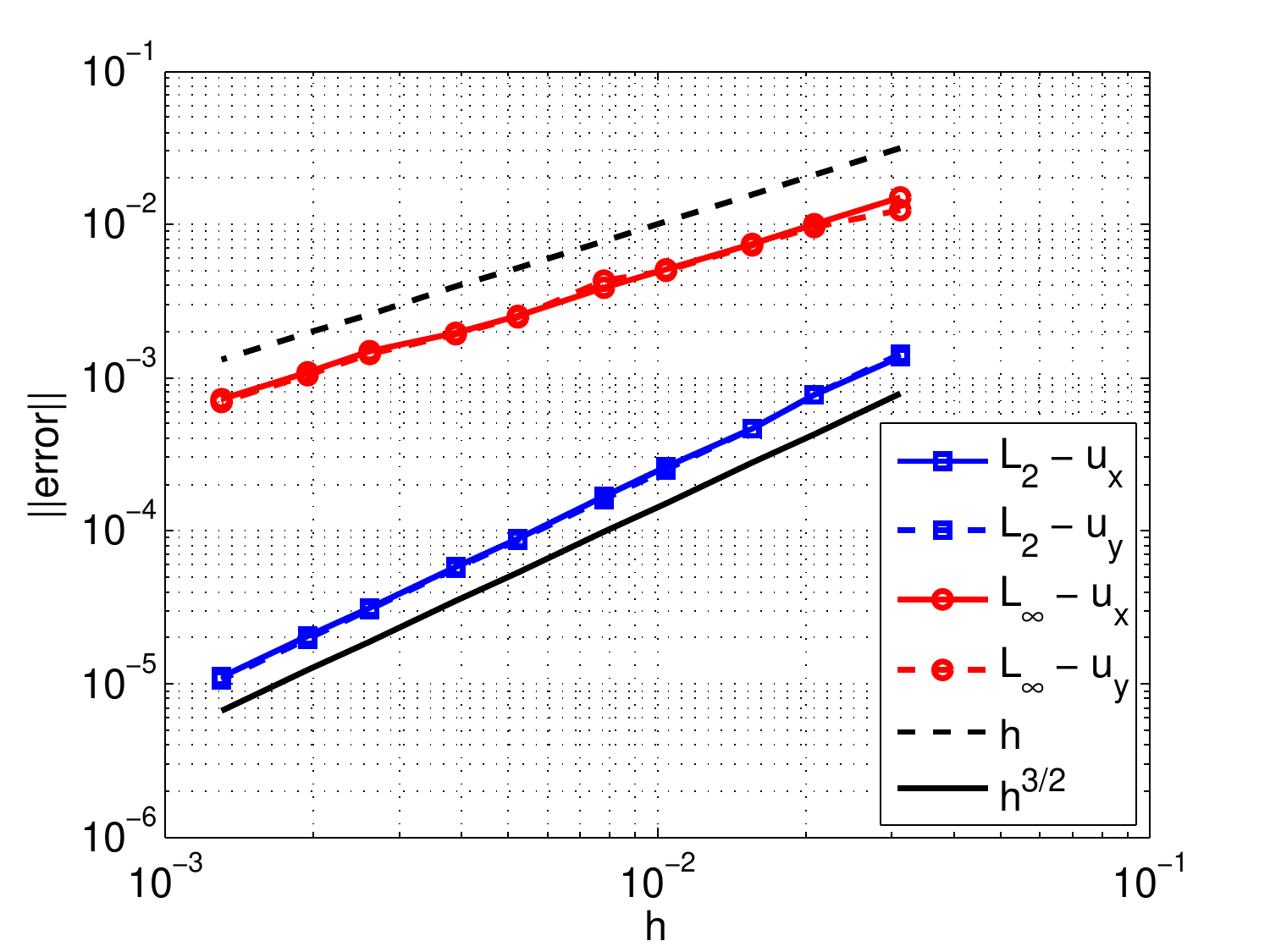}\\
  (a) Solution. \hspace*{1.6in} (b) Gradient.\hfill
 \end{center}
 \caption{Example 1$_s$ - Error behavior of the solution and its gradient in the $L_2\/$ and $L_\infty\/$ norms.}
 \label{fig:case1-2-convergence}
\end{figure}

Figure~\ref{fig:case1-2-convergence}(a) shows the convergence of the solution error
in the $L_2\/$ and $L_{\infty}\/$ norms, for both the SB and MB versions. As expected,
the overall behavior indicates \snd\/ order convergence, despite the small
oscillations that are characteristic of this implementation of the method, as
explained in \S~\ref{sub:case2}. Note that the MB version produces
significantly smaller errors (a factor of about 20) than the SB version. It
also exhibits a more robustly \snd\/ order convergence rate.
Moreover, figure~\ref{fig:case1-2-convergence}(b) shows the convergence of the error
of the gradient of the MB solution in the $L_2\/$ and $L_{\infty}\/$ norms. Here,
the gradient was computed using standard \snd\/ order accurate centered differences, as in
\eqref{eq:partialx} and \eqref{eq:partialy}. As we can observe, the gradient
converges to \fst\/ order in the $L_{\infty}\/$ norm and, apparently, as $h^{3/2}$ in the
$L_2\/$ norm.
%
\clearpage
%

\subsubsection{Example 2$_s$.}\label{subsub:case2}
\begin{figure}[htb!]
 \begin{center}
  \includegraphics[width=3.5in]{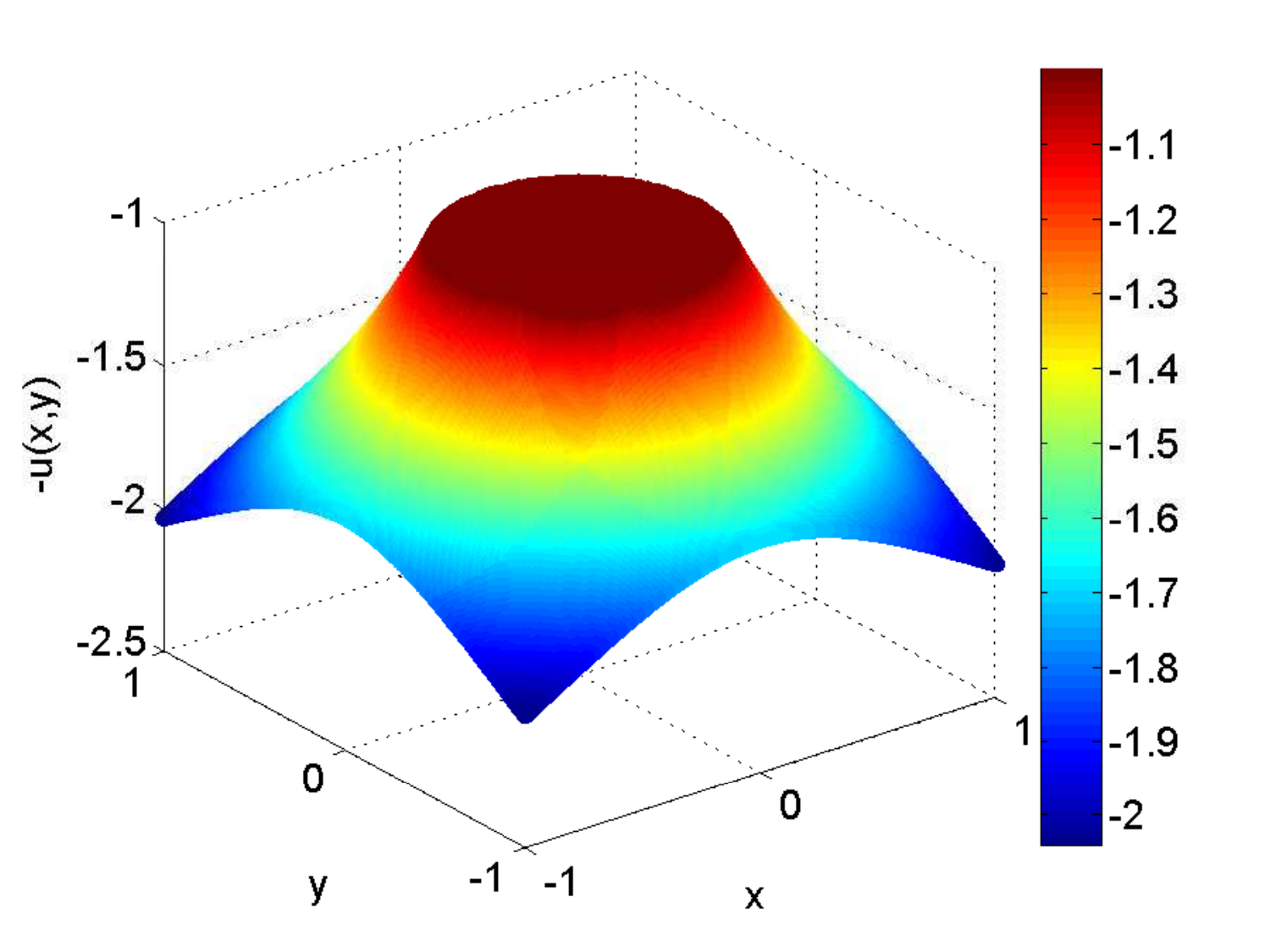}
 \end{center}
 \caption{Example 2$_s$. Numerical solution with $161 \times 161\/$ nodes.}
 \label{fig:case2-2-numerical}
\end{figure}

\begin{figure}[htb!]
 \begin{center}
  \includegraphics[width=2.6in]{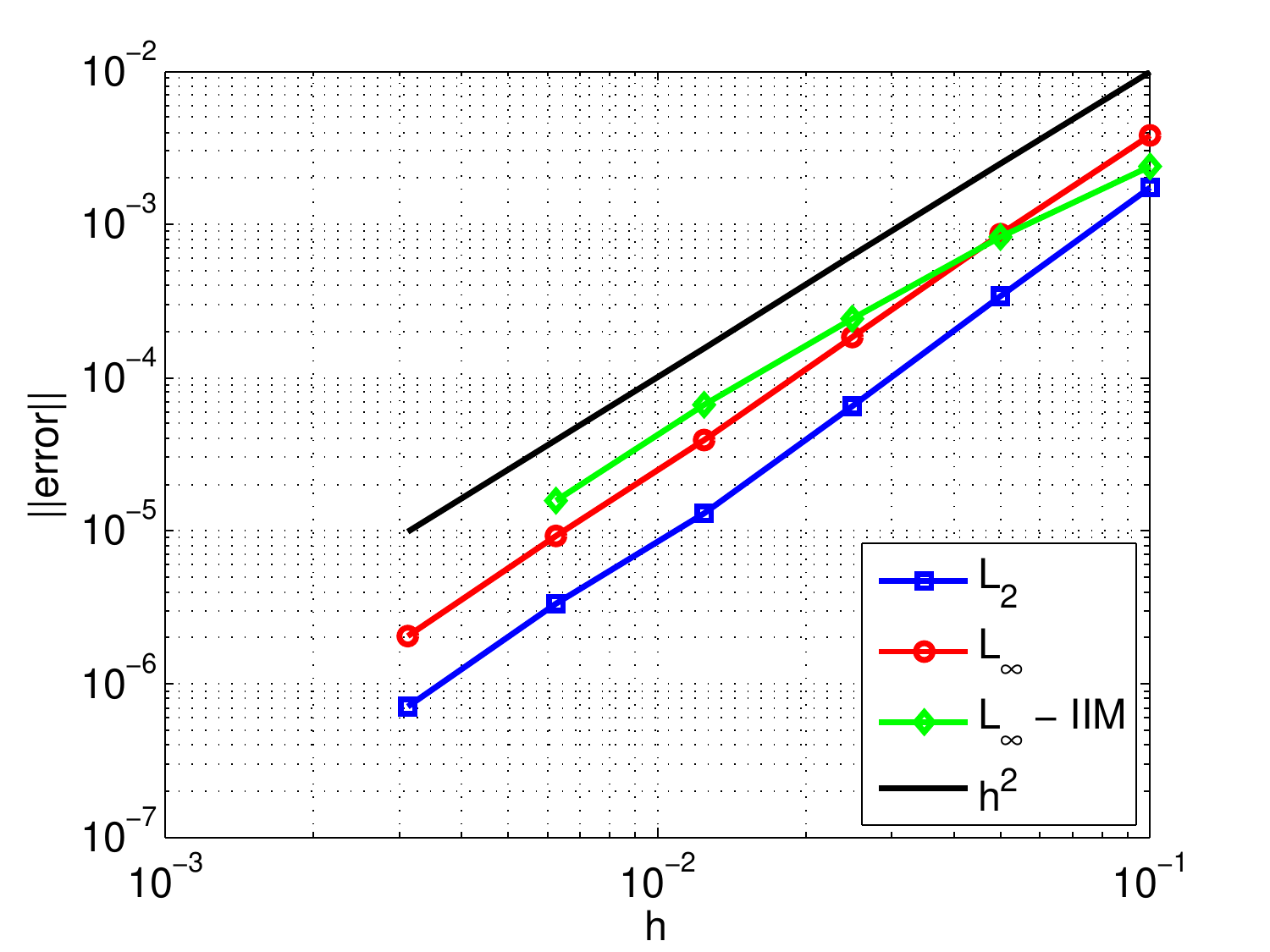}
  \includegraphics[width=2.6in]{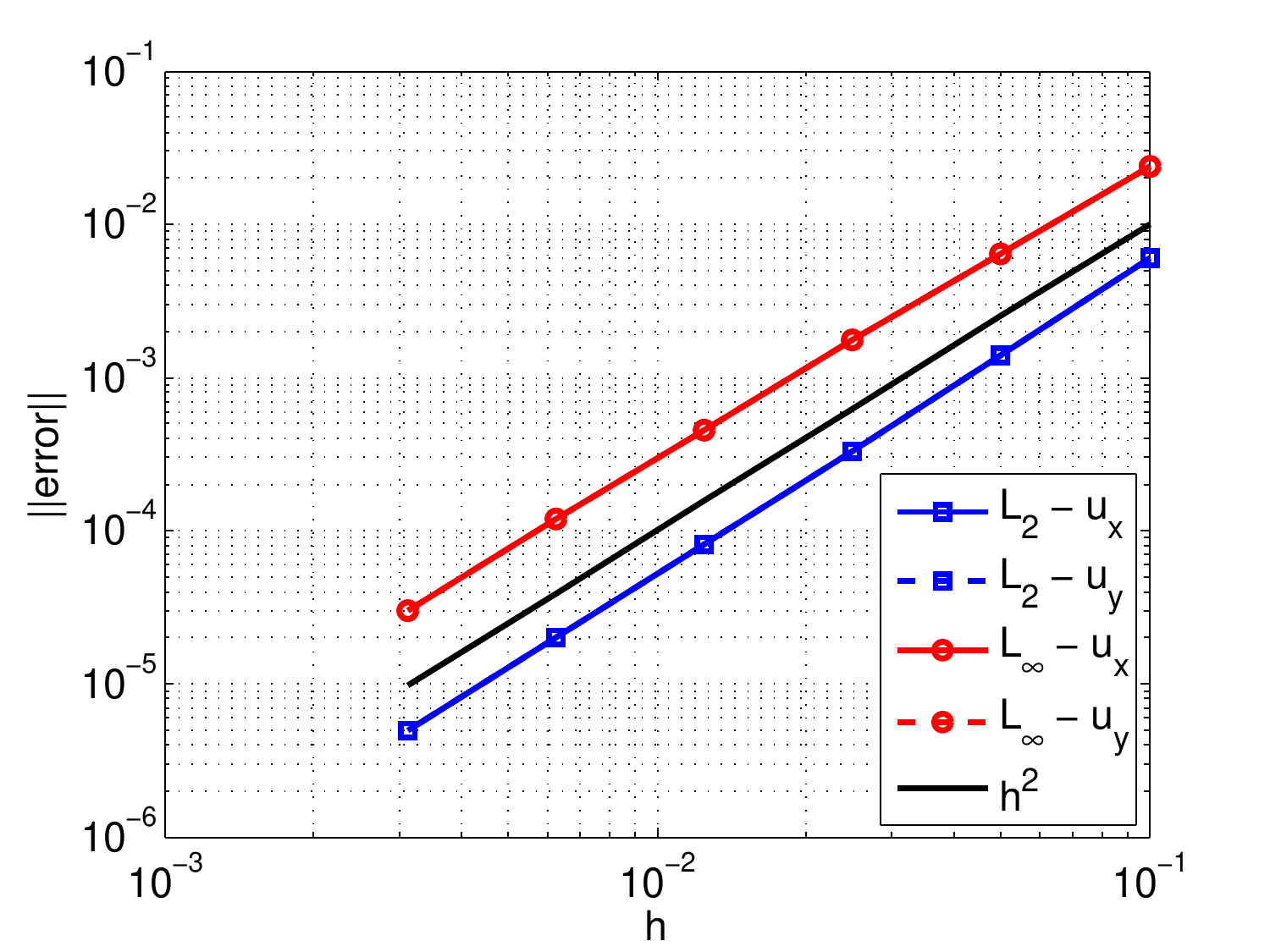}\\
  (a) Solution. \hspace*{1.6in} (b) Gradient.\hfill
 \end{center}
 \caption{Example 2$_s$ - Error behavior of the solution and its gradient in the $L_2\/$
          and $L_\infty\/$ norms. IIM refers to results obtained with the Immersed Interface
          Method in~\cite{leveque:94} (Copyright \copyright 1994 Society for Industrial and
          Applied Mathematics. Data compiled with permission. All rights reserved).}
 \label{fig:case2-2-convergence}
\end{figure}

\begin{itemize}
 \item 
 Domain: $(x,y) \in [-1,1] \times [-1,1]$.
 \item 
 Problem parameters:
 \begin{align*}
   f^+\p{x\/,\,y}           & = 0\/,\\
   f^-\p{x\/,\,y}           & = 0\/,\\
   [u]_{\Gamma}              & = \log\p{2\sqrt{x^2+y^2}}\/,\\
   \left[u_n\right]_{\Gamma} & = \dfrac{x\,n_x + y\,n_y}{x^2+y^2}\/.
 \end{align*}
 \item 
 Level set defining the interface: \hfill
 $\phi\p{x\/,\,y} = \sqrt{(x-x_0)^2 + (y-y_0)^2} - r_0\/$,\\
 where $x_0 = 0\/$, $y_0 = 0\/$, and $r_0 = 0.5\/$.
 \item 
 Exact solution:
 \begin{align*}
   u^+\p{x\/,\,y} & = 1 + \log\p{2\sqrt{x^2+y^2}}\/,\\
   u^-\p{x\/,\,y} & = 1.
 \end{align*}
\end{itemize}
\newpage
As mentioned earlier, this example corresponds to example 1
in~\cite{leveque:94}, where the same problem is solved using the Immersed
Interface Method. Hence, this provides a good opportunity to compare the CFM
with the well established IIM. Figure~\ref{fig:case2-2-numerical} shows the
numerical solution with a fine grid ($161 \times 161$ nodes). Once again, the
overall quality of the solution is very satisfactory.
Figure~\ref{fig:case2-2-convergence} shows the behavior of the error in the
$L_2$ and $L_{\infty}$ norms and the solution converges to \snd\/ order.
However, in this example the gradient also appears to converge to \snd\/ order.

Figure~\ref{fig:case2-2-convergence}(a) also includes the convergence of the
solution error in the $L_{\infty}$ norm obtained with the IIM --- we plot the errors
listed in table~1 of \cite{leveque:94}. Both methods produce similar convergence
rates. However, in this example at least, the CFM produces slightly smaller
errors --- by a factor of about 1.7.
%
%
\clearpage
%
%
\subsubsection{Example 3$_s$.}\label{subsub:case3}
\begin{figure}[htb!]
 \begin{center}
  \includegraphics[width=3.5in]{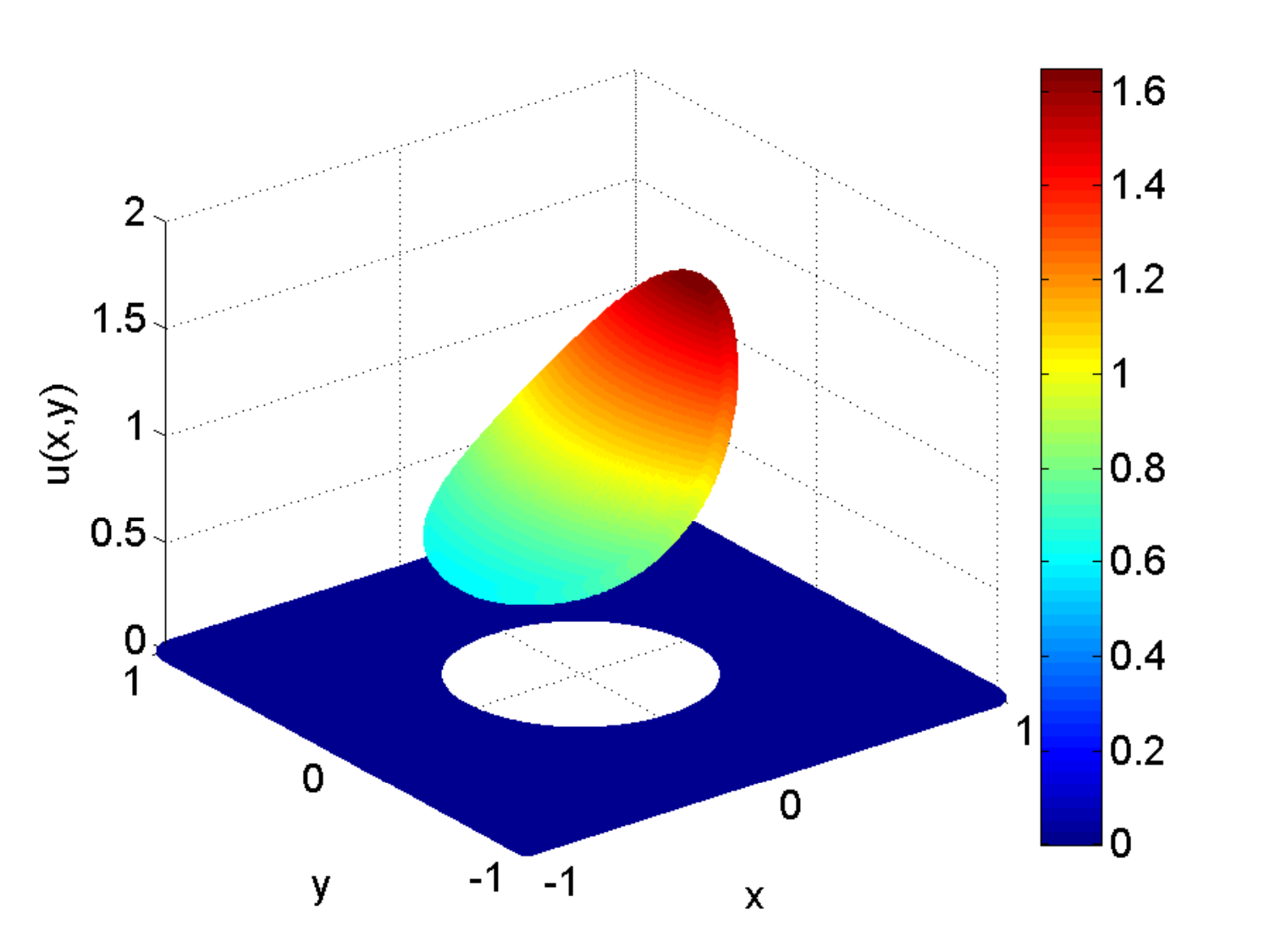}
 \end{center}
 \caption{Example 3$_s$. Numerical solution with $161 \times 161\/$ nodes.}
 \label{fig:case3-2-numerical}
\end{figure}

\begin{figure}[htb!]
 \begin{center}
  \includegraphics[width=2.6in]{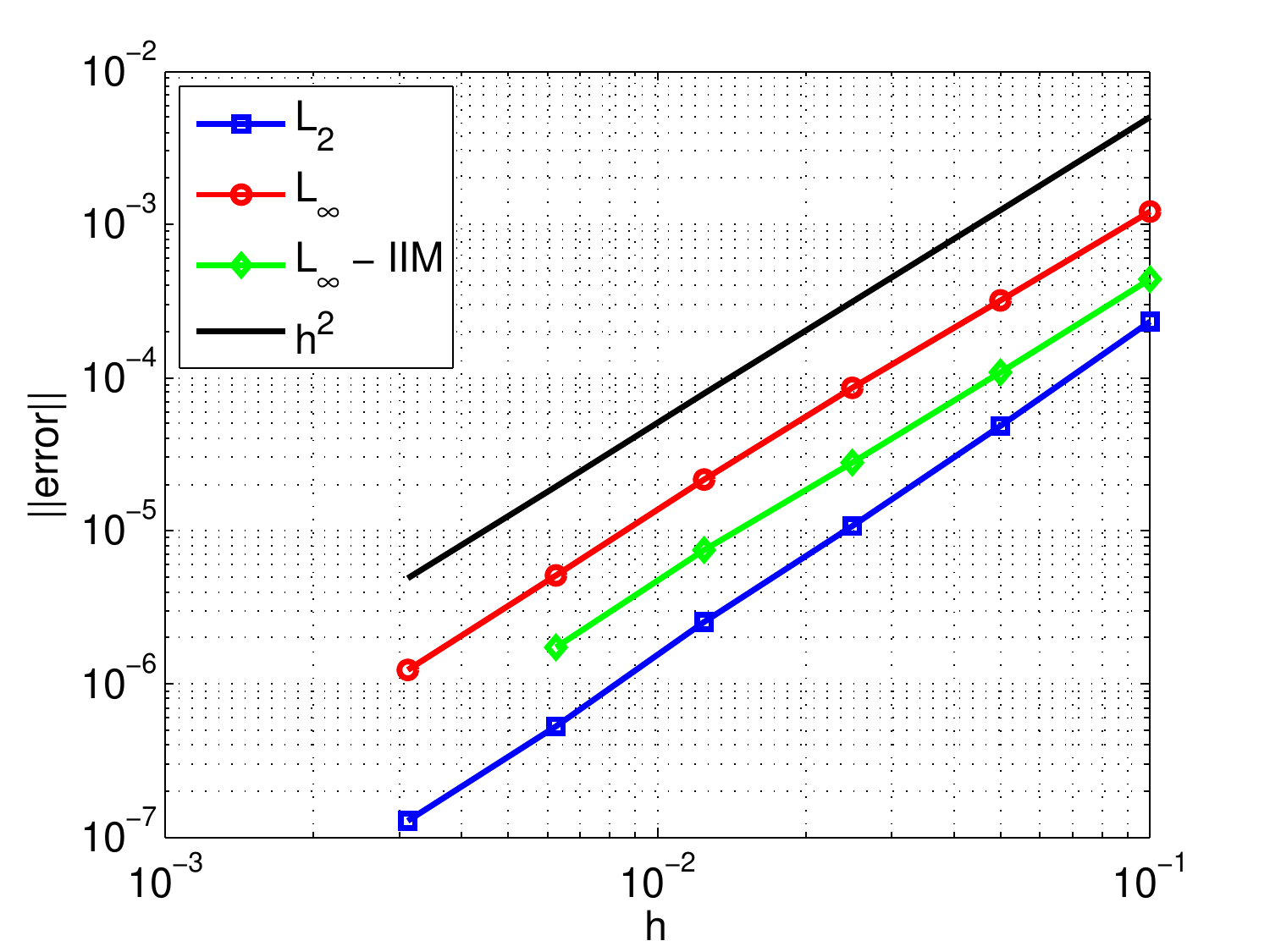}
  \includegraphics[width=2.6in]{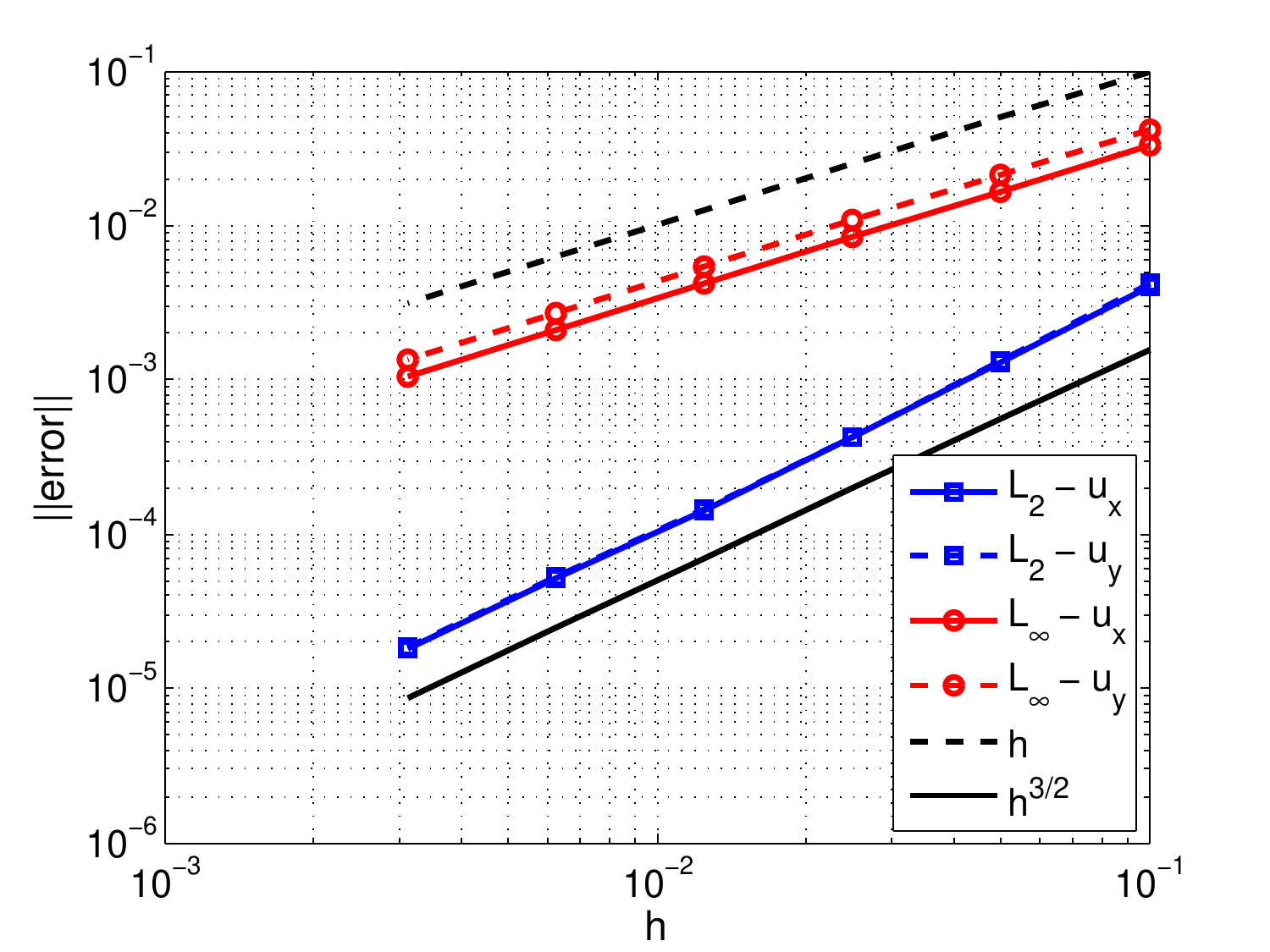}\\
  (a) Solution. \hspace*{1.6in} (b) Gradient.\hfill
 \end{center}
 \caption{Example 3$_s$ - Error behavior of the solution and its gradient in the $L_2\/$
          and $L_\infty\/$ norms. IIM refers to results obtained with the Immersed Interface
          Method in~\cite{leveque:94} (Copyright \copyright 1994 Society for Industrial and
          Applied Mathematics. Data compiled with permission. All rights reserved).}
 \label{fig:case3-2-convergence}
\end{figure}

\begin{itemize}
 \item 
 Domain: $(x,y) \in [-1,1] \times [-1,1]$.
 \item 
 Problem parameters:
 \begin{align*}
   f^+\p{x\/,\,y}           & = 0,\\
   f^-\p{x\/,\,y}           & = 0,\\
   [u]_{\Gamma}              & = -\exp(x)\cos(y)\/,\\
   \left[u_n\right]_{\Gamma} & = \exp(x)\left[-\cos(y)\,n_x +
                                \sin(y)\,n_y\right]\/.
 \end{align*}
 \item 
 Level set defining the interface: \hfill
 $\phi\p{x\/,\,y} = \sqrt{(x-x_0)^2 + (y-y_0)^2} - r_0\/$, \\
 where $x_0 = 0\/$, $y_0 = 0\/$, and $r_0 = 0.5\/$.
 \item 
 Exact solution:
 \begin{align*}
   u^+\p{x\/,\,y} & = 0\/, \\ u^-\p{x\/,\,y} & = \exp(x)\cos(y)\/.
 \end{align*}
\end{itemize}
This example corresponds to example 3 in [3] for the IIM method.
Figure~\ref{fig:case3-2-numerical} shows the numerical solution with a fine
grid ($161 \times 161\/$ nodes) while figure~\ref{fig:case3-2-convergence}
presents convergence results for the errors in the $L_2\/$ and $L_\infty\/$
norms. In addition, figure~\ref{fig:case3-2-convergence}(a) also includes the
$L_{\infty}$ norm of the error obtained with the IIM --- we plot the errors
listed in table~3 of \cite{leveque:94}\footnotemark[\value{footnote}].
In this case the IIM produces slightly smaller errors than the CFM --- by a
factor of about 2.8. Therefore, based on the results from examples 2$_s$ and
3$_s$, we may conclude that the IIM and the CFM produce results of comparable
accuracy --- each method generating slightly smaller errors in different cases.
%
%

%
%

\section*{Acknowledgements}

The authors would like to acknowledge the National Science Foundation support --- this research was partially supported by grant DMS-0813648. In addition, the first author acknowledges the support by Coordena\c{c}\~ao de Aperfei\c{c}oamento de Pessoal de N\'ivel Superior (CAPES -- Brazil) and the Fulbright Commission through grant BEX 2784/06-8. The second author also acknowledges support from the NSERC Discovery program. Finally, the authors would also like to acknowledge many helpful conversations with Prof.~B.~Seibold at Temple University, and the many helpful suggestions made by the referees of this paper.



\section*{References}

\bibliographystyle{elsarticle-num}
\bibliography{biblio}

\begin{thebibliography}{10}
\expandafter\ifx\csname url\endcsname\relax
  \def\url#1{\texttt{#1}}\fi
\expandafter\ifx\csname urlprefix\endcsname\relax\def\urlprefix{URL }\fi
\expandafter\ifx\csname href\endcsname\relax
  \def\href#1#2{#2} \def\path#1{#1}\fi

\bibitem{peskin:77}
C.~S. Peskin, Numerical analysis of blood flow in the heart, Journal of
  Computational Physics 25~(3) (1977) 220--252.
\newblock \href {http://dx.doi.org/10.1016/0021-9991(77)90100-0}
  {\path{doi:10.1016/0021-9991(77)90100-0}}.

\bibitem{sussman:94}
M.~Sussman, P.~Smereka, S.~Osher, A level set approach for computing solutions
  to incompressible two-phase flow, Journal of Computational Physics 114~(1)
  (1994) 146--159.
\newblock \href {http://dx.doi.org/10.1006/jcph.1994.1155}
  {\path{doi:10.1006/jcph.1994.1155}}.

\bibitem{leveque:94}
R.~J. LeVeque, Z.~Li, The immersed interface method for elliptic equations with
  discontinuous coefficients and singular sources, SIAM Journal on Numerical
  Analysis 31~(4) (1994) 1019--1044.
\newblock \href {http://dx.doi.org/10.1137/0731054}
  {\path{doi:10.1137/0731054}}.

\bibitem{leveque:97}
R.~J. LeVeque, Z.~Li, Immersed interface methods for {Stokes} flow with elastic
  boundaries or surface tension, SIAM Journal on Scientific Computing 18~(3)
  (1997) 709--735.
\newblock \href {http://dx.doi.org/10.1137/S1064827595282532}
  {\path{doi:10.1137/S1064827595282532}}.

\bibitem{johansen:98}
H.~Johansen, P.~Colella, A {Cartesian} grid embedded boundary method for
  {Poisson's} equation on irregular domains, Journal of Computational Physics
  147~(1) (1998) 60--85.
\newblock \href {http://dx.doi.org/10.1006/jcph.1998.5965}
  {\path{doi:10.1006/jcph.1998.5965}}.

\bibitem{fedkiw:99}
R.~P. Fedkiw, T.~Aslam, S.~Xu, The ghost fluid method for deflagration and
  detonation discontinuities, Journal of Computational Physics 154~(2) (1999)
  393--427.
\newblock \href {http://dx.doi.org/10.1006/jcph.1999.6320}
  {\path{doi:10.1006/jcph.1999.6320}}.

\bibitem{fedkiwetal:99}
R.~P. Fedkiw, T.~Aslam, B.~Merriman, S.~Osher, A non-oscillatory {Eulerian}
  approach to interfaces in multimaterial flows (the ghost fluid method),
  Journal of Computational Physics 152~(2) (1999) 457--492.
\newblock \href {http://dx.doi.org/10.1006/jcph.1999.6236}
  {\path{doi:10.1006/jcph.1999.6236}}.

\bibitem{kang:00}
M.~Kang, R.~P. Fedkiw, X.-D. Liu, A boundary condition capturing method for
  multiphase incompressible flow, Journal of Scientific Computing 15 (2000)
  323--360.
\newblock \href {http://dx.doi.org/10.1023/A:1011178417620}
  {\path{doi:10.1023/A:1011178417620}}.

\bibitem{liu:00}
X.-D. Liu, R.~P. Fedkiw, M.~Kang, A boundary condition capturing method for
  {Poisson's} equation on irregular domains, Journal of Computational Physics
  160~(1) (2000) 151--178.
\newblock \href {http://dx.doi.org/10.1006/jcph.2000.6444}
  {\path{doi:10.1006/jcph.2000.6444}}.

\bibitem{lai:00}
M.-C. Lai, C.~S. Peskin, An immersed boundary method with formal second-order
  accuracy and reduced numerical viscosity, Journal of Computational Physics
  160~(2) (2000) 705--719.
\newblock \href {http://dx.doi.org/10.1006/jcph.2000.6483}
  {\path{doi:10.1006/jcph.2000.6483}}.

\bibitem{li:01}
Z.~Li, M.-C. Lai, The immersed interface method for the {Navier-Stokes}
  equations with singular forces, Journal of Computational Physics 171~(2)
  (2001) 822--842.
\newblock \href {http://dx.doi.org/10.1006/jcph.2001.6813}
  {\path{doi:10.1006/jcph.2001.6813}}.

\bibitem{nguyen:01}
D.~Q. Nguyen, R.~P. Fedkiw, M.~Kang, A boundary condition capturing method for
  incompressible flame discontinuities, Journal of Computational Physics
  172~(1) (2001) 71--98.
\newblock \href {http://dx.doi.org/10.1006/jcph.2001.6812}
  {\path{doi:10.1006/jcph.2001.6812}}.

\bibitem{lee:03}
L.~Lee, R.~J. LeVeque, An immersed interface method for incompressible
  {Navier-Stokes} equations, SIAM Journal on Scientific Computing 25~(3) (2003)
  832--856.
\newblock \href {http://dx.doi.org/10.1137/S1064827502414060}
  {\path{doi:10.1137/S1064827502414060}}.

\bibitem{gibou:07}
F.~Gibou, L.~Chen, D.~Nguyen, S.~Banerjee, A level set based sharp interface
  method for the multiphase incompressible {Navier-Stokes} equations with phase
  change, Journal of Computational Physics 222~(2) (2007) 536--555.
\newblock \href {http://dx.doi.org/10.1016/j.jcp.2006.07.035}
  {\path{doi:10.1016/j.jcp.2006.07.035}}.

\bibitem{gong:08}
Y.~Gong, B.~Li, Z.~Li, Immersed--interface finite-element methods for elliptic
  interface problems with nonhomogeneous jump conditions, SIAM Journal on
  Numerical Analysis 46~(1) (2008) 472--495.
\newblock \href {http://dx.doi.org/10.1137/060666482}
  {\path{doi:10.1137/060666482}}.

\bibitem{dolbow:09}
J.~Dolbow, I.~Harari, An efficient finite element method for embedded interface
  problems, International Journal for Numerical Methods in Engineering 78~(2)
  (2009) 229--252.
\newblock \href {http://dx.doi.org/10.1002/nme.2486}
  {\path{doi:10.1002/nme.2486}}.

\bibitem{bedrossian:10}
J.~Bedrossian, J.~H. von Brecht, S.~Zhu, E.~Sifakis, J.~M. Teran, A second
  order virtual node method for elliptic problems with interfaces and irregular
  domains, Journal of Computational Physics 229~(18) (2010) 6405 -- 6426.
\newblock \href {http://dx.doi.org/DOI: 10.1016/j.jcp.2010.05.002}
  {\path{doi:DOI: 10.1016/j.jcp.2010.05.002}}.

\bibitem{trefthen:97}
L.~N. Trefethen, D.~Bau, Numerical Linear Algebra, SIAM: Society for Industrial
  and Applied Mathematics, 1997.

\bibitem{mayo:84}
A.~Mayo, The fast solution of {Poisson's} and the biharmonic equations on
  irregular regions, SIAM Journal on Numerical Analysis 21~(2) (1984) 285--299.
\newblock \href {http://dx.doi.org/10.1137/0721021}
  {\path{doi:10.1137/0721021}}.

\bibitem{udaykumar:99}
H.~S. Udaykumar, R.~Mittal, W.~Shyy, Computation of solid-liquid phase fronts
  in the sharp interface limit on fixed grids, Journal of Computational Physics
  153~(2) (1999) 535--574.
\newblock \href {http://dx.doi.org/10.1006/jcph.1999.6294}
  {\path{doi:10.1006/jcph.1999.6294}}.

\bibitem{gibou:02}
F.~Gibou, R.~P. Fedkiw, L.-T. Cheng, M.~Kang, A second-order-accurate symmetric
  discretization of the {Poisson} equation on irregular domains, Journal of
  Computational Physics 176~(1) (2002) 205--227.
\newblock \href {http://dx.doi.org/10.1006/jcph.2001.6977}
  {\path{doi:10.1006/jcph.2001.6977}}.

\bibitem{jomaa:05}
Z.~Jomaa, C.~Macaskill, The embedded finite difference method for the {Poisson}
  equation in a domain with an irregular boundary and {Dirichlet} boundary
  conditions, Journal of Computational Physics 202~(2) (2005) 488--506.
\newblock \href {http://dx.doi.org/10.1016/j.jcp.2004.07.011}
  {\path{doi:10.1016/j.jcp.2004.07.011}}.

\bibitem{gibou:05}
F.~Gibou, R.~Fedkiw, A fourth order accurate discretization for the {Laplace}
  and heat equations on arbitrary domains, with applications to the {Stefan}
  problem, Journal of Computational Physics 202~(2) (2005) 577--601.
\newblock \href {http://dx.doi.org/10.1016/j.jcp.2004.07.018}
  {\path{doi:10.1016/j.jcp.2004.07.018}}.

\bibitem{chen:07}
H.~Chen, C.~Min, F.~Gibou, A supra-convergent finite difference scheme for the
  {Poisson} and heat equations on irregular domains and non-graded adaptive
  {Cartesian} grids, Journal of Scientific Computing 31~(1) (2007) 19--60.
\newblock \href {http://dx.doi.org/10.1007/s10915-006-9122-8}
  {\path{doi:10.1007/s10915-006-9122-8}}.

\bibitem{sethian:92}
J.~A. Sethian, J.~Strain, Crystal growth and dendritic solidification, Journal
  of Computational Physics 98~(2) (1992) 231--253.
\newblock \href {http://dx.doi.org/10.1016/0021-9991(92)90140-T}
  {\path{doi:10.1016/0021-9991(92)90140-T}}.

\bibitem{chen:97}
S.~Chen, B.~Merriman, S.~Osher, P.~Smereka, A simple level set method for
  solving {Stefan} problems, Journal of Computational Physics 135~(1) (1997)
  8--29.
\newblock \href {http://dx.doi.org/10.1006/jcph.1997.5721}
  {\path{doi:10.1006/jcph.1997.5721}}.

\bibitem{osher:88}
S.~Osher, J.~A. Sethian, Fronts propagating with curvature--dependent speed:
  Algorithms based on hamilton-jacobi formulations, Journal of Computational
  Physics 79~(1) (1988) 12 -- 49.
\newblock \href {http://dx.doi.org/DOI: 10.1016/0021-9991(88)90002-2}
  {\path{doi:DOI: 10.1016/0021-9991(88)90002-2}}.

\bibitem{nave:10}
J.-C. Nave, R.~R. Rosales, B.~Seibold, A gradient-augmented level set method
  with an optimally local, coherent advection scheme, Journal of Computational
  Physics 229~(10) (2010) 3802 -- 3827.
\newblock \href {http://dx.doi.org/10.1016/j.jcp.2010.01.029}
  {\path{doi:10.1016/j.jcp.2010.01.029}}.

\end{thebibliography}







\end{document}